\documentclass[12pt,reqno,a4paper]{amsart}

\usepackage[top=1in, bottom=1in, left=1in, right=1in]{geometry}
\usepackage[OT2,T1]{fontenc}
\usepackage[english]{babel}
\usepackage{newtxtext}
\usepackage{newtxmath}
\usepackage{bbm}
\usepackage{tikz-cd}
\usepackage{enumitem}

\theoremstyle{plain}
\newtheorem{introtheorem}{Theorem}

\newtheorem{theorem}{Theorem}[section]
\newtheorem{lemma}[theorem]{Lemma}

\newtheorem{proposition}[theorem]{Proposition}
\theoremstyle{definition}
\newtheorem{definition}[theorem]{Definition}
\theoremstyle{remark}
\newtheorem{remark}[theorem]{Remark}

\newtheorem*{remark*}{Remark}
\newtheorem*{remarks*}{Remarks}

\numberwithin{equation}{section}

\usepackage[unicode=true,pdfusetitle, bookmarks=true,bookmarksnumbered=false,bookmarksopen=false, breaklinks=true,pdfborder={0 0 0},pdfborderstyle={},backref=false,colorlinks=true]{hyperref}
\definecolor{myblue}{rgb}{0.09,0.52,0.74}
\definecolor{mygreen}{rgb}{0.05,0.6,0.2}
\hypersetup{pdfborder={0 0 0},pdfborderstyle={},colorlinks=true,linkcolor=myblue,citecolor=mygreen,urlcolor=blue}
\usepackage{amsmath}
\usepackage{mathtools}
\usepackage{mathrsfs}
\usepackage{comment}

\newcommand{\Z}{\mathbf{Z}}
\newcommand{\R}{\mathbf{R}}
\newcommand{\Q}{\mathbf{Q}}
\newcommand{\F}{\mathbf{F}}
\renewcommand{\C}{\mathbf{C}}

\newcommand{\SR}{\mathsf{R}}
\newcommand{\ms}{\mathsf}
\newcommand{\mb}{\mathbbm}

\newcommand{\mr}{\mathrm}
\renewcommand{\geq}{\geqslant}
\renewcommand{\ge}{\geqslant}
\renewcommand{\leq}{\leqslant}
\renewcommand{\le}{\leqslant}
\renewcommand{\epsilon}{\varepsilon}
\DeclareMathOperator{\Res}{Res}
\newcommand{\car}{\operatorname{char}}
\DeclareMathOperator{\res}{res}
\DeclareMathOperator{\sign}{sign}
\DeclareMathOperator{\proj}{proj}

\DeclareMathOperator{\supp}{supp}
\DeclareMathOperator{\id}{id}
\newcommand{\ii}{\mathrm{i}}
\renewcommand{\bar}[1]{\overline{#1}}

\makeatletter
\newcommand*{\brmod}{%
  \nonscript\mskip-\medmuskip\mkern5mu%
  \mathbin{\operator@font rmod}\penalty900\mkern5mu%
  \nonscript\mskip-\medmuskip
}
\makeatother

\DeclarePairedDelimiter\abs{\lvert}{\rvert}%

\makeatletter
\let\oldabs\abs
\def\abs{\@ifstar{\oldabs}{\oldabs*}}
\makeatother
\newcommand{\norm}[1]{\| #1\| }

\newcommand{\CA}{\mathcal{A}}
\newcommand{\CB}{\mathcal{B}}
\newcommand{\CC}{\mathcal{C}}
\newcommand{\CD}{\mathcal{D}}
\newcommand{\CE}{\mathcal{E}}

\newcommand{\CG}{\mathcal{G}}
\newcommand{\CH}{\mathcal{H}}
\newcommand{\CI}{\mathcal{I}}
\newcommand{\CL}{\mathcal{L}}
\newcommand{\CM}{\mathcal{M}}
\newcommand{\CN}{\mathcal{N}}
\newcommand{\CO}{\mathcal{O}}
\newcommand{\CP}{\mathcal{P}}

\newcommand{\CR}{\mathcal{R}}
\newcommand{\CS}{\mathcal{S}}
\newcommand{\CT}{\mathcal{T}}
\newcommand{\CU}{\mathcal{U}}
\newcommand{\CY}{\mathcal{Y}}

\newcommand{\BA}{\boldsymbol{A}}
\newcommand{\BB}{\boldsymbol{B}}
\newcommand{\BC}{\boldsymbol{C}}
\newcommand{\BD}{\boldsymbol{D}}
\newcommand{\BG}{\boldsymbol{G}}
\newcommand{\BH}{\boldsymbol{H}}
\newcommand{\BU}{\boldsymbol{U}}
\newcommand{\BK}{\boldsymbol{K}}

\newcommand{\Bl}{\boldsymbol{\ell}}
\newcommand{\Bk}{\boldsymbol{k}}
\newcommand{\Br}{\boldsymbol{r}}
\newcommand{\Bzero}{\boldsymbol{0}}

\newcommand{\mfp}{\mathfrak{p}}

\definecolor{amethyst}{rgb}{0.6, 0.4, 0.8}
\definecolor{darkcyan}{rgb}{0.0, 0.55, 0.55}

\makeatletter
\@namedef{subjclassname@1991}{\textup{2020} Mathematics Subject Classification}
\makeatother

\begin{document}

\date{\today}
\title[Random reciprocal polynomials of large degree]{Irreducibility and Galois groups of random reciprocal polynomials of large degree}

\author[D.~Hokken]{David Hokken}
\address{\normalfont DH: Mathematisch Instituut\\
        Universiteit Utrecht\\
        Postbus 80.010, 3508 TA Utrecht, Nederland}
\email{{\tt d.p.t.hokken@uu.nl}}

\author[D.~Koukoulopoulos]{Dimitris Koukoulopoulos}
\address{\normalfont DK: D\'epartement de math\'ematiques et de statistique\\
Universit\'e de Montr\'eal\\
CP 6128 succ.~Centre-Ville\\
Montr\'eal, QC H3C 3J7\\
Canada}
\email{{\tt dimitris.koukoulopoulos@umontreal.ca}}


\begin{abstract} \noindent
Let $A = a_0T^m + \sum_{j=1}^{m-1} a_j (T^{m-j}+T^{m+j}) + T^{2m}+1  \in \Z[T]$ be a monic reciprocal polynomial of degree $2m$ sampled randomly by selecting its coefficients $a_0,a_1,\dots,a_{m-1}$ independently according to a given probability measure $\mu$ on $\Z$. For a wide range of measures $\mu$, we prove that $A$ is irreducible with probability $\ge 1-Cm^{-c}$ for some constants $c,C>0$. In addition, we prove that with the same probability the Galois group $\CG_A$ of $A$ is either the full hyperoctahedral group $\CC_2 \wr \CS_m$ or one of two of its index-$2$ subgroups. 
The main condition that $\mu$ must satisfy is of Fourier-theoretic nature, and holds for example when $\mu$ is the uniform measure on a set of at least $35$ consecutive integers, or on an arbitrary, sufficiently large subset of an interval $[-H,H]$, with $H$ larger than some absolute constant. Our most general result allows for each $a_j$ to be sampled by its own probability measure $\mu_j$.

Our approach builds on earlier work of Bary-Soroker, Kozma and the second author, who proved for essentially the same $\mu_j$ that the `standard' monic polynomial $a_0 + \cdots + a_{m-1}T^{m-1} + T^m$, conditioning on $a_0 \neq 0$, is irreducible and has as Galois group either the symmetric group $\CS_m$ or the alternating group $\CA_m$ with high probability. For reciprocal polynomials, we can study the discriminant of $A$ and rule out with high probability that $\CG_A \leq (\CC_2 \wr \CS_m) \cap \CA_{2m}$, i.e., that $\CG_A$ is contained in the maximal alternating subgroup. Furthermore, we establish a bound on the probability that $\CG_A \leq \CC_2 \times \CS_m$ by examining the Frobenius at suitable primes. In the process, we establish a \L uczak--Pyber theorem for the group $\CC_2 \wr \CS_m$, which may be of independent interest.
\end{abstract}

\maketitle
\setcounter{tocdepth}{1}
\tableofcontents
\newpage

\section{Introduction}
\label{sec:introduction}

Factoring polynomials, as well as understanding symmetries between their zeros, is a classical subject in number theory. Both of these aspects are captured by the Galois group of the polynomial. Here, we focus on polynomials with integer coefficients, say $A(T)\in\Z[T]$, and we will denote by $\CG_A$ the Galois group of (the splitting field of) $A$ over $\Q$. If $A$ is of degree $n$, we may view $\CG_A$ as a subgroup of $\CS_n$, the symmetric group on $n$ letters. Regarding the size of the Galois group, the general philosophy is that $\CG_A$ is ``almost always as large as it can be for a typical choice of $A$''. In other words, if we select $A$ randomly among all polynomials of degree $n$, then we expect that $\CG_A=\CS_n$ almost surely. Of course, we must specify what \emph{random} means. A classical case is the following: we fix a finite set $\CN$ of integers and we select all coefficients of $A$ independently and uniformly at random from $\CN$. Following Van der Waerden's pioneering work \cite{VdW1, VdW2}, much of the literature (see, e.g., \cite{ABD, BBM, Bhargava, CD, DDS, Gallagher}) focuses on the \emph{large box model}, where the degree $n$ is fixed and we take a growing set of coefficients, usually $\CN=[-H,H]\cap\Z$ with $H\to\infty$. The state-of-the-art in this direction is Bhargava's result \cite{Bhargava} that $\mb{P}(\CG_A \neq \CS_n) = \mb{P}(A \textup{ is reducible}) + \mb{P}(\CG_A = \CA_n) + o(H^{-1})$ and that $\mb{P}(\CG_A = \CA_n) = O(H^{-1})$, where $\CA_n$ is the alternating group on $n$ letters.

Instead, here we will focus on the \emph{restricted coefficient model} where $\CN$ is fixed and the degree $n$ tends to infinity. This was first considered by Odlyzko and Poonen \cite{OP}. Taking $\CN = \{0,1\}$ in the uniform measure, they stated as a folklore conjecture that the probability that $A$ is irreducible tends to $1$ as $n \to \infty$, conditioning on $A(0) = 1$. Konyagin \cite{Konyagin} showed that the probability is $\gg 1/\log{n}$. Breuillard and Varj\'u \cite{BV} proved the folklore conjecture under the Riemann hypothesis for Dedekind zeta functions. In 2020, Bary-Soroker, Kozma and the second author \cite{BKK} improved Konyagin's unconditional lower bound to $\gg 1$, a result that continues to hold when we instead sample the coefficients of $A$ independently according to a fixed probability measure $\mu$ on $\Z$ with rather generic Fourier-analytical properties, conditioning on $A(0) \neq 0$. Furthermore, again under a wide range of choices of $\mu$, they proved that $\CG_A$ contains $\CA_n$ with high probability. When $\mu$ is the uniform measure on $\{\pm 1\}$, Bary-Soroker, Kozma, Poonen and the first author \cite{BHKP} showed that $\limsup_{n \to \infty} \mb{P}(A \textup{ is irreducible}) = 1$.

In the present paper, we study the irreducibility and Galois group of polynomials $A$ in the restricted coefficient model with an additional constraint: we shall assume that $A$ is \emph{reciprocal}.

\begin{definition}[Reciprocal polynomial]
\label{def:reciprocal}
Let $A$ be a polynomial over a field $K$. We say that $A$ is \emph{reciprocal} if either $A=0$, or $A$ has even degree and satisfies the relation
\begin{equation}
\label{eq:reciprocal_symmetry}
A(T) = T^{\deg{A}} A(1/T).
\end{equation}
Equivalently, $A$ is reciprocal if either $A$ is a constant or there are $m\in\Z_{\ge1}$ and $a_0,a_1,\dots,a_m\in K$ with $a_m \neq 0$ such that
\begin{equation}
\label{eq:reciprocal polynomial}
		A(T) = a_m + a_{m-1}T + \cdots + a_{1} T^{m-1} + a_0 T^m + a_{1} T^{m+1} + \cdots + a_{m-1}T^{2m-1} + a_m T^{2m}.
\end{equation} 
\end{definition}

\begin{remark}
\label{rem:odd degree reciprocal}
If $A$ has odd degree and satisfies \eqref{eq:reciprocal_symmetry}, then we may write $A(T)=(T+1)B(T)$ with $B$ a reciprocal polynomial.
\end{remark}

For $m \ge 2$, the Galois group of squarefree $A$ with $A(0) \neq 0$ as in \eqref{eq:reciprocal polynomial} is always smaller than $\CS_{2m}$, because of the relations between the zeros of $A$ implied by \eqref{eq:reciprocal_symmetry}. In fact, $\CG_A$ is contained in the \emph{hyperoctahedral group} $\CC_2 \wr \CS_m$, whose construction we recall in \S \ref{sec:hyperoctahedral_group}. The following is one of our results about random reciprocal polynomials.

\begin{introtheorem}
\label{thm:simple_interval}
Let $\CN$ be a set of at least $35$ consecutive integers and let $m$ be a positive integer. Let $a_m=1$ and sample integers $a_0, a_1, \ldots, a_{m-1}$ independently and uniformly at random from $\CN$ and let these be the coefficients of the monic reciprocal polynomial $A$ as in \eqref{eq:reciprocal polynomial}. Then
\begin{equation*}
\lim_{m \to \infty} \mb{P}\big(A \textup{ is irreducible in } \Z[T]\big) = 1.
\end{equation*}
Moreover, if $\CG_A$ denotes the Galois group of $A$ over $\Q$, then
\begin{equation*}
\lim_{m \to \infty} \mb{P}\Big(\CG_A  \in \big\{ \CC_2 \wr \CS_m, \CC_2\wr \CA_m, G_2\big\}\Big)  = 1,
\end{equation*}
where $G_2$ denotes the group defined in \eqref{eq:G2}.
\end{introtheorem}

This introductory section is dedicated to presenting our main results (see Theorem~\ref{thm:not_large_not_sparse} and the general Theorem~\ref{thm:master} below) and comparing our study to the literature. A detailed discussion of our proof strategy and of the structure of the paper is given in \S \ref{sec:main_results_and_proof_strategy}. The proofs of the main theorems can be found in \S \ref{sec:proof_of_main_theorems}; they are based on a set of propositions that are stated in \S \ref{sec:main_results_and_proof_strategy} and proved in \S \S \ref{sec:divisors of reciprocal polynomials}--\ref{sec:galois_group_of_A}.

Reciprocal polynomials are ubiquitous in number theory and related areas: for example, they are the `optimisers' for Lehmer's longstanding conjecture regarding the minimal Mahler measure of a polynomial, and they appear as the characteristic polynomials of symplectic matrices. In the large box model, Bhargava's theorem \cite{Bhargava} mentioned in the first paragraph has already been adapted to the reciprocal setting by Anderson, Bertelli and O'Dorney \cite{ABD} (see also the earlier paper \cite{DDS} in this direction). They obtained the following: if $m \geq 3$ is fixed and $A$ is sampled uniformly at random from the set of monic reciprocal polynomials in $\Z[T]$ all of whose coefficients are at most $H$ in absolute value, then
\begin{equation*}
\mb{P}\Big(\CG_A \neq \CC_2 \wr \CS_m\Big) = \mb{P}\Big(\CG_A = (\CC_2 \wr \CS_m) \cap \CA_{2m}\Big) + O\big(H^{-1}\big) \asymp_m \frac{\log{H}}{H}, \qquad H \to \infty.
\end{equation*}
Here, the implicit constants depend on $m$. Notably, the leading term of the asymptotic for $\mb{P}(\CG_A \neq \CC_2 \wr \CS_m)$ comes from the maximal alternating subgroup $(\CC_2 \wr \CS_m) \cap \CA_{2m}$ and contains no contribution of reducible polynomials. In contrast, it is originally a conjecture by Van der Waerden \cite{VdW2} that in Bhargava's result about `standard' polynomials (where there are no further relations between the coefficients), the alternating group should not appear and instead $\mb{P}(\CG_A \neq \CS_n) = \mb{P}(A \textup{ is reducible}) + o(H^{-1})$ as $H$ goes to infinity. From the examples mentioned so far, one already observes that the alternating group plays a special role in probabilistic Galois theory. In the large box model, this is discussed in more detail in work by Bary-Soroker, {Ben-Porath} and Matei \cite{BBM}. Bary-Soroker and Goldgraber \cite{BG} proved that the alternating group can be ruled out in a hybrid model where both the degree and the coefficient box grow; in the setting of polynomials over $\F_q[t]$, Entin \cite{Entin} similarly showed that the Galois group over $\F_q(t)$ is the alternating group with low probability only. 

Other than reciprocal polynomials discussed in \cite{ABD} and here, polynomials with dependent coefficients have been considered from the Galois-theoretic perspective in a few other studies. In the model of fixed degree polynomials, Pham and Xu \cite{PX} generalised known results regarding irreducibility of polynomials over $\Q$ to a setting of significantly relaxed independence, uniformity and support assumptions. Another perspective is that of random polynomials coming from random matrices: Eberhard \cite{Eberhard} considered the characteristic polynomial $\phi$ of an $n \times n$ matrix with entries drawn independently according to a nontrivial measure $\mu$ on $\Z$ with finite support. He showed that $\phi$ is irreducible with high probability if $\mu$ is uniform modulo the product of four primes, generalising an earlier paper by Bary-Soroker and Kozma \cite{BK}. In addition, assuming the Riemann hypothesis for Dedekind zeta functions, Eberhard proved that $\CG_{\phi}$ contains $\CA_n$ with high probability, generalising the result of Breuillard and Varj\'u \cite{BV}. The latter result has been extended to the setting of symmetric matrices as well \cite{FJSS}. In \cite{Kim}, the authors give an application in quantum chaos --- related to the support of semiclassical measures for quantum cat maps --- of the following result they prove: when ordered by any norm, $100\%$ of the integral symplectic $2m \times 2m$ matrices $M$ have the property that the (reciprocal) characteristic polynomials of $M^k$ with $k \geq 1$ all have maximal Galois group.

Our approach builds on the aforementioned work of Bary-Soroker, Kozma and the second author \cite{BKK}, who employed rather different methods than those utilized in the large box model or in a setting where the extended Riemann hypothesis is assumed, drawing instead ideas from the anatomy of integers, permutations and polynomials, from $p$-adic Fourier analysis, and from asymptotic group theory. They showed analogues of Theorem~\ref{thm:simple_interval} and the other results presented below for standard monic polynomials (again, where all coefficients are drawn independently). Our proofs, outlined in \S \ref{sec:main_results_and_proof_strategy}, follow the same line of reasoning. We mention the main novelties in our work. In the Fourier-analytic part of the paper (\S\S \ref{sec:character_theory_and_special_residues}--\ref{sec:L_1_bounds}), the dependencies between the coefficients of $A$ form a serious technical obstacle. We use character theory of finite abelian groups to address this (see \S \ref{sec:character_theory_and_special_residues}). To obtain the Galois group results, we combine \cite{BKK} with an analysis of (1) the asymptotic subgroup structure of $\CC_2 \wr \CS_m$ by means of group cohomology, part of which was already described in \cite{ABD} (see \S \ref{sec:hyperoctahedral_group}); (2) the discriminant of random reciprocal polynomials to handle the group $(\CC_2 \wr \CS_m) \cap \CA_{2m}$, building on the first author's earlier paper \cite{Hokken} (see \S \ref{sec:discriminant_and_exceptional}); and (3) the Frobenius automorphism to handle the group $\CC_2 \times \CS_m$ (see \S \ref{sec:galois_group_of_A}).

Let us now discuss our other results. First, we state some conventions and introduce some notation, mostly following \cite{BKK}. In general, by saying that the polynomial $A$ in \eqref{eq:reciprocal polynomial} is a \emph{random monic reciprocal polynomial} we mean that $a_m=1$ and that the coefficients $a_0, a_1, \ldots, a_{m-1}$ are independent random variables, with each $a_j$ selected according to a fixed probability measure $\mu_j$ on the integers. 

\begin{definition}
\label{def:CRpm}
We denote by $\CR(m)$ the set of monic reciprocal polynomials $A \in \Z[T]$ of degree $2m$. In other words, $\CR(m)$ is the set of polynomials as in \eqref{eq:reciprocal polynomial} with integer coefficients and $a_m = 1$. Similarly, if $p$ is prime, then the set of such polynomials $A$ in the ring $\F_p[T]$ will be denoted by $\CR_p(m)$.
\end{definition}

Any sequence $\mu_0, \mu_1, \ldots, \mu_{m-1}$ of probability measures on the integers induces a probability measure on $\CR(m)$ via
\begin{equation}
	\label{eq:P on reciprocals}
\mb{P}_{\CR(m)}(A) \coloneqq \prod_{j=0}^{m-1} \mu_j(a_j),
\end{equation}
where it is understood that the $a_j$ refer to the coefficients of $A$ in \eqref{eq:reciprocal polynomial}. 
We will prove results such as those in Theorem~\ref{thm:simple_interval} for a wide range of varying measures $\mu_j$, subject to certain Fourier-analytic conditions. Our most general result, Theorem~\ref{thm:master} below, also specifies a lower bound for the rate of convergence that applies (for example) in the setting of Theorem~\ref{thm:simple_interval}. 

Since our general Theorem~\ref{thm:master} discussed below is so similar to the case $s=1$ in \cite[Theorem 7]{BKK}, our other results are, minor details aside, the exact analogues of results in \cite{BKK} --- such as Theorem~\ref{thm:not_large_not_sparse} below. To state this result, we first introduce the following standard notation. Given a probability measure $\mu$, we denote by $\supp(\mu)$ its support. Furthermore, when $\mu$ is supported on the integers, we set
\begin{equation}
\label{eq:norm_definition}
\norm{\mu}_q \coloneqq 
\begin{cases}
(\sum_{a \in \Z} \mu(a)^q)^{1/q} & \text{if } 1 \le q < \infty, \\
\sup_{a \in \Z} \mu(a) & \text{if } q = \infty.
\end{cases}
\end{equation}
Given a prime $p$ and a measure $\mu$, we will sometimes consider its projection $\bar{\mu}_p$ to $\F_p$, i.e., the measure
\begin{equation}
\label{eq:def mu mod p}
\bar{\mu}_p \colon \F_p \to [0,1], \quad \bar{\mu}_p(a) = \sum_{\substack{b \in \Z \\ b \equiv a \bmod{p}}} \mu(b).
\end{equation}
When the measure is $\mu=\mu_j$, we write $\bar{\mu}_{j,p}$ instead to avoid confusion. The norms $\norm{\bar{\mu}_p}_q$ for $q \in [1, \infty]$ are defined as in \eqref{eq:norm_definition}, with the domain of the sum and of the supremum equal to $\F_p$ instead of $\Z$.

\begin{introtheorem}
\label{thm:not_large_not_sparse}
Fix $H \ge 3$. Let $m \ge 3$ be an integer and suppose $\mu_j = \mu$ is a fixed probability measure for all $j = 0, 1, \ldots, m-1$, with the following properties:
\begin{enumerate}
  \item (bounded support) $\supp(\mu) \subset [-H,H]$;
  \item (support not too sparse) $\norm{\mu}_2^2 \le \min \{H^{-4/5}, m^{1/16}/H\}/(\log{H})^2$.
\end{enumerate}
Then there are absolute constants $c, C, H_0 > 0$ such that $H \geq H_0$ implies
\begin{equation*}
\mb{P}_{\CR(m)}(A \textup{ is irreducible}) \ge 1- Cm^{-c}.
\end{equation*}
\end{introtheorem}

Theorems~\ref{thm:simple_interval} and \ref{thm:not_large_not_sparse} are both consequences of Theorem~\ref{thm:master}, which is our most general result. (We will refrain from rehashing other results of \cite{BKK}, such as their Theorems 3 and 4, in the setting of reciprocal polynomials; the interested reader could derive them from the general Theorem~\ref{thm:master} below in the same way as Theorems 3 and 4 in \cite{BKK} follow from their Theorem 7.) For this, define the function $e(x) \coloneqq \exp(2 \pi \ii x)$ and recall that the Fourier transform of a measure $\mu$ on the integers is defined as
\begin{equation*}
\hat{\mu}\colon \R/\Z \to \C, \qquad \hat{\mu}(\theta) = \sum_{a \in \Z} \mu(a) e(a\theta).
\end{equation*}

Our most general result is then the following.
\begin{introtheorem}
\label{thm:master}
Let $P$ be the product of four distinct primes. Consider an integer $m \ge P^4$ and probability measures $\mu_0, \mu_1, \ldots, \mu_{m-1}$  on the integers, all supported in $[-H,H]$ with $H \geq 2$. Let $\alpha \in (0, 1/2)$ and $B\ge1$ be two constants. Suppose that the following hold: 
\begin{enumerate}
  \item (bounded support) $H \le Bm^{1/2-\alpha}$;
  \item (controlled Fourier transform modulo four primes) We have
  \begin{equation*}
  \frac{1}{\sqrt{Q}}\sum_{k=0}^{Q-1} \abs{\hat{\mu}_j(k/Q + \ell/R)} \le 1-m^{-1/10}
  \end{equation*}
  for all $j=0, \ldots, m-1$ and all integers $Q,R,\ell$ such that $QR=P$ and $Q>1$.
\end{enumerate}
Then there are constants $c, C > 0$ depending at most on  $\alpha$ and $B$  such that
\begin{equation}
\label{eq:A is irreducible}
\mb{P}_{A \in \CR(m)}(A \textup{ is irreducible}) \ge 1- Cm^{-c}.
\end{equation}
Moreover, there are constants $c', C' > 0$ depending at most on $\alpha$ and $B$ such that
\begin{equation}
\label{eq:A has large Galois group}
\mb{P}_{A \in \CR(m)}\Big(\CG_A \in \big\{\CC_2 \wr \CS_m, \CC_2 \wr \CA_m, G_2\big\}\Big) \ge 1- C'm^{-c'},
\end{equation}
where $G_2$ denotes the group defined in \eqref{eq:G2}.
\end{introtheorem}

A more precise bound, given specific choices of the $\mu_j$, on the probability that $A$ is irreducible and $\CG_A \leq (\CC_2 \wr \CS_m) \cap \CA_{2m}$ can be derived using Proposition~\ref{prop:square_discriminant}. The probability that $A$ is irreducible and $\CG_A \leq \CC_2 \times \CS_m$ is bounded in Lemma~\ref{lem:Galois_C2xSm_2}.

Surprisingly, our irreducibility results are of virtually the same strength as those in \cite{BKK}, despite concerning only a sparse subsequence of all polynomials of degree $n$. This is thanks to the rigid divisor structure of reciprocal polynomials (see Lemma~\ref{lem:factorisation options} below), which dictates that a reducible reciprocal polynomial either has a nontrivial reciprocal divisor or factors in an exceptional, structured way that occurs only rarely\footnote{This is markedly different from the direct analogue over the integers: the divisors of a \emph{palindromic integer} are typically not palindromic themselves.} (see \S \ref{sec:discriminant_and_exceptional}).
On the other hand, all our results are ``all or (almost) nothing'': for a given sequence of measures $\mu_0, \mu_1, \ldots, \mu_{m-1}$ we can either show that $A$ is irreducible with high probability, or we can only rule out divisors of degree $\ll m^{1/10}$, following Konyagin's classical argument \cite{Konyagin}. In contrast, \cite[Theorem~7]{BKK} shows that in the standard polynomial case, divisors of degree $\le \theta m$ for some constant $\theta > 0$ (depending on the measures) can also be precluded in many cases. The source of this difference lies in a trick employed in the proof of the $L^1$ bounds in \cite[Lemma~6.3]{BKK}, the analogue of which fails to hold in our setting. This is explained in more detail in Remark~\ref{rem:why not s>1}.

\begin{remark}
Reciprocal polynomials, sometimes with slightly different definitions, are occasionally called palindromic, self-reciprocal, {(self-)}inversive or symmetric polynomials.
\end{remark}

\subsection*{Acknowledgements}
Many thanks to Gunther Cornelissen for helpful conversations and feedback on earlier versions of this manuscript, as well as to the referee for a very careful reading of the paper, and for catching and correcting a subtle issue in an earlier version of Lemma \ref{lem:subgroups projecting onto S_m or A_m}. 

We used ChatGPT to check the paper for possible errors. This revealed various minor issues that we corrected. The most important correction was to Lemmas \ref{lem:minor residues aux} and \ref{lem:L1 bound}.
 
D.~H. is supported by the Dutch Research Council (NWO) through the grant OCENW.M20.233.
This project was originally conceived during a visit of D.~H. to Université de Montréal in October 2023, supported by the Courtois Chair. He would like to thank his hosts for the pleasant stay.

D.~K. is supported by the Courtois Chair II in fundamental research, by the Natural Sciences and Engineering Research Council of Canada (RGPIN-2024-05850), by the Fonds de recherche du Qu\'ebec - Nature et technologies (2025-PR-345672) and by the Simons Foundation.

\subsection*{Notation}
We utilize the standard Vinogradov and big-O asymptotic notation, with $f \ll g$, $g \gg f$ and $f = O(g)$ on some domain $X$ all meaning that there exists a constant $C=C(f,g,X)>0$ such that $\abs{f(x)} \leq Cg(x)$ for all $x \in X$. If the constant depends on more than just the datum $f$, $g$ and $X$, this is indicated by a subscript such as $f \ll_{\delta} g$. The notation $f \asymp g$ means both $f \ll g$ and $g \ll f$ hold.

\section{Overview and proof strategy}
\label{sec:main_results_and_proof_strategy}

\subsection*{Overview and structure of the paper} 
Globally, we follow the same strategy as in the work of Bary-Soroker, Kozma and the second author \cite{BKK}. We will borrow their notation unless stated explicitly. We briefly summarise their approach to prove irreducibility for `standard' monic random polynomials $A(T)=a_0+a_1T+\cdots+a_{m-1}T^{m-1}+T^m$:

\begin{enumerate}
  \item The first step is to show that $A$ cannot have a divisor of small degree, say smaller than $m^{1/10}$, using a classical approach due to Konyagin \cite{Konyagin}. This approach combines the random walk-type properties of the polynomial $A$ to exclude cyclotomic divisors, with an analysis based on the Mahler measure for the noncyclotomic divisors. The result is then obtained with a union bound.
  \item Next, for the divisors of degree $k>m^{1/10}$, they start with the basic observation that if $A$ has a divisor of degree $k$, then so does $A_p \coloneqq (A \bmod{p}) \in \F_p[T]$. They then show that $A_p$ is approximately uniformly distributed among the polynomials of degree $m$ in $\F_p[T]$, in the sense that, for $D \in \F_p[T]$, the approximation $\mb{P}(D \textup{ divides } A_p) \approx p^{-\deg{D}}$ holds very well \emph{on average}.
    \item The divisor structure of polynomials over $\F_p$ is well-understood: the probability that a `typical' random polynomial of degree $m$ in $\F_p[T]$ has a divisor of degree $k$ is approximately $k^{-1+\log{2}} \approx k^{-0.3}$. Since $A_p$ is approximately uniform, this also holds for $A_p$. However, the probability is too large to simply sum over all $k>m^{1/10}$.
  \item To circumvent this problem, they use not just one, but multiple primes as in the predecessor paper of Bary-Soroker and Kozma \cite{BK}. (Naively one would expect the factorisations of $A_p$, as $p$ varies, to be independent; this indeed holds by the Chinese remainder theorem for the primes $p$ dividing $H$ when the coefficients of $A$ are selected uniformly at random from the set $\{1,2,\ldots,H\}$.) Thus one needs to control the \emph{joint} distribution of $(A_p)_{p \in \CP}$ for some set $\CP$ of four primes, as $k^{4(-1+\log{2})} \approx k^{-1.2}$; summing this over all $k>m^{1/10}$ yields a quantity that goes to $0$ as $m$ goes to infinity. 
\end{enumerate}

To adjust this strategy in our context, we introduce two pivotal concepts. First, the \emph{trace polynomial} of $A$, denoted $A_{\SR}$, is the unique polynomial of degree $m$ with the property
\begin{equation*}
A(T) = T^m A_{\SR}(T+T^{-1}).
\end{equation*}
The mapping $A \mapsto A_{\SR}$ is multiplicative, setting up a correspondence between the divisors of $A_{\SR}$ and the \emph{reciprocal} divisors of $A$; see \S \ref{sec:the_trace_polynomial} below. Secondly, the \emph{reversal polynomial} of a polynomial $F$ is 
\begin{equation*} 
F_{\ms{rev}}(T) \coloneqq T^{\deg{F}} F(T^{-1}),
\end{equation*}
so that $F = F_{\ms{rev}}$ when $F$ is reciprocal. (The notation comes from \cite{CC}.)

A fundamental fact, proved in \S \ref{sec:divisors of reciprocal polynomials}, is the following: if $A$ is reducible, then either $A$ has a nonconstant reciprocal divisor of degree $\le m$, or $A = \pm I\cdot I_{\ms{rev}}$ for some monic irreducible nonreciprocal polynomial $I$. We deal with small-degree reciprocal divisors (step 1) in \S \ref{sec:small degree factors}, and with the `exceptional factorisation' $A = \pm I \cdot I_{\ms{rev}}$ in \S \ref{sec:discriminant_and_exceptional}. The Fourier analysis (steps 2 and 4) is carried out in \S\S \ref{sec:Euclidean division}--\ref{sec:L_1_bounds}, and \S \ref{sec:anatomy} concerns the anatomy of reciprocal polynomials (step 3). These sections cover the proof of irreducibility of $A$ with high probability. Lastly, \S \S \ref{sec:hyperoctahedral_group}--\ref{sec:galois_group_of_A} cover the Galois group. Some preliminary Galois-theoretic work is done in \S \ref{sec:discriminant_and_exceptional}, where we study the discriminant of $A$.

We compare our proof of irreducibility with that in \cite{BKK}, other than the exceptional factorisation that we already mentioned. Let $D \in \CR(k)$ be a reciprocal polynomial of degree $2k$. We handle the case $k \le m^{1/10}$ with a simple adaptation of the approach in \cite{BKK}. Step 2, where $k > m^{1/10}$, requires considerable effort: it turns out that the main contribution to $\mb{P}(D \textup{ divides } A_p)$ does not come from just the zero residue as in \cite{BKK}, but rather from a linear subspace of $\F_p[T]/(D)$ of dimension $k$, which is of character-theoretic origin. This happens because of the dependencies between the coefficients of reciprocal polynomials, and is a significant hurdle in establishing the $L^1$ bounds (\S \ref{sec:L_1_bounds}). The main results that help overcome this are covered in \S\S \ref{sec:Euclidean division}--\ref{sec:character_theory_and_special_residues}. The proof of the $L^{\infty}$ bounds (\S \ref{sec:L infinity bounds}) is also quite different, with our approach boiling down to linear algebra and some classical results about Chebyshev polynomials. For the anatomy of reciprocal polynomials (step 3), we can rely --- after a trick --- on the work carried out in \cite[\S \S 8-10]{BKK}. 

For the Galois theory, we have a different approach from that in the paper \cite{BKK}, although we do require the result of their work. Starting in \S \ref{sec:hyperoctahedral_group}, we carry out an analysis of the subgroup structure of the \emph{hyperoctahedral group} $\CC_2 \wr \CS_m$, which is the largest Galois group $A$ can have. We then show in \S \ref{sec:galois_group_of_trace_pol} that the Galois group of $A$ is either one of five `large' groups --- those that are compatible with the results of \cite{BKK} and have bounded index in $\CC_2 \wr \CS_m$ as $m$ goes to infinity --- or is contained, up to conjugation, in one `small' group: $\CC_2 \times \CS_m$. We can exclude two of the five large groups by using our results on the discriminant, which we already consider in \S \ref{sec:discriminant_and_exceptional}. Lastly, in \S \ref{sec:galois_group_of_A}, we use an idea about the modulo $p$ factorisation of reciprocal polynomials that we believe is new to preclude the small group $\CC_2 \times \CS_m$ and its subgroups as possible Galois groups.

As mentioned in the introduction, the present section provides the precise statements proved in \S \S \ref{sec:divisors of reciprocal polynomials}--\ref{sec:galois_group_of_A}. From these statements, we deduce in the remaining \S \ref{sec:proof_of_main_theorems} the main Theorems~\ref{thm:simple_interval}, \ref{thm:not_large_not_sparse} and \ref{thm:master} that were presented in \S \ref{sec:introduction}.

\subsection*{Preliminaries}

A preliminary result that greatly facilitates the analysis carried out in this paper is the following basic proposition regarding the divisor structure of a reciprocal polynomial.

\begin{lemma}[Divisors of reducible reciprocal polynomials]
  \label{lem:factorisation options}
  Suppose that $A \in \Z[T]$ is a monic reducible reciprocal polynomial of degree $2m$. Then either $A$ has a monic nontrivial reciprocal divisor in $\Z[T]$ of degree $2k\le m$, or there exists a monic irreducible  polynomial $I \in \Z[T]$ of degree $m$ and a sign $u\in\{\pm1\}$ such that $A = u I \cdot I_{\ms{rev}}$.
\end{lemma}

Lemma~\ref{lem:factorisation options} is a direct consequence of the following more general proposition.

\begin{proposition}[Prime factorization of reciprocal polynomials]
\label{prop:factorisation options}
Let $K$ be a field. Suppose that $A(T)\in K[T]$ is a nonzero reciprocal polynomial. Then there exist irreducible polynomials $I_1,\dots,I_r$ and $J_1,\dots,J_s$ over $K$, integers $a,b\ge0$ and an element $u\in K^\times$ such that:
\begin{enumerate}[label=(\alph*)]
	\item each $I_j$ is reciprocal;
	\item each $J_j$ is coprime to $J_{j,\ms{rev}}$;
	\item $A(T)= u \cdot (T-1)^{2a}(T+1)^{2b} \prod_{i=1}^r I_i(T) \prod_{j=1}^s J_j(T)J_{j,\ms{rev}}(T)$.
\end{enumerate}
\end{proposition}

Proposition~\ref{prop:factorisation options} and Lemma~\ref{lem:factorisation options} are shown in \S \ref{sec:divisors of reciprocal polynomials}.

\subsection*{Small degree reciprocal factors and the exceptional factorisation}

As explained in the first subsection, we will proceed by showing that $A$ is unlikely to have a divisor of degree $k$, for any $k\in\{1,2,\dots,m\}$. Since only $p^k$ of the $p^{2k}$ monic polynomials in $\F_p[T]$ of degree $2k$ are reciprocal, Proposition~\ref{prop:factorisation options} allows us to consider a much smaller set of possible factors.

We start by proving that $A$ is unlikely to have a reciprocal factor of small degree.

\begin{proposition}[Small degree factors]
\label{prop:small degree factors}
Let $m \in \Z_{\ge 1}$ and let $\mu_0, \ldots, \mu_{m-1}$ be a sequence of probability measures on the integers. Suppose
\begin{itemize}
    \item[(1)] (bounded support) $\supp(\mu_j) \subset [-\exp(m^{1/3}), \exp(m^{1/3})]$ for all $j \ge 0$,
    \item[(2)] (anti-concentration) $\norm{\mu_j}_\infty \le 1-m^{-1/10}$ for all $j \ge 0$.
\end{itemize}
Then
\begin{equation*}
\mb{P}_{A\in\CR(m)}\Big(A \textup{ has a reciprocal divisor } D \in \CR(k) \textup{ with } 1 \le k \le m^{1/10}\Big) \ll m^{-7/20}.
\end{equation*}
\end{proposition}
The proof of Proposition~\ref{prop:small degree factors} is given in \S \ref{sec:small degree factors}. The method goes back to Konyagin \cite{Konyagin}. Here, we closely follow the adaptation of Konyagin's ideas laid out in \cite[\S 7]{BKK}.

Next, we show that $A$ is unlikely to have the \emph{exceptional factorisation} $A = \pm I\cdot I_{\ms{rev}}$, where $I \in \Z[T]$ is a monic, irreducible, nonreciprocal polynomial.

\begin{proposition}[Exceptional factorisation]
\label{prop:exceptional}
Let $m \in \Z_{\ge 2}$. Let $\mu_0, \ldots, \mu_{m-1}$ be a sequence of probability measures supported on $[-H,H]$ for some $H \geq 2$. Let $\epsilon > 0$ such that $\norm{\mu_j}_{\infty} \leq 1-\epsilon$ for all $j = 0, 1, \ldots, m-1$. Then
\begin{equation}
\label{eq:exceptional_factorisation}
\mb{P}_{A\in\CR(m)}\bigg(  \begin{array}{l} \exists \textup{ monic, irreducible, nonreciprocal } I \in \Z[T] \\ 
	\text{such that}\ A = \pm I\cdot I_{\ms{rev}} 
\end{array}
	\bigg) \ll \frac{H (\log{m})^{3/2}}{\epsilon \sqrt{m}}.
\end{equation}
\end{proposition}

Proposition~\ref{prop:exceptional} is shown in \S \ref{sec:discriminant_and_exceptional}. In a nutshell, the proof proceeds by observing that the expression $A(1)A(-1)$ is a square up to sign if $A$ factors as $A = \pm I\cdot I_{\ms{rev}}$, and then showing that this is unlikely to occur by writing $A(1)A(-1)$ as the difference of two (squared) random walks, with steps given by the coefficients of $A$.

\subsection*{Large degree reciprocal factors: approximate equidistribution and the anatomy of reciprocal polynomials}

\subsubsection*{Perfect equidistribution modulo four primes}

Let us see how what we laid out so far helps to prove irreducibility in a particularly `fair' setting. Let $\CP = \{p_1, \ldots, p_r\}$ be a set of $r \ge 4$ primes. Let $H$ be a positive integer divisible by each of the $p_j$, and suppose $A$ has coefficients drawn independently and uniformly at random from $\{1, \ldots, H\}$. Fix $p \in \CP$ and write $A_p \coloneqq (A \bmod{p}) \in \F_p[T]$. Then $A_p$ is perfectly equidistributed among the reciprocal polynomials of degree $2m$ in $\F_p[T]$. By the details of the correspondence between reciprocal polynomials and their trace polynomials (Lemma~\ref{lem:trace pol iso of vs} below), the trace polynomial $(A_p)_\SR = (A_\SR)_p$ is then perfectly equidistributed among the polynomials of degree $m$ in $\F_p[T]$. But then one can leverage the earlier work \cite{BK} to prove that the polynomial $A_\SR$ is irreducible with high probability. As a result, by Lemma~\ref{lem:factorisation options} and the reciprocal polynomial-trace polynomial correspondence, we conclude that, with high probability, the polynomial $A$ is irreducible or factors as $A = \pm I\cdot I_{\ms{rev}}$. But this last option is unlikely by Proposition~\ref{prop:exceptional}, thus concluding the proof that $A$ is irreducible with high probability. In this approach, note that it is essential that $r \geq 4$: this is an intrinsic limit to the method of both \cite{BK} and \cite{BKK}, and will thus remain essential in this paper as well.

\subsubsection*{Approximate equidistribution and independence}

Before proceeding, we will first introduce some notation, borrowing from \cite{BKK} whenever possible. We denote by boldface letters tuples of polynomials in $\F_{\CP}[T] \coloneqq \prod_{p \in \CP} \F_p[T]$. For example, the tuple $\BD = (D_p)_{p \in \CP}$ consists of polynomials $D_p \in \F_p[T]$, for each $p \in \CP$. In the case of boldface $\BA$, this will typically denote the specific tuple $(A \bmod{p})_{p \in \CP}$. 
We then set
\begin{equation*}
\norm{D_p}_p \coloneqq p^{\deg{D_p}}, \qquad \norm{\BD}_{\CP} \coloneqq \prod_{p \in \CP} \norm{D_p}_p.
\end{equation*}
The phrase `$\BD$ divides $\BA$' will mean that $D_p$ divides $A_p$ in $\F_p[T]$, for all $p \in \CP$. Suppose $\BD$ is a \emph{reciprocal} tuple, that is, each $D_p$ is a reciprocal polynomial. Furthermore, assume that each $D_p$ is of degree $\leq 2m$. In the setting of perfect equidistribution outlined above, observe that 
\begin{equation*}
\mb{P}\big(D_p \text{ divides } A_p\big) = \mb{P} \big((D_p)_{\SR} \text{ divides } (A_p)_{\SR}\big) = \frac{1}{\norm{(D_p)_{\SR}}_p}  = \frac{1}{\norm{D_p}^{1/2}_p}
\end{equation*}
for each $p \in \CP$. In other words,
\begin{equation*}
\mb{P}(\BD \text{ divides } \BA) = \frac{1}{\norm{\BD}^{1/2}_{\CP}}.
\end{equation*}
The power $1/2$ in the preceding line is somewhat surprising at first sight, but it stems from the correspondence between the reciprocal polynomial $A$ and its trace polynomial $A_{\SR}$, which is of half the degree of $A$ and has a divisor if and only if $A$ has a reciprocal divisor.

Now, in order to establish irreducibility results in a setting away from perfect equidistribution, we follow the ways of Fourier analysis --- as in \cite{BKK}. Let us first sketch an, at this point, tempting route that however fails. In \cite{BKK}, Bary-Soroker, Kozma and the second author established for \emph{standard polynomials} $B = b_0 + \cdots + b_{m-1}T^{m-1} + T^m$ the following. If
\begin{equation*}
\mb{P}(\BG \text{ divides } \BB) \approx \frac{1}{\norm{\BG}_{\CP}}
\end{equation*}
\emph{on average} over all tuples $\BG = (G_p)_{p \in \CP}$ with $\deg{G_p} \le m$ for all $p$, and the measures according to which the $b_j$ are selected are not too concentrated modulo each $p$, then $B$ is unlikely to have a factor of large degree. (The precise conditions are not important at this point.) Thus, in our setting, it seems reasonable to first pass to the trace polynomial by `setting' $\BB = \BA_{\SR}$, and then import the results of \cite{BKK}. However, there is one major problem: the coefficients of $A_{\SR}$ are generally speaking not independent (see Lemma~\ref{lem:trace polynomial coefficients}), that is, there exists a choice of measures $\mu_0, \ldots, \mu_{m-1}$ such that the induced probability distribution on the trace polynomial $\BA_{\SR}$ cannot be realised by selecting the coefficients of $\BA_{\SR}$ independently according to some other sequence of measures $\mu_0', \ldots, \mu_{m-1}'$.
But the independence assumption on the coefficients of the polynomial $B$ underlying $\BB$ is indispensable in \cite{BKK}. More precisely, the approach taken in \cite{BKK} fails without this assumption at the very start of the Fourier-analytic approach, in their Lemma 4.1. 

Thus, instead of \emph{first} passing to the trace polynomial and then carrying out the Fourier analysis, here we perform the Fourier analysis on the reciprocals immediately. Of course, dependence between the coefficients cannot be avoided, since the coefficients of $A$ are already dependent by the reciprocity assumption. Overcoming this (comparatively mild) dependency is the main technical novelty in this work, with the core of it described in \S \ref{sec:character_theory_and_special_residues}. As a consequence, we show that
\begin{equation*}
\mb{P}(\BD \text{ divides } \BA) \approx \frac{1}{\norm{\BD}_\CP^{1/2}}
\end{equation*}
holds on average over all vectors $\BD = (D_p)_{p \in \CP}$ with $\deg{D_p} \le m$ for all $p$, under the assumption that the Fourier transforms of the measures $\mu_j$ are sufficiently well-behaved.
Since
\begin{equation*}
\mb{P}(\BD \text{ divides } \BA) = \mb{P}(\BD_{\SR} \text{ divides } \BA_{\SR}) \quad \text{and} \quad \frac{1}{\norm{\BD}_\CP^{1/2}} = \frac{1}{\norm{\BD_{\SR}}_\CP}
\end{equation*}
this implies
\begin{equation*}
\mb{P}(\BD_\SR \text{ divides } \BA_\SR) \approx \frac{1}{\norm{\BD_\SR}_\CP}
\end{equation*} 
on average as well, allowing us to import the results from \cite{BKK} --- which at this point \emph{is} possible (see \S \ref{sec:anatomy}).

In fact, we show that $\BA$ is approximately equidistributed among the `admissible’ residue classes modulo $\BD$. The precise result is the following. Let $\mb{P}_{\CR_p(m)}$ denote the probability measure induced by $\mb{P}_{\CR(m)}$ after reduction modulo $p$, that is
\begin{equation*}
\mb{P}_{\CR_p(m)}(A) \coloneqq \prod_{j=0}^{m-1} \bar{\mu}_{j,p}(a_j),
\end{equation*}
where we recall that $\bar{\mu}_{j,p}$ is defined in \eqref{eq:def mu mod p}. Similarly, set
\begin{equation*}
\mb{P}_{\CR_\CP(m)}(\BA) \coloneqq \prod_{j=0}^{m-1} \Big( \sum_{\substack{a \in \Z \\ a \equiv a_{j,p} \bmod{p} \ \forall p\in\CP}} \mu_j(a) \Big)
\end{equation*}
where we write $a_{j,p}$ for the coefficient of the monomial $T^{m-j}$ of $A_p$. We shall denote $\BA \equiv \BC \bmod{\BD}$ to mean that $D_p$ divides $A_p - C_p$ for all $p \in \CP$. By $R_m(\BD) \subset \F_\CP[T]/(\BD) \coloneqq \prod_{p \in \CP} \F_p[T]/(D_p)$ we (essentially) denote the set of all possible residues that \emph{could} occur as the residue of a tuple of reciprocal polynomials $\BA \in \CR_{\CP}(m) \coloneqq \prod_{p \in \CP} \CR_p(m)$ modulo a reciprocal polynomial $\BD \in \CR_{\CP}(\Bk) \coloneqq \prod_{p \in \CP} \CR_{p}(k_p)$, where $\Bk = (k_p)_{p \in \CP}$; the more precise definition, given in \eqref{eq:def RmD tuple} and Definition~\ref{def:RmD}, is slightly different and requires some work. 
Lastly, set
\begin{equation}
\label{eq:Delta def}
\Delta^\SR_\CP(m;k) \coloneqq \sum_{\substack{\BD \in \CR_{\CP} \\ \deg{\BD} \le 2k \\2 \in \CP \implies T^2+1 \nmid D_2}} \max_{\BC \in R_m(\BD)} \Bigg\vert \mb{P}_{\BA \in \CR_{\CP}(m)}\big(\BA \equiv \BC \bmod{\BD}\big) - \frac{1}{\norm{\BD}_{\CP}^{1/2}} \Bigg\vert,
\end{equation}
where
\[
\deg{\BD}\coloneqq \max_{p\in\CP}\deg(D_p). 
\]
We also adopt the convention $\Delta_p^{\SR}(m;k) \coloneqq \Delta_{\{p\}}^{\SR}(m;k)$. The quantity $\Delta^\SR_\CP$ measures the extent to which the probability measure $\mb{P}_{\CR(m)}$ fails to be equidistributed modulo the primes in $\CP$ for reciprocal moduli up to a given degree. The requirement that $T^2+1$ does not divide $D_2$ whenever $2 \in \CP$ is technical; see \S \ref{sec:L infinity bounds}. It is equivalent to the condition that $T$ does not divide the trace polynomial of $D_2$. This is reminiscent of the (stronger) condition that $T$ does not divide $D_p$ for \emph{any} $p \in \CP$, which is stipulated in \cite[(2.6)]{BKK}.

We then show that $\Delta^\SR_\CP$ is very small when the measures $\mu_0, \mu_1, \ldots$ have good Fourier-theoretic properties.
\begin{proposition}[Controlled Fourier transform implies approximate equidistribution] 
\label{prop:Delta bound}
Let $\CP = \{p_1, \ldots, p_r\}$ be a set of distinct primes, and set $P = p_1 \cdot \ldots \cdot p_r$. 
In addition, consider an integer $m \ge P^4$ and a sequence $\mu_0, \ldots, \mu_{m-1}$ of probability measures on the integers for which there is $\gamma \ge 1/2$ such that
\begin{equation*}
\sum_{k=0}^{Q-1} \abs{\hat{\mu}_j(k/Q + \ell/R)} \le \big(1-m^{-1/10}\big) \cdot Q^{1-\gamma}
\end{equation*}
for all $j=0, \ldots, m-1$ and all integers $Q,R,\ell$ such that $QR=P$ and $Q>1$. Then
\begin{equation*}
\Delta^\SR_\CP\big(m; \gamma m + m^{0.88}\big) \ll_r  e^{-m^{1/10}}  .
\end{equation*}
\end{proposition}
Proposition~\ref{prop:Delta bound} will be proved in \S \ref{sec:L_1_bounds}; the proof combines results from \S \S \ref{sec:Euclidean division}--\ref{sec:L_1_bounds}.

\subsubsection*{Conclusion: precluding large degree reciprocal factors}

As in \cite{BKK}, we define the constant
\begin{equation*}
\lambda_0 \coloneqq \frac{1}{4-4\log{2}} = 0.81\ldots.
\end{equation*}

The final step in the proof of irreducibility is to rule out reciprocal factors of large degree under the assumption of approximate equidistribution. We prove the following analogue of \cite[Proposition~2.2]{BKK}.

\begin{proposition}[Large degree factors]
\label{prop:large degree factors}
Let $m \in \Z_{\ge 1}$ and $\epsilon \in (0, 1/100]$. Let $\CP$ be a set of four primes. Suppose $\mu_0, \ldots, \mu_{m-1}$ is a sequence of probability measures on the integers satisfying the following:
\begin{enumerate}
  	\item (mod-$\CP$ approximate equidistribution) $\Delta^\SR_\CP(m;m/2 + m^{\lambda_0 + \epsilon}) \leq m^{-30}$.
  	\item   (mod-$\CP$ anti-concentration) $\sup_{0 \leq j < m} \norm{\bar{\mu}_{j,p}}_{\infty} \leq 1-m^{-\epsilon/200}$ for all $p \in \CP$.
\end{enumerate}
Then there are constants $c, C > 0$ depending at most on $\epsilon$ such that
\begin{equation*}
\mb{P}_{A \in \CR(m)} \Big(A \textup{ has a reciprocal divisor } D \in \CR(k) \textup{ with } k \in [m^{1/10}, m/2] \Big) \le Cm^{-c}.
\end{equation*}
\end{proposition}

Proposition~\ref{prop:large degree factors} is proved in \S \ref{sec:anatomy}. The idea behind the proof is anatomy of polynomials: an analysis of the typical multiplicative structure of the polynomial $A$ reduced modulo the primes in $\CP$. As alluded to above, at this stage we \emph{can} use the work of \cite{BKK}, because the anatomical part of their argument is insensitive to the dependencies between the coefficients of $A_{\SR}$. In this way, we bypass almost all of the hard work that had to be carried out in \cite{BKK} at this point (their \S \S 8-9) and swiftly prove irreducibility after translating their results to our setting.

It is not always clear in the anatomy part of \cite{BKK} that the coefficients could be dependent variables. For example, their Proposition~2.2 explicitly mentions the sequence of measures $\mu_0, \mu_1, \ldots$, implying that the coefficients are chosen independently. In \S \ref{sec:anatomy}, we state a version of that proposition for a general probability measure $\mb{P}_{\CM(m)}$ on the set $\CM(m)$ of monic polynomials of degree $m$.
In a similar vein, even though the probability measure $\mb{P}=\mb{P}_{\CR(m)}$ in our Proposition~\ref{prop:large degree factors} is induced by the sequence $\mu_0, \ldots, \mu_{m-1}$, this is not necessary. A more general version of Proposition~\ref{prop:large degree factors}, where condition (2) is replaced by a suitable alternative, is given in Proposition~\ref{prop:large degree factors alternative}.

\subsection*{Galois theory}

In \S\S \ref{sec:hyperoctahedral_group}--\ref{sec:galois_group_of_A}, we study the Galois group $\CG_A$ of $A$ over $\Q$. Some preparatory work is already done in \S \ref{sec:discriminant_and_exceptional}, where we bound the probability that the polynomial discriminant $\Delta(A)$ of $A$ is a nonzero square.

\subsubsection*{The Galois group}

The zeros of a squarefree, monic reciprocal polynomial $A$ come in pairs: if $\alpha$ is a zero of $A$, then $\alpha^{-1}$ is also a zero of $A$. Any element $\tau$ of the Galois group $\CG_A$ of $A$ over $\Q$, being a field automorphism, must respect this property as $\tau(\alpha^{-1}) = \tau(\alpha)^{-1}$. Let $\alpha_1, \alpha_1^{-1}, \ldots, \alpha_m, \alpha_m^{-1}$ denote the zeros of $A$. The argument above shows that $\CG_A$ permutes the \emph{blocks} $\{\alpha_j, \alpha_j^{-1}\}$. In general, the group of symmetries of the set $\{-m, \ldots, -2, -1, 1, 2, \ldots, m\}$ that permutes the pairs $\{k, -k\}$ is the \emph{hyperoctahedral group} $\CC_2 \wr \CS_m$. Thus the Galois group of $A$ is contained in $\CC_2 \wr \CS_m$. This group is realised as a \emph{permutational wreath product}, and its construction is explained in \S \ref{sec:hyperoctahedral_group}. 

Our main result states that under the condition of approximate equidistribution modulo \emph{one} prime, the probability that $A$ is irreducible and has Galois group containing at least one of two index-$2$ subgroups of the hyperoctahedral group tends to $1$ as $m$ tends to infinity. More precisely, we have the following.

\begin{proposition}
\label{prop:Galois}
Consider probability measures $\mu_0, \mu_1, \ldots, \mu_{m-1}$  on the integers, all supported on $[-H,H]$ with $H \geq 2$. Fix a prime $p$ and real numbers $B\ge1$, $\alpha \in (0,1/2)$ and $\epsilon \in (0, 1/100]$. Suppose the following hold:
\begin{enumerate}
	\item (bounded support) $H \le Bm^{1/2-\alpha}$; 
  \item (mod-$p$ approximate equidistribution) $\Delta^{\SR}_p(m, m/2 + m^{\lambda_0+\epsilon}) \le m^{-10}$;
  \item (mod-$p$ anti-concentration) $\sup_{0 \le j < m} \norm{\bar{\mu}_{j,p}}_{\infty} \le 1-(\log{m})^{-2}$.
\end{enumerate}
Then there are constants $c, C > 0$ depending at most on $B$, $\alpha$ and $\epsilon$ such that
\begin{equation*}
\mb{P}_{A \in \CR(m)}\Big(\CG_A \not\in \{\CC_2 \wr \CA_m, \CC_2 \wr \CS_m, G_2\} \textup{ and } A \textup{ is irreducible}\Big) \le Cm^{-c}.
\end{equation*}
Here, the group $G_2$ is defined in \eqref{eq:G2}.
\end{proposition}
Proposition~\ref{prop:Galois} is proved in \S \ref{sec:galois_group_of_A}.
We use the results from \cite{BKK} to prove in \S \ref{sec:galois_group_of_trace_pol} that $A_{\SR}$ has, with high probability, either the symmetric group $\CS_m$ or the alternating group $\CA_m$ as Galois group. Then we need to understand what the possibilities for $\CG_A$ are, given that $\CG_{A_{\SR}} \in \{\CS_m, \CA_m\}$. In \S \ref{sec:hyperoctahedral_group}, we use group cohomology to show that $\CG_A$ is either one of five `large' groups (of index at most $4$ in $\CC_2 \wr \CS_m$) or conjugate to a subgroup of the `small' group $\CC_2 \times \CS_m$. Two of the five large groups can be excluded by studying the discriminant of $A$ (see below), and the small group can be excluded by analyzing the factorisation pattern of $A$ modulo a suitable prime. Another byproduct of the cohomological analysis is a \L uczak--Pyber theorem for the group $\CC_2 \wr \CS_m$ --- see Theorem~\ref{thm:Luczak--Pyber}.

\subsubsection*{The discriminant: ruling out alternating subgroups}

In \cite{BKK}, the authors show that the typical Galois group of a standard, monic polynomial $B = b_0 + \ldots + b_{m-1} T^{m-1} + T^m$ is $\CS_m$ or $\CA_m$ with high probability, when the $b_j$ are sampled independently according to measures satisfying Fourier-theoretic conditions that are essentially the same as ours in Proposition~\ref{prop:Delta bound}. 
The group $\CS_m$ is maximal possible; to rule out the alternating group and its subgroups as possible Galois groups of some squarefree polynomial $F \in \Z[T]$ of degree $m$, one can consider its discriminant $\Delta(F)$, since
\begin{equation*}
\Delta(F) \in \Z \text{ is a nonzero square if and only if } \CG_F \leq \CA_m.
\end{equation*}

For a standard random polynomial $B$, it is an open problem to quantify the probability that $\Delta(B)$ is a square; it is believed to be a highly unlikely event. For reciprocal polynomials however, the discriminant is a much more accessible invariant, and we prove the following:

\begin{proposition}
\label{prop:square_discriminant}
Let $\epsilon > 0$. Suppose the probability measures $\mu_0, \ldots, \mu_{m-1}$ have support contained in $[-H,H]$ with $H \ge 2$ and satisfy $\norm{\mu_j}_{\infty} \le 1-\epsilon$ for all $j \in [0,m-1]\cap \Z$. Then 
\begin{equation}
\label{eq:square_discriminant}
\mb{P}_{A \in \CR(m)}\big(\Delta(A) \textup{ is a nonzero square}\big) \ll \frac{H (\log{m})^{3/2}}{\epsilon \sqrt{m}}.
\end{equation}
\end{proposition}

The proof of Proposition~\ref{prop:square_discriminant} is given in \S \ref{sec:discriminant_and_exceptional}. Observe that the right-hand side of \eqref{eq:square_discriminant} is the same as the right-hand side of \eqref{eq:exceptional_factorisation}. Both bounds follow from a bound on the probability that $\pm A(1)A(-1)$ is a square. 

\begin{remark}
To establish this bound we follow \cite{Hokken}, in which the first author considered the case where $\mu_j=\mu$ for all $j$ with the law $\mu(1)=\mu(-1)=1/2$; these are known as \emph{Littlewood polynomials} or \emph{random Rademacher polynomials}. When $m \equiv 0, 3 \bmod{4}$, then \cite[Corollary~1.2]{Hokken} implies the asymptotic $\mb{P}_{\CR(m)}(\Delta(A) \textup{ is a nonzero square}) \sim c\log(m)/\sqrt{m}$ for an explicit constant $c>0$, whereas the probability is $0$ when $m \equiv 1,2 \bmod{4}$ due to arithmetic effects \cite[Lemma~7.1]{Hokken}. In particular, the bound \eqref{eq:square_discriminant} cannot be sharpened in general by more than a factor of $\sqrt{\log{m}}$. Moreover, a lower bound on the probability that $\Delta(A)$ is a non-zero square must take arithmetic effects into account.
\end{remark}

\section{From Propositions~\ref{prop:factorisation options}--\ref{prop:square_discriminant} to Theorems~\ref{thm:simple_interval}--\ref{thm:master}}
\label{sec:proof_of_main_theorems}

In this section, we deduce Theorems~\ref{thm:simple_interval}--\ref{thm:master} from Propositions~\ref{prop:factorisation options}--\ref{prop:square_discriminant}. We commence by proving Theorem~\ref{thm:master}, which we then show implies the other main theorems. The proofs of Propositions~\ref{prop:factorisation options}--\ref{prop:square_discriminant} will take up the remaining sections of this paper. 

\subsection*{Proof of Theorem~\ref{thm:master}}

Define the events
\begin{align*}
\CE_1 &\coloneqq \{A \textup{ has a reciprocal divisor } D \in \CR(k) \textup{ with } 1 \le k < m^{1/10}\}, \\ 
\CE_2 &\coloneqq \{ A \textup{ has a reciprocal divisor } D \in \CR(k) \textup{ with } k \in [m^{1/10}, m/2]\}, \\ 
\CE_3 &\coloneqq \{ \exists \textup{ monic, irreducible, nonreciprocal } I \in \Z[T] : A = \pm I\cdot I_{\ms{rev}} \}.
\end{align*}
To show \eqref{eq:A is irreducible}, by Lemma~\ref{lem:factorisation options} it suffices to show that $\CE_1$, $\CE_2$ and $\CE_3$ occur with very small probability. Possibly after adjusting $C$, we may assume $m$ is sufficiently large throughout. As in the proof of \cite[Theorem~7]{BKK}, we have the following: for any $j$, any divisor $Q \mid P$ with $Q>1$, and any $b \in \Z/Q\Z$, Fourier inversion and condition (2) imply
\begin{equation*}
\sum_{a \equiv b \bmod{Q}} \mu_j(a) \le \frac{1}{Q} \sum_{k = 0}^{Q-1} \abs{\hat{\mu}_j(k/Q)} \leq \frac{1-m^{-1/10}}{\sqrt{Q}} < \frac{1}{\sqrt{2}}.
\end{equation*}
In particular, for any $j$ and any prime divisor $p \mid P$ we have
\begin{equation}
\label{eq:anti-concentration}
\norm{\mu_j}_{\infty} \leq \norm{\bar{\mu}_{j,p}}_{\infty} < \frac{1}{\sqrt{2}}.
\end{equation}

All this implies the following. For $\CE_1$, condition (1) of Theorem~\ref{thm:master} implies condition (1) of Proposition~\ref{prop:small degree factors} as $m$ can be taken sufficiently large. Condition (2) of Proposition~\ref{prop:small degree factors} holds by \eqref{eq:anti-concentration}. Thus there exists an absolute constant $C_1 > 0$ such that $\mb{P}(\CE_1) \leq C_1 m^{-7/20}$.
For $\CE_2$, we may first apply Proposition~\ref{prop:Delta bound} with $\CP$ the set of prime divisors of $P$ and $\gamma = 1/2$ to find that condition (1) of Proposition~\ref{prop:large degree factors} holds for any $\epsilon_{\mr{Proposition~\ref{prop:large degree factors}}} \in (0, 1/100]$. Pick $\epsilon_{\mr{Proposition~\ref{prop:large degree factors}}} = 1/100$. Condition (2) of Proposition~\ref{prop:large degree factors} is also satisfied by \eqref{eq:anti-concentration}. Thus there are absolute constants $c_2, C_2 > 0$ for which $\mb{P}(\CE_2) \leq C_2 m^{-c_2}$.
For $\CE_3$, we again employ \eqref{eq:anti-concentration} to apply Proposition~\ref{prop:exceptional} with $\epsilon_{\mr{Proposition~\ref{prop:exceptional}}} = 1-1/\sqrt{2}$. By condition (1) of Theorem~\ref{thm:master}, this implies there exists $C_3 = C_3(\alpha, B)>0$ such that $\mb{P}(\CE_3) \leq C_3 m^{-\alpha/2}$.
We conclude that
\begin{equation*}
\mb{P}_{A \in \CR(m)}(A \textup{ is reducible}) \le \mb{P}(\CE_1) + \mb{P}(\CE_2) + \mb{P}(\CE_3) \le Cm^{-c}
\end{equation*}
for some constants $c, C>0$ depending at most on $\alpha$ and $B$.

To go from irreducibility to the Galois group, we require Proposition~\ref{prop:Galois}. Its condition (1) holds by condition (1) of Theorem~\ref{thm:master}, its condition (2) by condition (2) of Theorem~\ref{thm:master} combined with Proposition~\ref{prop:Delta bound} and the choice $\epsilon_{\mr{Proposition~\ref{prop:Galois}}} = 1/100$, and its condition (3) by \eqref{eq:anti-concentration}. Thus, combining the conclusion of Proposition~\ref{prop:Galois} with the already proved \eqref{eq:A is irreducible}, the result follows. \qed

\begin{remark}
\label{rem:thmCthm3}
(a) Our Theorem~\ref{thm:master} is quite similar to Theorem 7 in \cite{BKK}. However, observe that our condition on the support of the measures $\mu_j$ is stricter, which is needed to preclude the exceptional factorisation (in Proposition~\ref{prop:exceptional}) and show that the discriminant is unlikely to be a square (in Proposition~\ref{prop:square_discriminant}). Furthermore, in the notation of \cite[Theorem~7]{BKK}, we do not cover the case $s>1$, for technical reasons explained in Remark~\ref{rem:why not s>1}.

\medskip

(b) The conclusion of Theorem 7 in \cite{BKK} contains a typo: it should say that $A$ is unlikely to have divisors of small degree, not that it is unlikely not to have divisors of small degree.
\end{remark}

\subsection*{Proof of Theorem~\ref{thm:simple_interval}}
We will deduce this from Theorem~\ref{thm:master}. As $H = \max_{a \in \CN} \abs{a}$ is bounded, condition (1) of Theorem~\ref{thm:master} is satisfied for any $\alpha \in (0, 1/2)$. To prove that the Fourier condition (2) of Theorem~\ref{thm:master} is satisfied when all measures $\mu_j$ equal the fixed, uniform measure $\mu$ on a set $\CN$ of at least $35$ consecutive integers, we refer to \cite[\S 3.1]{BKK}. \qed

\subsection*{Proof of Theorem~\ref{thm:not_large_not_sparse}}
Again, we deduce this from Theorem~\ref{thm:master}. Possibly after adjusting $C$, we may assume that $m$ is sufficiently large. Since $\mu_j = \mu$ is fixed, this also implies that we may assume that condition (1) of Theorem~\ref{thm:master} is satisfied. To prove that condition (2) of Theorem~\ref{thm:master} holds as well, we refer to \cite[\S 3.6]{BKK}. \qed

\section{Prime factorization of reciprocal polynomials}
\label{sec:divisors of reciprocal polynomials}

In this section, we prove Proposition \ref{prop:factorisation options}. We begin with the following simple lemma.

\begin{lemma}\label{lem:irreducible vs reciprocity}
	Let $K$ be a field and let $I$ be a monic irreducible polynomial over $K$. 
			\begin{enumerate}
				\item[(a)] If $I=I_{\ms{rev}}$, then either $I$ is reciprocal or $I(T)=T+1$.
				\item[(b)] If $I\neq I_{\ms{rev}}$, then either $I$ and $I_{\ms{rev}}$ are coprime, or $I(T)=T-1$ and $\car(K)\neq2$.
			\end{enumerate} 
\end{lemma}

\begin{proof}
We start with (a). The only monic polynomial $I$ of degree $1$ satisfying the relation $I=I_{\ms{rev}}$ is the polynomial $I(T)=T+1$. Let us assume now that $\deg{I}\ge2$. Since $I=I_{\ms{rev}}$, we will have that $I$ is reciprocal as long as we can show that $\deg{I}$ is even. Indeed, according to Remark \ref{rem:odd degree reciprocal}, if $\deg{I}$ is odd, then the relation $I=I_{\ms{rev}}$ implies that $I$ can be factored as $I(T)=(T+1)J(T)$ for some reciprocal polynomial $J$. This is impossible because we have assumed that $I$ is irreducible and has degree $\ge2$. We have thus established our claim that $\deg{I}$ is even. This completes the proof of part (a) of the lemma.

For (b), let us assume that $I$ and $I_{\ms{rev}}$ are not coprime; otherwise there is nothing to prove. Since both of these polynomials are irreducible, we must have that $I_{\ms{rev}}=uI$ for some $u\in K^\times$. Since $I$ is monic, we have $I(T)=T^m+c_{m-1}T^{m-1}+\cdots+c_0$ for some $c_j\in K$. The relation $I_{\ms{rev}}=uI$ then implies that $c_0=u$ and that $1=uc_0$. We thus find that $u=c_0=\pm1$. Since we have assumed that $I \neq I_{\ms{rev}}$, we must have $u=-1\neq 1$. In particular, $K$ has characteristic different than $2$, and thus the relation $I=-I_{\ms{rev}}$ implies that $I(1)=0$. Since $I$ is monic irreducible, we must have $I(T)=T-1$. This completes the proof of the lemma.
\end{proof}

\begin{proof}[Proof of Proposition~\ref{prop:factorisation options}]
Let $J$ be an irreducible factor of $A$ such that $J\neq J_{\ms{rev}}$. The polynomial $J_{\ms{rev}}$ is also irreducible. Since $J$ divides $A$, the polynomial $J_{\ms{rev}}$ divides $A_{\ms{rev}}=A$ as well. Moreover, since $A$ is reciprocal, $J$ and $J_{\ms{rev}}$ must divide $A$ to the exact same power. If $J$ and $J_{\ms{rev}}$ are coprime, we thus find that there exists an integer $\nu\ge0$ such that 
\[
A=(JJ_{\ms{rev}})^\nu B.
\]
The polynomial $B$ must be reciprocal because $A$ and $(JJ_{\ms{rev}})^\nu$ are reciprocal. 
These observations together with Lemma \ref{lem:irreducible vs reciprocity} readily imply that we may write 
\begin{equation}
	\label{eq:A factorization}
	A(T)= u \cdot (T-1)^k (T+1)^\ell \prod_{i=1}^r I_i(T) \prod_{j=1}^s J_j(T)J_{j,\ms{rev}}(T), 
\end{equation}
where each $I_i$ is reciprocal irreducible, each $J_j$ is irreducible and coprime to $J_{j,\ms{rev}}$, and $k,\ell\in\Z_{\ge0}$ and $u\in K^\times$. It remains to check that $k$ and $\ell$ are even. For the former claim, we may assume that $\car(K)\neq2$ (as otherwise $T-1=T+1$ and we may take $k=0$). If $B(T)$ denotes the polynomial on the right-hand side of \eqref{eq:A factorization}, then $B_{\ms{rev}}=(-1)^k B$. Since $B=A$ is reciprocal, and we also know here that $\car(K)\neq2$, we conclude that $k$ must be even. 
Finally, we must have that $\ell$ is even because $A$ and the $I_i$ are all reciprocal polynomials and thus have even degrees (cf.~Definition \ref{def:reciprocal}).
\end{proof}

\begin{proof}[Proof of Lemma~\ref{lem:factorisation options}] 
Since $A$ is monic and has integer coefficients here, the polynomials $I_i$ and $J_j$ in Proposition \ref{prop:factorisation options} can be assumed to be monic and to lie in $\Z[T]$. Comparing leading coefficients, and using again our assumption that $A$ is monic, we find that the $u$ of Proposition~\ref{prop:factorisation options} satisfies $u=\pm1$. 

If $a+b+r+s\ge2$, the result follows by selecting a factor of minimal degree among
\begin{equation*}
		\underbrace{(T-1)^2,\dots,(T-1)^2}_{a\ \text{times}}, \quad \underbrace{(T+1)^2,\dots,(T+1)^2}_{b\ \text{times}}, \quad I_1,\dots,I_r, \quad \text{and} \quad J_1J_{1,\ms{rev}}, \dots, J_sJ_{s,\ms{rev}}. 
\end{equation*}
	
Let us now examine the case $a+b+r+s=1$. The subcase when  $r=1$ and $a=b=s=0$ can be discarded, because $A$ would then be irreducible. In each of the remaining subcases, we may easily check that we may write $A=uI\cdot I_{\ms{rev}}$ for some monic irreducible polynomial $I$: we take $I=T-1$, $I=T+1$ or $I=J_1$, according to whether $a=1$, $b=1$ or $s=1$. This completes the proof.
\end{proof}

We will also make use of the following notion --- a reciprocal version of the greatest common divisor.

\begin{definition}
\label{def:reciprocal gcd}
For two (not-necessarily reciprocal) polynomials $A,B \in K[T]$, not both equal to $0$, let $(A,B)_{\Br}$ denote the \emph{greatest common reciprocal divisor} of $A$ and $B$: that is, the monic reciprocal polynomial of maximal degree dividing both $A$ and $B$. 
\end{definition}

\begin{remarks*}
(a) The polynomial $(A,B)_{\Br}$ is unique by the requirement that it is monic, just as the usual greatest common divisor $(A,B)$ is unique if it is required to be monic. 

\medskip

(b) Using Proposition \ref{prop:factorisation options}, we find that $(A,B)_{\Br} = (A,B)$ if $A$ and $B$ are reciprocal. However, this relation is not true in general.
\end{remarks*}

\section{The trace polynomial}
\label{sec:the_trace_polynomial}

Let $K$ be a field. In this preliminary section, we define the \emph{trace polynomial} $A_\SR$ of a reciprocal polynomial $A \in K[T]$. When viewing $A$ as an element of a particular vector space defined below, the correspondence $A \mapsto A_{\SR}$ defines a linear bijection that is also multiplicative. Under this map, the divisors of $A_\SR$ correspond to the \emph{reciprocal} divisors of $A$. We will use this correspondence in several places, in particular when studying the \emph{anatomy} of reciprocal polynomials in \S \ref{sec:anatomy}. 

We start by observing that the set of polynomials of degree $\le \ell$ in $K[T]$ is a vector space over $K$, whereas the set of reciprocal polynomials of degree $\le 2\ell$ is not: it is not closed under addition. For this reason, we introduce the following notion:

\begin{definition}[Shifted reciprocal polynomials]\label{def:ell-shifted reciprocal}
	Let $K$ be a field. 
		\begin{enumerate}[label=(\alph*)]
			\item We write $\CR^{\ms{sh}}_K(\ell)$ for the set of \emph{$\ell$-shifted reciprocal polynomials}, which are defined to be all polynomials of the form
			\[
			C(T) = c_\ell + c_{\ell-1}T + \cdots + c_1 T^{\ell-1}+ c_0 T^{\ell} + c_1 T^{\ell+1} + \cdots + c_{\ell-1}T^{2\ell-1} + c_\ell T^{2\ell} \in K[T].
			\]
			\item We write $\CR^{\ms{sh}}_K=\bigcup_{\ell\ge0} \CR^{\ms{sh}}_K(\ell)$ for the set of \emph{shifted reciprocal polynomials}.
		\end{enumerate}
\end{definition}

\begin{remarks*}
	(a) If $C(T)\in \CR^{\ms{sh}}_K$ is nonzero, then there exists a unique integer $\ell\ge0$ such that $C(T)\in \CR^{\ms{sh}}_K(\ell)$. 
	
	\medskip
	
	(b) $\CR^{\ms{sh}}_K(\ell)$ is a $K$-vector space of dimension $\ell+1$. A basis for it is given by
	\begin{equation}
		\label{eq:reciprocal basis}
		\{T^{\ell}, T^{\ell-1} + T^{\ell+1}, \ldots, T + T^{2\ell-1}, 1+T^{2\ell}\}.
	\end{equation}
	
	(c)  For any integers $0 \le i \le j \le \ell$ we have the inclusions
	\begin{equation*}
		T^{\ell-i} \CR_K(i) \subset T^{\ell-j} \CR_K^{\ms{sh}}(j) \subset \CR_K^{\ms{sh}}(\ell) ,
	\end{equation*}
where the last inclusion is one of vector spaces. 
\end{remarks*}

Let $K[T]_{\ell}$ denote the set of monic polynomials in $K[T]$ of degree $\ell$ and write $K[T]_{\le \ell}$ for the set of (not-necessarily monic) polynomials in $K[T]$ of degree $\le \ell$, the latter including the zero polynomial. Also denote by $\CR_K$ the set of all reciprocal polynomials over $K$. We construct a correspondence between $K[T]_{\le\ell}$ and $\CR_K^{\ms{sh}}(\ell)$. 

\begin{definition}[Reciprocal maps]\label{def:reciprocal maps}
	Let $K$ be a field. 
	\begin{enumerate}[label=(\alph*)]
		\item For each integer $\ell\ge0$, we define the map $(\boldsymbol{-})^{\SR, \ell} \colon K[T]_{\le \ell} \to \CR^{\ms{sh}}_K(\ell)$ by
		\[
			G(T) \mapsto G^{\SR, \ell}(T) \coloneqq T^{\ell} G(T + T^{-1}) .
		\]
		\item We define the \emph{reciprocal mapping} $(\boldsymbol{-})^{\SR} \colon K[T] \to \CR_K$ by
		\[
		G(T) \mapsto G^{\SR}(T) \coloneqq 
			\begin{cases} G^{\SR, \, \deg{G}}(T) = T^{\deg{G}} G(T + T^{-1})&\text{if}\ G\neq0,\\
					0 &\text{otherwise}.
			\end{cases}
		\]
	\end{enumerate}
\end{definition}

\begin{lemma}
	\label{lem:trace pol iso of vs}
	The map $(\boldsymbol{-})^{\SR, \ell}$
	is an isomorphism of $K$-vector spaces. Moreover, for any $j\le \ell$, it restricts to an isomorphism between the $K$-vector spaces $K[T]_{\le j}$ and $T^{\ell-j}\CR_K^{\ms{sh}}(j)$, and to a bijection between the sets $K[T]_j$ and $T^{\ell-j}\CR_K(j)$.
\end{lemma}
\begin{proof}
The map $(\boldsymbol{-})^{\SR, \ell}$ is a linear isomorphism because it is an injective, linear homomorphism between vector spaces of the same dimension. The claims made in the last sentence of the lemma statement are immediate consequences of the statement made in the first sentence, so we are done.
\end{proof}

\begin{definition}
\label{def:inverse general reciprocal map}
The inverse of $(\boldsymbol{-})^{\SR, \ell}$ is denoted by $(\boldsymbol{-})_{\SR, \ell}$.
\end{definition}

Next, we define the trace polynomial of a reciprocal polynomial.

\begin{definition}[The trace polynomial]\label{def:trace polynomial} Let $K$ be a field. Let $(\boldsymbol{-})_{\SR}$ denote the inverse function of $(\boldsymbol{-})^{\SR}$. Given $A\in \CR_K$, we call $A_{\SR}$ the \emph{trace polynomial} of $A$. 
\end{definition}

From Definitions~\ref{def:reciprocal maps} and \ref{def:trace polynomial}, it follows immediately that the reciprocal mapping and its inverse commute with multiplication, that is $(FG)^{\SR} = F^{\SR} G^{\SR}$ and $(AB)_{\SR} = A_{\SR} B_{\SR}$. 

The trace polynomials of the polynomials of the form $T^{2j}+1$ are the \emph{Chebyshev polynomials}.

\begin{definition}[Chebyshev polynomials] \label{def:chebyshev polynomials}Let $K$ be a field and $j\ge0$. We then define the $j$-th \emph{Chebyshev polynomial} $C_j(T) \in K[T]$ through the recursion
\[
		C_0(T) \coloneqq 2, \quad C_1(T) \coloneqq T, \quad \text{and} \quad C_{j+1}(T) \coloneqq TC_j(T) - C_{j-1}(T) \, \text{ for } j \ge 1.
\]
We extend this definition to negative indices by setting
\[
C_{-j}\coloneqq C_j.
\]
\end{definition}

\begin{remarks*}
	(a) For reference, when the characteristic of $K$ is not $2$, our Chebyshev polynomials are a common rescaling of the standard Chebyshev polynomials $\mr{T}_j$. In characteristic $0$, these are defined as the unique polynomial satisfying $\mr{T}_j(\cos{x}) = \cos(jx)$; they are related through $C_j(2T) = 2\mr{T}_j(T)$. 
	
	\medskip
	
	(b) A simple inductive argument yields
	\begin{equation*}
		C_j(T+T^{-1}) = T^j + T^{-j}\quad\text{for}\ j=0,1,2,\dots
	\end{equation*}
	In other words, we indeed have $C_j^{\SR} = T^{2j}+1$ for nonnegative $j$.

	\medskip 

	(c) The set $\{1, C_1(T), C_2(T), \ldots, C_{\ell}(T)\}$ is a basis of $K[T]_{\le \ell}$ and maps to the basis \eqref{eq:reciprocal basis} of $\CR^{\ms{sh}}_K(\ell)$ under $(\boldsymbol{-})^{\SR, \ell}$.
	
	\medskip
	
	(d) We may easily check that
	\begin{equation}
		\label{eq:chebyshev recursion general}
		TC_j(T) = C_{j+1}(T) + C_{j-1}(T) \, \text{ for all } j\in\Z.
	\end{equation} 
\end{remarks*}

One can spell out explicitly how the coefficients of $A_{\SR}$ depend on those of $A$:

\begin{lemma}
\label{lem:trace polynomial coefficients}
Given $A(T) = 1 + a_{m-1} T + \cdots + a_0T^m + \cdots + a_{m-1}T^{2m-1} + T^{2m} \in \CR_K(m)$ with trace polynomial $A_{\SR}(T) = b_0  + \cdots + b_{m-1} T^{m-1} + T^m$, the coefficients of $A$ and $A_{\SR}$ satisfy the relations
\begin{equation} 
	\label{eq:trace polynomial coefficients}
	b_{i} = a_{i} + \sum_{j = 1}^{\lfloor \frac{m-i}{2} \rfloor} (-1)^j \left(\binom{i+j}{j} + \binom{i+j-1}{j-1} \right) a_{i+2j}
\end{equation}
and
\begin{equation*}
a_i = \sum_{j=0}^{\lfloor \frac{m-i}{2} \rfloor} b_{i+2j} \binom{i+2j}{i+j}.
\end{equation*}
\end{lemma}
The proof is elementary but lengthy; we omit it here.

A consequence of the relation between the coefficients of $A$ and $A_{\SR}$, described by Lemma~\ref{lem:trace polynomial coefficients}, is that the probability measure $\mb{P}_{\CR(m)}$ on the set of monic reciprocal polynomials in $\Z[T]$ of degree $2m$ induces a probability measure $\mb{P}_{\SR,\CM(m)}$ on the set $\CM(m)$ of monic polynomials in $\Z[T]$ of degree $m$:
\begin{definition}
\label{def:induced_measure_on_trace_pols}
If $B \in \Z[T]$ is a monic polynomial of degree $m$, then we define
\begin{equation*}
\label{eq:induced_measure_on_trace_pols}
\mb{P}_{\SR,\CM(m)}(B) \coloneqq \mb{P}_{\CR(m)}\big(B^{\SR}\big).
\end{equation*}
If $P$ is some property, then we will further use the notation
\[
\mb{P}_{\SR,B\in\CM(m)}\Big(B\ \text{has property}\ P\Big) \coloneqq  \mb{P}_{A\in \CR(m)}\Big(A_\SR\ \text{has property}\ P\Big).
\]
\end{definition}
By \eqref{eq:trace polynomial coefficients}, the coefficients of $B$ are generally speaking dependent random variables, so we cannot apply the results of \cite{BKK} directly to show that the trace polynomial of $A$ is irreducible with high probability.

\begin{definition}[Semi-irreducibility]
\label{def:semi-irreducibility} 
Let $K$ be a field and let $A \in \CR_K$ be nonconstant. We say that $A$ is \emph{semi-irreducible} if it has no nontrivial reciprocal factors (of even degree, cf. Definition~\ref{def:reciprocal}).
\end{definition}

\begin{lemma}	\label{lem:semi-irreducibility}
	Let $K$ be a field and let $A\in\CR_K$ be nonconstant. The following are equivalent:
	\begin{enumerate}[label=(\alph*)]
		\item $A$ is semi-irreducible;
		\item $A_{\SR}$ is irreducible;
		\item either $A$ is irreducible, or $A=uI\cdot I_{\ms{rev}}$ for a nonreciprocal irreducible polynomial $I$ and for some $u\in K^\times$.
	\end{enumerate}
\end{lemma}

\begin{proof} It is easy to see that (a) and (b) are equivalent: if we have a factorization of the trace polynomial $A_\SR=BC$, then we obtain a factorization $A=B^\SR C^\SR$ into two reciprocal factors of even degree. The converse is also true.
	
Let us now prove that (a) and (c) are equivalent. Clearly, (c) implies (a). Hence, it remains to show that (a) implies (c). Recall Proposition \ref{prop:factorisation options} and the notation therein. Since $A$ has no nontrivial reciprocal factors, we must have that $a+b+r+s = 1$. If $a=1$, we can take $I = T-1$. If $b=1$, we can take $I=T+1$. If $r=1$, then $A$ is irreducible. Lastly, if $s=1$, we can take $I = J_1$. This completes the proof.
\end{proof}

\begin{lemma}\label{lem:coprimality in terms of trace polynomials}
	Let $K$ be a field and let $A,B\in K[T]\setminus\{0\}$. Then $A$ and $B$ are coprime if, and only if, $A^{\SR}$ and $B^{\SR}$ are coprime.
\end{lemma}

\begin{proof}
	Assume that $A$ and $B$ are coprime, and let $D=\gcd(A^{\SR},B^{\SR})$. Using Proposition \ref{prop:factorisation options}, we find that $D$ is reciprocal. In particular $D_{\SR}$ is a common factor of $A$ and $B$, and hence $D_{\SR}=1$. We thus also find that $D=1$.
	
	Conversely, if $D$ is a nontrivial common factor of $A$ and $B$, then $D^{\SR}$ is a nontrivial common factor of $A^{\SR}$ and $B^{\SR}$. This completes the proof.
\end{proof}

\section{Small degree factors}
\label{sec:small degree factors}

In this section, we prove Proposition~\ref{prop:small degree factors}, which states that $A$ is unlikely to have reciprocal factors of degree $\le 2m^{1/10}$. The proof is a direct analogue of \cite[Proposition~2.1]{BKK}; we will thus skip some details.
Both proofs follow the work of Konyagin \cite{Konyagin}, who considered polynomials with coefficients in $\{0,1\}$. He presented two arguments, of which we follow the first and simpler one. The argument proceeds in two steps, bounding cyclotomic and non-cyclotomic divisors separately.

As in the other results in \cite{BKK} and in the present work, we must assume that the probability measures are contained in some bounded interval $[-H,H]$, but here we may take the comparatively large $H \coloneqq \lfloor \exp(m^{1/3}) \rfloor$. If $A$ has a reciprocal factor of degree $\leq 2k_0$, then it must have a semi-irreducible factor of degree $\le 2k_0$ (see Definition~\ref{def:semi-irreducibility}). Hence it suffices to control the probability that the latter event occurs with $k_0=m^{1/10}$. 

\begin{definition}
We let $\CD(k_0)$ denote the set of monic semi-irreducible reciprocal polynomials $1+d_1T + \cdots + d_1T^{2k-1} + T^{2k}$  of degree $\le 2k_0$ for which each $d_j$ is bounded by $(H+1)^{2k}$ in absolute value.
\end{definition}
An analogous argument to the one leading up to \cite[Equation~7.2]{BKK} implies that if $A$ has a factor of degree $\le 2k_0$, then it must have a factor in $\CD(k_0)$.

We start by bounding the probability that $A$ has a cyclotomic divisor of small degree. Here and elsewhere in the paper we will use the following inequality, due to Kolmogorov and Rogozin \cite[Theorem~2]{Rogozin}. 
\begin{lemma}[Kolmogorov--Rogozin inequality]
\label{lem:Kolmogorov--Rogozin}
Suppose $X_1, \ldots, X_\ell$ are discrete, independent random variables. Write $X = X_1 + X_2 + \cdots + X_\ell$. Then 
\begin{equation}
\label{eq:Kolmogorov--Rogozin}
\sup_a \mb{P}(X=a) \ll \Big( \sum_{j=1}^\ell \big(1-\sup_b \mb{P}(X_j=b)  \big) \Big)^{-1/2}.
\end{equation}
In particular, if there is an $\epsilon>0$ such that $\sup_{b} \mb{P}(X_j=b) \leq 1-\epsilon$ for all $j$, then
\begin{equation}
\label{eq:Kolmogorov--Rogozin_discrete}
\sup_{a} \mb{P}(X = a) \ll \frac{1}{\sqrt{\epsilon \ell}}.
\end{equation}
\end{lemma}
\begin{remark*}
When $\sup_{b} \mb{P}(X_j=b) = 1$ for all $j$, the right-hand side of \eqref{eq:Kolmogorov--Rogozin} is understood as $1$.	
\end{remark*}

\begin{lemma}
\label{lem:small_cyclotomic}
Fix $\epsilon>0$. For each $j \in [0,m-1]\cap\Z$, let $\mu_j$ be a probability measure such that $\norm{\mu_j}_{\infty} \le 1-\epsilon$. 
Then for fixed $z \in \C$, we have
\begin{equation*}
\mb{P}_{A \in \CR(m)} (A(z) = 0) \ll \frac{1}{\sqrt{\epsilon m}}.
\end{equation*}
In particular, given some positive integer $k_0$, we find
\begin{equation*}
\sum_{d \,:\, \phi(d) \le 2k_0} \mb{P}_{A \in \CR(m)}(\Phi_d \mid A) \ll \frac{k_0}{\sqrt{\epsilon m}}.
\end{equation*}
\end{lemma}
\begin{proof}
We follow \cite[Lemma~7.3]{BKK}.
Note that $A(0) \neq 0$ for a reciprocal polynomial, so we may assume $z \neq 0$. For integers $j$ in the interval $[1,m-1]$, define the random variable $X_j = a_{j} (z^{m-j} + z^{m+j})$. Set $X_0 = a_0 z^m$. Then $X_0, X_1, \ldots, X_{m-1}$ is a sequence of discrete random variables. Let $J$ denote the set of values of $j \in [0,m-1] \cap \Z$ for which $z^{m-j} + z^{m+j}$ is nonzero. Since $z$ is nonzero, the cardinality of $J$ is strictly smaller than $m$ only if $z$ is a root of unity, say $z = e(\theta)$ for some $\theta \in \Q$. In that case, the relation $z^{m-j} + z^{m+j} = 0$ is equivalent to $\cos(2\pi j\theta) = 0$, i.e. to $j\theta \in \{1/4, 3/4\} \bmod{1}$. In particular, the denominator of $\theta$ must be a multiple of $4$, and $J$ must have cardinality at least $m/2$ (with equality if $\theta \in \{1/4, 3/4\}$ and $m$ is even). 

Now, let $X = \sum_{j\in J} X_j$. Then
\begin{align*}
\mb{P}\big(A(z) = 0\big) = \mb{P}\big( X = -(1+z^{2m})\big) 
	&\le \sup_{u \in \supp(X)} \mb{P}(X = u) \\
	&\ll \Big( \sum_{j \in J} \big(1-\sup_{b \in \supp{X_j}}\mb{P}(X_j = b)\big) \Big)^{-1/2},
\end{align*}
by Lemma~\ref{lem:Kolmogorov--Rogozin}. Hence
\begin{equation*}
\mb{P}(A(z) = 0) \ll \frac{1}{\sqrt{\epsilon \cdot \#J}} \ll \frac{1}{\sqrt{\epsilon m}}.
\end{equation*}
Since $\Phi_d$ is irreducible for all $d$, we have $\Phi_d \mid A$ if and only if $A(e(1/d))=0$. There are $\ll k_0$ values of $d$ with $\phi(d) \leq 2k_0$ by \cite{Bateman}, which concludes the proof. 
\end{proof}

For the noncyclotomic divisors, we proceed by means of the \emph{Mahler measure}: if $F \in \Z[T]$ is a polynomial of degree $d$ with leading coefficient $f_0$ and zeros $z_1, \ldots, z_{d}$ (listed with multiplicity), then the (exponential) Mahler measure of $F$ is the quantity
\begin{equation*}
M(F) \coloneqq \abs{f_0} \prod_{j=1}^{d} \max\big\{1, |z_j| \big\}.
\end{equation*}
Choose $C>0$ such that $\log(x)/\log{\log{x}}$ is strictly increasing and larger than $1$ on $x>C$. Let $c \ge 1200$ be minimal such that $M(F) \ge 1+1/c$ for all irreducible, noncyclotomic polynomials $F$ with $F(0) \neq 0$ of degree $\le C$ and with leading coefficient $\pm 1$; it is a classical result that such $c$ exists, see e.g. \cite[Lemma~1.3]{BZ}. Now define 
\begin{equation*}
L(x) = 
\begin{cases}
c & \text{if } x \le C, \\ 
c(\log(x)/\log{\log{x}})^3 &\text{if } x>C,
\end{cases}
\end{equation*}
and observe that $L$ is a nondecreasing function on $x>0$. In particular, the work of Dobrowolski \cite{Dobrowolski} implies that 
\begin{equation}
	\label{eq:Dobrowolski}
	M(D) \ge 1 + 1/L(2k)
\end{equation}
for all semi-irreducible reciprocal polynomials $D$ of degree $2k$ that are noncyclotomic and different from $(T\pm 1)^2$. Indeed, if $D$ is irreducible, this follows directly by \cite{Dobrowolski}. On the other hand, if $D = \pm I\cdot I_{\ms{rev}}$, then we have $M(D) = M(I)M(I_{\ms{rev}}) \ge (1+1/L(k))^2 > 1 + 1/L(2k)$ since $L$ is nondecreasing.

Fix a noncyclotomic polynomial $D \in \CD(k_0)$ of degree $2k$ other than $(T+1)^2$ and $(T-1)^2$. Let $J=D$ if $D$ is irreducible, and let $J=I$ if
$D=\pm I I_{\ms{rev}}$, where we may assume $I$ to be monic. Set $d=\deg J$, so that $d=2k$ in the first case and $d=k$ in the second. Then $J$ is monic, irreducible and noncyclotomic. Let $z_1,\ldots,z_d$ be its distinct roots. Following \cite{BKK}, pick a prime number
\begin{equation*}
	p = p_D \in ((1+L(2k))\log(4Hm), 2(1+L(2k))\log(4Hm)]
\end{equation*}
for which $z_1^p$ has algebraic degree $d$. That such $p$ exists follows from \cite[Lemma~3]{Dobrowolski}. The conjugates of $z_1^p$ are $z_2^p,\ldots,z_d^p$, which must thus be distinct. Moreover, we have
\[
M(J)\ge 1+\frac{1}{L(d)}
\ge 1+\frac{1}{L(2k)}.
\]

\begin{lemma}
\label{lem:small_noncyclotomic}
Let $m,k\in\Z_{\ge1}$. Let $D \in \CD(k_0)$ be distinct from $(T+1)^2$ and $(T-1)^2$ and noncyclotomic, and let $p=p_D$ be as above. Consider integers $c_j$ for each $j \in [0,m]\cap\Z$ with $j\not\equiv 0\pmod p$. Then there exists at most one reciprocal polynomial $A(T) = T^{2m}+a_{m-1}T^{2m-1}+\cdots +a_0T^m+\cdots + a_{m-1}T + 1$ such that $D|A$, $|a_j|\le H$ for all $j$, and $a_j = c_j$ for all $j\not\equiv 0\pmod p$.
\end{lemma}

\begin{proof}
Assume $A$ and $B$ are two distinct polynomials as in the lemma statement. The $2m$ coefficients of $A-B$ are bounded by $2H$ in absolute value, and the Mahler measure of $A-B$ is bounded above by the sum of the absolute values of these coefficients \cite[Lemma~1.7]{BZ}, so that $M(A-B) < 4Hm$. However, we will prove that $M(A-B) \ge M(J)^p$ (with $J$ as above), which is $> 4Hm$ by definition of $p$ and by \eqref{eq:Dobrowolski}. This yields a contradiction.

To prove the claim that $M(A-B) \ge M(J)^p$, observe that $D$ divides
\begin{equation*}
A(T)-B(T) = T^m \bigg( g_0 + \sum_{0 < j < m/p} g_j (T^{-pj} + T^{pj}) \bigg).
\end{equation*}
The Mahler measure of $A-B$ is the same as that of $G(T) = T^{(p-1)m}(A(T)-B(T))$, which is still divisible by $D$ and hence by $J$, and is a polynomial in $T^p$. Let $\zeta$ be a primitive $p$-th root of unity. Observe that there are no integers $i \neq j$ and $a$ such that $\zeta^a z_i = z_j$, as then we would have $z_i^p = z_j^p$. 
It follows that $G$ is divisible by $\prod_{\ell=1}^{d} \prod_{j=0}^{p-1} (T-\zeta^j z_\ell)$. 
Hence
\[
M(A-B) \ge \prod_{\ell=1}^{d} \prod_{j=0}^{p-1} \max\big\{1, |\zeta^j z_\ell | \big\} = \prod_{\ell=1}^{d} \max\big\{1, |z_\ell| \big\}^p = M(J)^p. \qedhere
\]
\end{proof}

This brings us to the proof of Proposition~\ref{prop:small degree factors}, which is again very similar to the proof of \cite[Proposition~2.1]{BKK}.

\begin{proof}[Proof of Proposition~\ref{prop:small degree factors}]
Set $k_0 \coloneqq \lfloor m^{1/10} \rfloor$. By Lemma~\ref{lem:small_cyclotomic}, which we may apply with $\epsilon = m^{-1/10}$, the probability that $\Phi_d \mid A$ for some cyclotomic $\Phi_d \in \CD(k_0)$ is $\ll m^{1/10}/\sqrt{m^{9/10}} = m^{-7/20}$. The probability that $(T+1)^2$ or $(T-1)^2$ divides $A$ equals the probability that $A(1)$ or $A(-1)$ vanishes, which is $\ll m^{-9/20}$, again by invoking Lemma~\ref{lem:small_cyclotomic}.

It remains to bound the probability of $A$ having a noncyclotomic divisor $D\in\CD(k_0)$ other than $(T \pm 1)^2$. Let $\CE$ denote this event. Then
\[
\mb{P}(\CE) \le \# \CD(k_0) \cdot \sup_{\substack{D \in \CD(k_0) \\ D \text{ noncyclotomic} \\ D \neq (T \pm 1)^2}} \mb{P}(D \mid A).
\]
Observe that all primes $p = p_D$ in the proof of Lemma~\ref{lem:small_noncyclotomic} satisfy $p \le p_{\mr{max}} \le 2(1+L(2m)) \log(4Hm) \ll (\log{m})^3 m^{1/3}$ for some prime $p_{\mr{max}}$. By Lemma~\ref{lem:small_noncyclotomic}, we then find
\begin{equation*}
\sup_{\substack{D \in \CD(k_0) \\ D \text{ noncyclotomic} \\ D \neq (T \pm 1)^2}} \mb{P}(D \mid A) \le \sup_{0 \le j \le m-1} \norm{\mu_j}_{\infty}^{\lfloor m/p_{\mr{max}} \rfloor} \le (1-m^{-1/10})^{\lfloor m/p_{\mr{max}} \rfloor} \le \exp(-m^{0.55}).
\end{equation*}
Lastly, since each coefficient of $D$ is bounded by $(H+1)^{2k_0}$ in absolute value, a straightforward count shows that $\# \CD(k_0) \ll \exp(m^{0.54})$. This completes the proof.
\end{proof}

\section{Exceptional factors and the discriminant}
\label{sec:discriminant_and_exceptional}

Lemma \ref{lem:semi-irreducibility} shows that $A$ may be reducible even if its trace polynomial is not --- that is, even if it has no nontrivial reciprocal divisor. In this case, we have the factorization 
\begin{equation*}
	A = \pm I\cdot I_{\ms{rev}}
\end{equation*} 
for some monic nonreciprocal irreducible polynomial $I \in \Z[T]$ of degree $m$. In this section, we show that this is unlikely to occur by studying the quantity $A(1)A(-1)$. If $A = \pm I\cdot I_{\ms{rev}}$, then $A(1) = \pm I(1)I_{\ms{rev}}(1) = \pm I(1)^2$ and $A(-1) = \pm I(-1) I_{\ms{rev}}(-1) = \pm I(-1)^2$. Thus
\begin{equation}
\label{eq:exceptional_to_A(1)A(-1)}
\mb{P}_{\CR(m)}(A = \pm II_{\ms{rev}} \textup{ for some } I) \le \mb{P}_{\CR(m)}(A(1)A(-1) \neq 0 \textup{ is a square up to sign}).
\end{equation}
Bounding the right-hand side of \eqref{eq:exceptional_to_A(1)A(-1)} is also useful when ruling out subgroups of the alternating group $\CA_{2m}$ as possible Galois group of $A$, as we will do in \S \ref{sec:galois_group_of_A}. Let us first explain this. Recall that if $F \in \Z[T]$ is of degree $d$ and has zeros $\alpha_1, \ldots, \alpha_d$ (listed with multiplicities) and leading coefficient $c$, then the \emph{(polynomial) discriminant} $\Delta(F)$ of $F$ is the integer
\begin{equation}
\label{eq:discriminant}
\Delta(F) \coloneqq c^{2d-2}\prod_{1 \leq i < j \leq d} (\alpha_i - \alpha_j)^2.
\end{equation}
We then have the following basic lemma from Galois theory.
\begin{lemma}
\label{lem:square-disc-alt-gp}
Let $F \in \Z[T]$ be a squarefree polynomial of degree $d$. Then $\Delta(F) \in \Z$ is a square if and only if $\CG_F \leq \CA_{d}$.
\end{lemma}
For general polynomials $F \in \Z[T]$, the quantity $\Delta(F)$ is difficult to analyze. However, reciprocal polynomials have the following neat property.

\begin{lemma}
\label{lem:disc_of_reciprocal}
Let $A \in \CR(m)$. Then $\Delta(A) = (-1)^m A(1)A(-1)\Delta(A_{\SR})^2$. 
\end{lemma}

\begin{proof} 
See e.g. \cite[p.~85]{Dubickas}.
\end{proof}

A consequence of Lemma~\ref{lem:disc_of_reciprocal} is that $(-1)^m A(1) A(-1)$ must be a nonzero square if $\Delta(A)$ is a nonzero square. The converse does not necessarily hold, since $\Delta(A_{\SR})$, and hence $\Delta(A)$, may vanish even if $(-1)^m A(1) A(-1)$ is a nonzero square. Furthermore, the event that $(-1)^m A(1) A(-1)$ is a nonzero square is contained in the event that $A(1)A(-1)$ is a nonzero square up to sign. For convenience, we work with this latter event. In summary, we obtain the inequality
\begin{equation}
\label{eq:square_to_A(1)A(-1)}
\mb{P}_{\CR(m)}\big(\Delta(A) \neq 0 \textup{ is a square}\big) \le \mb{P}_{\CR(m)}\big(A(1)A(-1) \neq 0 \textup{ is a square up to sign}\big),
\end{equation}
which is the same as the right-hand side of \eqref{eq:exceptional_to_A(1)A(-1)}.

Thus our task is to analyze the product $A(1)A(-1)$. In the spirit of \cite{Hokken}, our approach will be to write $A(1)A(-1)$ as the difference of two squares, say $X^2-Y^2$, by splitting the even and odd part of the polynomial $A$. Moreover, the random variables $X$ and $Y$ are both sums of mutually independent random variables. Our main auxiliary result, then, is the following.

\begin{lemma}
\label{lem:X^2-Y^2 is powerful or square}
Fix $c,\epsilon \in (0,1)$. Let $m, m_1, m_2 \in \Z_{\ge 1}$ be positive integers with $m=m_1+m_2$ and $m_1, m_2 \ge cm$. Suppose $X_1, X_2, \ldots, X_{m_1}$ and $Y_1, Y_2, \ldots, Y_{m_2}$ are integer-valued, independent random variables supported in $[-H,H]$. Assume that $\sup_{k \in \Z} \mb{P}(X_j = k) \le 1-\epsilon$ and $\sup_{k \in \Z} \mb{P}(Y_j = k) \le 1-\epsilon$ hold for all $j$.
Write $X = X_1 + X_2 + \cdots + X_{m_1}$ and $Y = Y_1 + Y_2 + \cdots + Y_{m_2}$. 
Then
\begin{equation}
\label{eq:X^2-Y^2 is a square}
\mb{P}(X^2-Y^2 \neq 0 \textup{ is a square}) \ll \frac{H (\log m)^{3/2}}{\epsilon c\sqrt{m}}.
\end{equation}
\end{lemma}

Before proving Lemma~\ref{lem:X^2-Y^2 is powerful or square}, let us see how it implies Proposition~\ref{prop:exceptional} and Proposition~\ref{prop:square_discriminant}.

\begin{proof}[Proof of Proposition~\ref{prop:exceptional} and of Proposition~\ref{prop:square_discriminant}]
By \eqref{eq:square_to_A(1)A(-1)}, and \eqref{eq:exceptional_to_A(1)A(-1)}, it suffices to bound the probability that $A(1)A(-1)$ is nonzero and square up to sign. Write $A(T) = A_e(T^2) + T A_o(T^2)$ for the decomposition of $A$ into its even and odd parts. Then we have $A(1)A(-1) = A_e(1)^2-A_o(1)^2$.
Furthermore, note that $A_e(1) = 2 + 2a_{m-2} + 2a_{m-4} + \cdots$ and $A_o(1) = 2a_{m-1} + 2a_{m-3} + 2a_{m-5} + \cdots$ are both sums of $\ge -1 + m/2$ discrete, independent random variables $V_{j}$ with support contained in $[-3H, 3H]$ and with $\sup_{k \in \Z} \mb{P}(V_j=k) \le 1-\epsilon$ for all $j$. Thus \eqref{eq:square_discriminant} and Proposition~\ref{prop:exceptional} both follow by applying Lemma~\ref{lem:X^2-Y^2 is powerful or square} twice: first with $X = A_e(1)$ and $Y = A_o(1)$, and then with $X = A_o(1)$ and $Y = A_e(1)$.
\end{proof}

Apart from the Kolmogorov--Rogozin inequality (Lemma~\ref{lem:Kolmogorov--Rogozin}), our main probabilistic tool for the proof of Lemma~\ref{lem:X^2-Y^2 is powerful or square} is Hoeffding's inequality \cite[Theorem~2]{Hoeffding}.

\begin{lemma}[Hoeffding's inequality]
\label{lem:Hoeffding}
Let $X_1, \ldots, X_\ell$ be discrete and independent random variables with support contained in $[-H,H]$. Write $X = X_1 + \cdots + X_\ell$. If $t>0$, then
\begin{equation*}
\mb{P}\Big(\big\lvert X - \mb{E}X \big\rvert \ge t\Big) \le 2 \exp \bigg(-\frac{t^2}{2H^2\ell} \bigg).
\end{equation*}
\end{lemma}

We are now ready to prove Lemma~\ref{lem:X^2-Y^2 is powerful or square}. 

\begin{proof}[Proof of Lemma~\ref{lem:X^2-Y^2 is powerful or square}]
We may assume that $m$ is sufficiently large and that $H \le \sqrt{m}$, as otherwise the conclusion of the lemma is trivial. First, by the Kolmogorov--Rogozin inequality, the event that $Y = 0$ has probability $\ll 1/\sqrt{\epsilon c  m}$, which is smaller than the right-hand side of \eqref{eq:X^2-Y^2 is a square}. Define the set
\begin{equation*}
S = \{(k,r,s) \in \Z^3 : k \neq 0, \text{ and } r>s>0 \text{ are coprime and of opposite parity}\}.
\end{equation*}
If $X^2-Y^2$ is a nonzero square and $Y \neq 0$, then $X = k(r^2+s^2)$ and $Y \in \{\pm 2krs, \pm k(r^2-s^2)\}$ for some tuple $(k,r,s) \in S$. Write $p(k,r,s) \coloneqq \mb{P}(X=k(r^2+s^2)) \mb{P}(Y\in \{\pm 2krs, \pm k(r^2-s^2)\})$. Then 
\begin{equation}
\label{eq:pkrs}
\mb{P}(X^2-Y^2 = \square \neq 0) \leq \mb{P}(Y = 0) + \sum_{(k,r,s) \in S} p(k,r,s).
\end{equation}
Set $t \coloneqq H \sqrt{m \log{m}}$. By Lemma~\ref{lem:Hoeffding}, we have
\begin{equation}
\label{eq:squares_tail}
\mb{P}\Big(\big\lvert X - \mb{E}X \big\rvert \ge t\Big) \leq 2 \exp \bigg( - \frac{m \log{m}}{2m_1} \bigg)  \ll \frac{1}{\sqrt{m}}
\end{equation}
since $m_1 \le m$. The right-hand side of \eqref{eq:squares_tail} is small compared to the right-hand side of \eqref{eq:X^2-Y^2 is a square}. Moreover, if we let $S_t$ be the set of triples in $S$ for which $\abs{k(r^2+s^2)-\mb{E}X} < t$, then
\begin{equation}
\label{eq:pkrs-main}
\sum_{(k,r,s) \in S_t} p(k,r,s) \le 4\abs{S_t} \sup_{a \in \Z} \mb{P}(X=a) \sup_{b \in \Z} \mb{P}(Y=b) \ll \frac{\# S_t}{\epsilon \sqrt{m_1m_2}}\ll \frac{\# S_t}{\epsilon cm}
\end{equation}
by Lemma~\ref{lem:Kolmogorov--Rogozin} and by our assumption that $m_1+m_2=m$ and that $m_1,m_2\ge cm$. 

It remains to bound the size of $S_t$. Let us denote by $S_t' \subset S$ the set of triples with the property $\abs{k(r^2+s^2)} \le t$. By the work of Stronina \cite{Stronina}, we have
\begin{equation}
\label{eq:Stronina}
\# S_t' = c_1 t \log{t} + c_2 t + o\left(\sqrt{t}\right)
\end{equation}
for explicit nonzero constants $c_1, c_2$. Let $x=|\mb{E}X|$ and note that $x\le Hm$ since $X$ is the sum of $\le m$ random variables, each of which is supported in $[-H,H]$. If now $x<2t$, then $S_t \subset S_{3t}'$, which is of size $\ll t \log{t} \ll H\sqrt{m (\log m)^3}$ by \eqref{eq:Stronina} and the assumption $H \le \sqrt{m}$. On the other hand, if $x\ge2t$, then again by \eqref{eq:Stronina} and by the Mean Value Theorem, we find
\[
\# S_t \le \# S_{x+t}' - \# S_{x-t}' \ll t\log(x+t) + \sqrt{x+t} .
\]
Since $2t\le x\le Hm$ here, we conclude that $\# S_t \ll H\sqrt{m (\log m)^3}$. Combining these results with \eqref{eq:pkrs-main} yields
\begin{equation*}
\sum_{(k,r,s) \in S_t} p(k,r,s) \ll \frac{H (\log m)^{3/2}}{\epsilon c \sqrt{m}}
\end{equation*}
which, together with \eqref{eq:squares_tail} and \eqref{eq:pkrs}, proves the lemma.
\end{proof}

\section{Interlude: Euclidean division}
\label{sec:Euclidean division}

In \S \ref{sec:small degree factors} and \S \ref{sec:discriminant_and_exceptional}, we ruled out reciprocal divisors of small degree and the exceptional factorisation $A = \pm I\cdot I_{\ms{rev}}$. The goal of \S \S \ref{sec:Euclidean division}--\ref{sec:anatomy} is to rule out reciprocal divisors of large degree. We follow the proof methods of \cite{BKK}, and in many cases our arguments are direct analogues of those presented there. However, the lack of independence of the coefficients does mean that more work is required in many places.

First off, in this short, intermediate section, we prove a fundamental result regarding division with remainder of a (shifted) reciprocal polynomial $C$ by a nonzero reciprocal polynomial $D$ over $\F_p[T]$. By classical Euclidean division, there exists a quotient $Q \in \F_p[T]$ and a remainder $R \in \F_p[T]$ such that $R=0$ or $\deg{R} < \deg{D}$; typically, $R$ is not a reciprocal polynomial, and it is not immediately clear how many distinct $R$ may occur as such a residue. However, in Definition~\ref{def:RmD} we describe an explicit, complete set $R_m(D)$ of remainders that are shifted reciprocal polynomials. The set $R_m(D)$ is of size $p^{\deg(D)/2}$. Proofs of the results in subsequent sections rely heavily on the description of $R_m(D)$.

Recall the notation $\CR_p(k)$ that was introduced in Definition~\ref{def:CRpm}; it denotes the set of monic reciprocal polynomials of degree $2k$ in the ring $\F_p[T]$. Recall also the notation $\CR_{\F_p}^{\ms{sh}}(m)$ for the set of $m$-shifted reciprocal polynomials over $\F_p$ from Definition~\ref{def:ell-shifted reciprocal}, meaning the set of all polynomials of the form 
\[
c_m + c_{m-1}T + \cdots + c_1 T^{m-1}+ c_0 T^{m} + c_1 T^{m+1} + \cdots + c_{m-1}T^{2m-1} + c_m T^{2m} \in \F_p[T].
\] 
For convenience, we shall write $\CR_{p}^{\ms{sh}}(m) \coloneqq \CR_{\F_p}^{\ms{sh}}(m)$.

\begin{lemma}[Euclidean division for shifted reciprocal polynomials]
\label{lem:reciprocal Euclid} Let $p$ be a prime, $k\in\Z_{\ge1}$ and $m\in\Z_{\ge0}$. In addition, let $C \in \CR_p^{\ms{sh}}(m)$ and $D \in \CR_p(k)$.
There are unique polynomials $Q \in \F_p[T]$ and $R \in T^{m-k+1} \CR_p^{\ms{sh}}(k-1)\cap \F_p[T] \subset \CR_p^{\ms{sh}}(m)$ such that $C = QD + R$. 
\end{lemma}

\begin{remark*}
If $\ell=\min\{k-1,m\}$, then a simple computation shows that $T^{m-k+1} \CR_p^{\ms{sh}}(k-1)\cap \F_p[T]=T^{m-\ell} \CR_p^{\ms{sh}}(\ell)$. The reason why we state Lemma \ref{lem:reciprocal Euclid} in this way is to motivate Definition \ref{def:RmD}, whose precise shape will be crucial in \S\S \ref{sec:character_theory_and_special_residues}--\ref{sec:L_1_bounds}. 

At any rate, for the purposes of establishing Theorem~\ref{thm:master}, we only need to control $D \in \CR_p(k)$ with $k$ at most slightly larger than $m/2$. Hence, we will always have that $k \le m+1$. However, in certain auxiliary results we will consider products of the form $DD'$ with $D \in \CR_p(k)$, $D' \in \CR_p(k')$, in which case $k+k'$ can exceed $m+1$.
\end{remark*}

\begin{proof}[Proof of Lemma~\ref{lem:reciprocal Euclid}]
Let $\ell=\min\{k-1,m\}$. By the above remark, we must show there are unique polynomials $Q \in \F_p[T]$ and $R \in T^{m-\ell} \CR_p^{\ms{sh}}(\ell)$ such that $C = QD + R$. If such $Q$ and $R$ exist, their uniqueness is immediate, because if $C=QD+T^{m-\ell}R_1=Q'D+T^{m-\ell}R_1'$ with $R_1,R_1'\in \CR_p^{\ms{sh}}(\ell)$, then $T^{m-\ell} R_1\equiv T^{m-\ell}R_1'\mod D$. We also know that $T\nmid D$ because $D\in \CR_p(k)$. Thus $R_1\equiv R_1'\mod D$. Since $R_1$ and $R_1'$ are either $0$ or have degrees $\le 2\ell<\deg(D)$, we conclude that $R_1=R_1'$ as needed.

Let us now prove the existence of $Q$ and $R$. If $k>m$, we may trivially take $R=C$ and $Q=0$. Now suppose $k\le m$. Recall from Definition~\ref{def:inverse general reciprocal map} that $(\boldsymbol{-})_{\SR, m}$ denotes the inverse of the map $(\boldsymbol{-})^{\SR, m}$, which by Lemma \ref{lem:trace pol iso of vs} is an isomorphism between the $K$-vector spaces $K[T]_{\le m}$ and $\CR_K^{\ms{sh}}(m)$. Hence by Euclidean division, there exist $Q', R'$ with $C_{\SR, m}=Q'D_{\SR} + R'$ and $\deg{R'} < \deg{D_{\SR}} = k\le m$ or $R' = 0$. Moreover, if $Q'\neq0$, then $\deg(Q')\le m-k$. Multiplying by $T^{m}$, we find that
\begin{align*}
C(T) = T^{m} C_{\SR, m}(T+T^{-1}) &= T^{m}Q'(T+T^{-1}) D_{\SR}(T+T^{-1}) + T^{m} R'(T+T^{-1}) \\
&= T^{m-k}Q'(T+T^{-1}) D(T) + T^{m} R'(T+T^{-1}).
\end{align*}
Hence the claim holds with $Q = (Q')^{\SR, m-k}$ and $R = (R')^{\SR, m}$.
\end{proof}

\begin{definition}\label{def:reciprocal remainder}
In the notation of Lemma~\ref{lem:reciprocal Euclid}, we say that $R$ is the \emph{reciprocal remainder} of $C$ modulo $D$, and write $R \coloneqq C \brmod{D}$.
\end{definition}

As $T^{m-k+1} \CR_p^{\ms{sh}}(k-1) \subset \CR_p^{\ms{sh}}(m)$ when $m \ge k-1$, Lemma~\ref{lem:reciprocal Euclid} has the following immediate consequence.
\begin{lemma}
	\label{lem:image coincidence}
	Let $p$ be a prime, let $k\in\Z_{\ge1}$, let $D \in \CR_p(k)$ and let $m \ge k-1$ be an integer. Then the images of $T^{m-k+1} \CR_p^{\ms{sh}}(k-1)$ and $\CR_p^{\ms{sh}}(m)$ in $\F_p[T]/(D)$ coincide. 
\end{lemma}

\begin{definition}
	\label{def:RmD}
	Let $p$ be a prime, let $k\in\Z_{\ge0}$, let $D \in \CR_p(k)$ and let $m \in\Z_{\ge0}$.
	If $k\ge1$, we denote by $R_m(D)$ the image of $T^{m-k+1}\CR_p^{\ms{sh}}(k-1)$ in $\F_p[T]/(D)$. If $k=0$, we let $R_m(D)=\{0\}$. 
\end{definition}
\begin{remark*}
Since $T$ and any (reciprocal) $D$ are coprime, there exists a multiplicative inverse of $T$ in $\F_p[T]/(D)$, so Definition~\ref{def:RmD} also makes sense when $m\le k-2$.
\end{remark*}

We will also need to prove the technical Lemma \ref{lem:technical Euclid J>0} below, which we will use in \S \ref{sec:L infinity bounds}. First, we need an auxiliary result concerning the Chebyshev polynomials $C_j$ (see Definition~\ref{def:chebyshev polynomials}). 

\begin{lemma}
	\label{lem:coprime chebyshevs}
	Let $p$ be a prime and let $a, b \in \Z_{>0}$.
	\begin{enumerate}[label=(\alph*)]
		\item Suppose $p$ is odd, and $t/s \neq p^{\lambda}$ for all $\lambda \in \Z$ and for all positive integers $t$ and $s$ satisfying $t \mid 4a$, $t\nmid 2a$, $s \mid 4b$ and $s \nmid 2b$. Then the Chebyshev polynomials $C_a$ and $C_b$ are coprime in $\F_p[T]$. In particular, if $a$ and $b$ have a different $2$-valuation, then $C_a$ and $C_b$ are coprime in $\F_p[T]$.
		\item Suppose $p=2$. Then $C_a/T$ and $C_b/T$ are coprime in $\F_p[T]$ if $2|(a+1)(b+1)$ and $\gcd(a,b)=1$. In particular, $C_a/T$ and $C_{a+1}/T$ are coprime.
	\end{enumerate}
\end{lemma}
\begin{proof}
	(a) If $T^{2a}+1$ and $T^{2b}+1$ are coprime, then their trace polynomials $C_a$ and $C_b$ are also coprime: indeed, if the latter have a common factor $I$, then $I^{\SR}$ divides both $T^{2a}+1$ and $T^{2b}+1$. Now, writing $\Phi_d$ for the $d$-th cyclotomic polynomial, the $\Q$-irreducible factorisation of $T^{2k}+1$ is given by
	\begin{equation}
		\label{eq:cheby cyclo factorisation}
		T^{2k}+1 = \frac{T^{4k}-1}{T^{2k}-1} = \prod_{\substack{d \mid 4k \\ d \nmid 2k}} \Phi_d.
	\end{equation}
	Recall that two polynomials $A$ and $B$ have a common factor if and only if their resultant $\Res(A,B)$ vanishes. Furthermore, the function $\Res$ is multiplicative in both arguments up to sign. Therefore $T^{2a}+1$ and $T^{2b}+1$ are coprime if and only if $\Res(\Phi_t, \Phi_s) \neq 0$ for all $t \mid 4a$, $t\nmid 2a$ and $s \mid 4b$, $s \nmid 2b$. Now, over $\Q$, we have
	\begin{equation*}
		\Res(\Phi_t, \Phi_s) =
		\begin{cases}
			0 & \text{if } t=s \\
			q^{\varphi(\min\{s,t\})} & \text{if } t/s = q^{\lambda} \text{ for some } \lambda \in \Z \setminus \{0\} \text{ and prime } q \\
			1 & \text{else},
		\end{cases}
	\end{equation*}
	where $\varphi$ is Euler's totient function (see, for example, \cite[\S 3.3.6]{Prasolov}). Over $\F_p$, it follows that $\Res(\Phi_t, \Phi_s) \neq 0$ if $t/s \neq p^{\lambda}$ for all $\lambda \in \Z$. Combining this with \eqref{eq:cheby cyclo factorisation} and our assumptions proves the first claim in part (a). The second claim immediately follows from this.
	
	\medskip
	
	(b) By the same reasoning as above, we find that the polynomials $C_a/T$ and $C_b/T$ are coprime if the polynomials $(T^{2a}+1)/(T^2+1)$ and $(T^{2b}+1)/(T^2+1)$ are coprime. Since we work modulo $2$ here, it suffices that the polynomials
	\begin{equation*}
		\frac{T^{2a}-1}{T^2-1} = \prod_{\substack{t \mid 2a \\ t\nmid 2}} \Phi_t \qquad \text{and} \qquad \frac{T^{2b}-1}{T^2-1} = \prod_{\substack{s \mid 2b \\ s\nmid 2}} \Phi_s 
	\end{equation*}
	are coprime. Let us write $t=2^vt_1$ and $s=2^ws_1$ with $t_1,s_1$ odd. By the discussion in part (a), we would only have a common factor if $t/s=2^\lambda$ for some $\lambda\in\Z$. Equivalently, $t_1=s_1$. Since we have assumed that $\gcd(a,b)=1$, we must then have that $t_1=s_1=1$, meaning that $t$ and $s$ are both powers of $2$. We further know that $a$ or $b$ is odd. If $a$ is odd, then the relation $t=2^v|2a$ implies that $t\in\{1,2\}$, which is impossible. Similarly, if $b$ is odd, we find that $s\in\{1,2\}$, which is again impossible. This completes the proof of the first claim of part (b). The second claim then follows readily.
\end{proof}

\begin{lemma}
	\label{lem:technical Euclid J>0}
	Fix positive integers $j$, $k$, and $m$ with $j+2k-1 \leq m$. Consider the linear subspace $\CR_p(j,k)$ of $\CR_p^{\ms{sh}}(m)$ consisting of all polynomials of the form 
	\begin{equation}
		\label{eq:J,k lacunary pol}
		C(T) = T^m \sum_{i=0}^{2k-1} c_i(T^{j+i} + T^{-j-i}) \quad \text{with } c_0, c_1, \ldots, c_{2k-1} \in \F_p.
	\end{equation}
	Fix $D \in \CR_p(k)$; in case $p=2$, assume in addition that $T^2+1 \nmid D$. Then
	\begin{equation*}
		\CR_p(j,k) \to \CR_p^{\ms{sh}}(m), \quad C \mapsto C \brmod{D}
	\end{equation*}
	is a linear map that surjects onto $T^{m-k+1}\CR_p^{\ms{sh}}(k-1)$.
\end{lemma}
\begin{proof}
	Let $i\geq 0$ be an integer. The Chebyshev polynomial $C_{j+i} \in \F_p[T]$ is monic of degree $j+i>0$. Hence the span $V$ of the polynomials $C_{j}, C_{j+1}, \ldots, C_{j+2k-1}$ over $\F_p$ is of dimension $2k$. Since $k\le m$ here, tracing the proof of Lemma~\ref{lem:reciprocal Euclid} we find that we may split the map that sends an element $C\in \CR_p(j,k)$ to its reciprocal remainder  $C \brmod{D}$ into three pieces:
	\begin{equation}
		\label{eq:Lacunary maps sequence}
		\begin{tikzcd}[row sep=tiny, column sep=3.5em]
			\CR_p(j,k) & V & \F_p[T]_{\leq k-1} & T^{m-k+1} \CR_p^{\ms{sh}}(k-1) \\
			C & C_{\SR, m} & C_{\SR, m} \bmod{D_{\SR}} & C \brmod{D}.
			\arrow["(\boldsymbol{-})_{\SR, m}", from=1-1, to=1-2]
			\arrow["\bmod{D_{\SR}}", from=1-2, to=1-3]
			\arrow["(\boldsymbol{-})^{\SR, m}", from=1-3, to=1-4]
			\arrow[maps to, from=2-1, to=2-2]
			\arrow[maps to, from=2-2, to=2-3]
			\arrow[maps to, from=2-3, to=2-4]
		\end{tikzcd}
	\end{equation}
	The third map is a linear isomorphism by Lemma~\ref{lem:trace pol iso of vs}. The first map is also a linear isomorphism by Lemma~\ref{lem:trace pol iso of vs}; indeed, domain and codomain are of equal dimension over $\F_p$, and given $C$ as in \eqref{eq:J,k lacunary pol}, we have 
	\begin{equation*}
		C_{\SR, m}(T) = \sum_{i=0}^{2k-1} c_iC_{j+i}(T) \in V.
	\end{equation*}
	It remains to show that the linear map in the middle of \eqref{eq:Lacunary maps sequence} is surjective. 
	
	Firstly, using \eqref{eq:chebyshev recursion general} and induction on $\ell$, we have the formula
	\begin{equation}
		\label{eq:Chebyshev powers}
		T^{\ell} C_a = \sum_{i=0}^{\ell} \binom{\ell}{i} C_{a+\ell-2i}
	\end{equation}
	for any integers $a$ and $\ell\geq0$. Let $\ell \in \{0, 1, \ldots, k-1\}$. In particular, we find that the polynomials $T^{\ell} C_{j+k-1}(T)$ and $T^{\ell} C_{j+k}(T)$ lie in $V$. 
	
	Secondly, note that the map $\bmod \, D_{\SR} \colon V \to \F_p[T]_{\leq k-1}$ is the same as first sending an element of $V$ to $\F_p[T]/(D_{\SR})$ by reduction modulo $D_{\SR}$, and then picking the unique representative of its image that is of degree $<k$. Now, if $p\neq 2$ then $C_{j+k-1}$ and $C_{j+k}$ are coprime over $\F_p$ by Lemma~\ref{lem:coprime chebyshevs}(a). Hence there exist polynomials $F, G \in \F_p[T]$ such that $FC_{j+k-1} + GC_{j+k} = 1$ in $\F_p[T]$, and the same equality also holds in $\F_p[T]/(D_{\SR})$. Denote by $F_{\ell}$ the unique lift of $T^{\ell} F \bmod{D_{\SR}} \in \F_p[T]/(D_{\SR})$ to $\F_p[T]$ that is of degree $<k$, and let $G_{\ell}$ be a similar lift of $T^{\ell}G$. Then the polynomials $F_{\ell}C_{j+k-1} + G_{\ell}C_{j+k} \in \F_p[T]$ with $\ell=0,1,\dots,k-1$ all lie in $V$ by the discussion in the preceding paragraph. 	Therefore the reduction $F_{\ell}C_{j+k-1} + G_{\ell}C_{j+k} = T^{\ell} \bmod{D_{\SR}} \in \F_p[T]_{\leq k-1}$ lies in the image of $V$ under the reduction map modulo $D_{\SR}$. This finishes the proof for odd $p$.
	
	Similarly, if $p=2$, then $(C_{j+k-1}, C_{j+k}) = T$ by Lemma~\ref{lem:coprime chebyshevs}(b). So there exist polynomials $F'$ and $G'$ such that $F'C_{j+k-1} + G'C_{j+k} = T$ in $\F_p[T]$. Since $T^2+1$ does not divide $D$, the trace polynomial $D_{\SR}$ is coprime with $T$, so there exists a polynomial $H \in \F_p[T]$ so that $HF'C_{j+k-1} + HG'C_{j+k} = 1 \in \F_p[T]/(D_{\SR})$. We now set $F = HF'$ and $G = HG'$ and repeat the same steps as in the case of $p \neq 2$ above to finish the proof for $p=2$ as well.
\end{proof}

\section{Character theory and special residues}
\label{sec:character_theory_and_special_residues}

In \S\S \ref{sec:character_theory_and_special_residues}--\ref{sec:L_1_bounds}, we work towards the proof of Proposition~\ref{prop:Delta bound}. The present section consists of preparatory work for the Fourier analysis carried out in \S\S \ref{sec:Fourier transform}--\ref{sec:L_1_bounds}.

Denote by $\CP = \{p_1, \ldots, p_r\}$ a set of $r$ (distinct) primes. For reasons explained in \S \ref{sec:main_results_and_proof_strategy}, we need to control the joint divisor structure distribution of $\BA = (A_p)_{p \in \CP}$ inside $\F_{\CP}[T] \coloneqq \prod_{p \in \CP} \F_p[T]$. We first introduce the theoretical framework, mostly following \cite{BKK}. Let $\Bk = (k_p) \in \Z^r_{\ge 0}$ be a tuple of nonnegative integers. Denote by $\BD = (D_p)_{p \in \CP}$ a tuple of reciprocal polynomials in $\F_{\CP}[T]$, with $D_p \in \CR_p(k_p)$; briefly, we will write $\BD \in \CR_\CP(\Bk)$ for this datum. Define $\F_\CP[T]/(\BD) \coloneqq \prod_{p \in \CP} \F_p[T]/(D_p)$. In similar fashion, denote $\BB \bmod{\BD} \coloneqq (B_p \bmod{D_p})_{p \in \CP} \in \F_\CP[T]/(\BD)$.  

As $\F_{\CP}[T]/(\BD)$ is an abelian group for addition, we may apply Fourier inversion to expand the indicator function of $\BA \equiv \BC \bmod{\BD}$ as a sum of characters.
Denote by $\F_p((1/T))$ the field of Laurent series $X(T) = \sum_{-\infty < j \le J} c_j T^j$ where $J \in \Z$ and $c_j \in \F_p$. For $X \in \F_p((1/T))$, set $\res(X) \coloneqq c_{-1}$. Let 
\begin{equation*}
\psi_{\BD} \colon \F_{\CP}[T]/(\BD) \to \R/\Z, \quad \psi_{\BD}(\BB) \coloneqq \sum_{p \in \CP} \frac{\res(B_p/D_p)}{p}.
\end{equation*}
Note that this is well-defined, because $\res(B_p/D_p)$ does not change when choosing another polynomial representative for $B_p \bmod{D_p}$. 

In addition, given integers $m,j\ge0$, let $\BU = \BU_{m,j,\BD} \in\F_{\CP}[T]$ be such that $U_p\equiv T^{m-j}+1_{j\ge1} T^{m+j}\mod{D_p}$
for each $p\in\CP$. Such a tuple $\BU$ always exists because $D_p$ is a monic reciprocal polynomial, and thus $T\nmid D_p$ for each $p\in\CP$. We then define the maps
\begin{equation*}
	\psi^{(m,j)}_{\BD} \colon \F_{\CP}[T]/(\BD) \to \R/\Z, 
		\quad 
\psi^{(m,j)}_{\BD}(\BB) \coloneqq 
\psi_{\BD}( \BU\BB) .
\end{equation*}
Note that this is well-defined no matter which representatives $U_p$ and $B_p$ we choose, because $\res(U_pB_p/D_p)$ stays stable mod $p$ for all choices. 

Lastly, recall the notation $e(x) \coloneqq \exp(2\pi \ii x)$.

\medskip

We conclude this subsection by recording the following simple compatibility lemma.

\begin{lemma}
	\label{lem:psi-compatibility}
	Let $\CP$ be a finite set of primes, let $\BD,\BD'\in\CR_{\CP}$ and let $\BB\in\F_{\CP}[T]$. We have
	\[
	\psi_{\BD}(\BB) \equiv \psi_{\BD \BD'}(\BB \BD') \mod 1
	\quad\text{and}\quad
	\psi^{(m,j)}_{\BD}(\BB) \equiv \psi^{(m,j)}_{\BD \BD'}(\BB \BD') \mod 1
	\]
	for all integers $m,j\ge0$.
\end{lemma}

\begin{proof}
	The first identity is immediate from the definition. For the second one, let $\BU\in \F_{\CP}[T]$ be such that $U_p\equiv T^{m-j}+1_{j\ge1}T^{m+j}\mod {D_pD_p'}$ for all $p\in\CP$. It then also follows that $U_p\equiv T^{m-j}+1_{j\ge1} T^{m+j} \mod{D_p}$. Hence, the second identity follows too.
	\end{proof}

\subsection*{Major and minor residues}
\label{sec:major residues}

Fix $m\in\Z_{\ge0}$, $\Bk \in \Z_{\ge 0}^r$ and $\BD \in \CR_\CP(\Bk)$. Denote by $R_m(\BD)$ the Cartesian product
\begin{equation}
\label{eq:def RmD tuple}
R_m(\BD) \coloneqq \prod_{p \in \CP} R_m(D_p) \subset \F_\CP[T]/(\BD),
\end{equation}
with $R_m(D_p)$ as in Definition~\ref{def:RmD}.
If $k_p \le m+1$ for all $p \in \CP$, Lemma~\ref{lem:image coincidence} ensures that $R_m(\BD)$ is also the image of $\CR^{\ms{sh}}_\CP(m)$ inside $\F_{\CP}[T]/(\BD)$. We have $\BA \bmod \BD \in R_m(\BD)$ for any $\BD$, but if $k_p > m+1$ for some $p$, then $\BA \bmod \BD$ will a priori (i.e., without knowing anything about the probability distribution on $\BA$) already lie in a proper subset of $R_m(\BD)$. 

The set $R_m(\BD)$ is a subgroup of the additive abelian group of $\F_\CP[T]/(\BD)$. Moreover, the usual multiplication by the ring $\F_{\CP}$ endows it with a module structure over $\F_\CP$. 

Let us now consider the maps
\begin{equation*}
\widetilde{\chi}_{\BB} \colon \F_{\CP}[T]/(\BD) \to \C^{\times}, \quad \BC \mapsto e(\psi_{\BD}(\BC\BB))
\end{equation*}
with $\BB$ varying over a complete set of residue classes modulo $\BD$. They form a complete set of characters of the group $\F_\CP[T]/(\BD)$. We then introduce the following important notions:

\begin{definition}[Major and minor residues]
\label{def:Major residues}
Let $m$, $\CP$ and $\BD$ be as above. 
We write $N_m(\BD)$ for the set of residues $\BB \in \F_\CP[T]/(\BD)$ for which the restriction $\chi_{\BB} \coloneqq \widetilde{\chi}_{\BB} \big\vert_{R_m(\BD)}$ is trivial. We call these the \emph{major residues} modulo $\BD$. 

On the other hand, we shall refer to the residues $\BB \bmod{\BD}$ that do not lie in $N_m(\BD)$ as the \emph{minor residues}.
\end{definition}

The set $N_m(\BD)$ is an abelian group. Moreover, since $R_m(\BD)$ is an $\F_{\CP}$-module, so is $N_m(\BD)$. 

\begin{lemma} Let $m$, $\CP$ and $\BD$ be as above. 
\label{lem:character groups}\ 
\begin{enumerate}[label=(\alph*)]
	\item The maps
\begin{equation*}
\chi_{\BB} \colon R_m(\BD) \to \C^{\times}, \quad \BC \mapsto e(\psi_{\BD}(\BC\BB)),
\end{equation*}
with $\BB \bmod{\BD}$ varying over any complete set of representatives of the cosets of $N_m(\BD)$ in $\F_{\CP}[T]/(\BD)$, form a complete set of characters for $R_m(\BD)$. 
\item The set $N_m(\BD)$ has cardinality $\|\BD\|_{\CP}^{1/2}$. 
 \end{enumerate}
\end{lemma}

\begin{proof}
	For a group $G$, denote by $\hat{G}$ its group of characters. If $H \le G$ is a subgroup, set
	\begin{equation*}
		H^{\perp} \coloneqq \{\chi \in \hat{G} : \chi = 1 \textup{ on } H\}.
	\end{equation*}
	
(a) We use a fundamental result from character theory: if $G$ is a finite abelian group and $H\le G$ a subgroup, then each character of $H$ can be extended to a character of $G$. In other words, the restriction map $\hat{G} \to \hat{H}$ defined by $\chi \mapsto \chi\vert_H$ is a surjective homomorphism. It has kernel $H^{\perp}$, so we obtain the isomorphism $\hat{G}/H^{\perp} \cong \hat{H}$. Applying these results to the group $G=\F_\CP[T]/(\BD)$ and the subgroup $H = R_m(\BD)$ yields the claim. This is because the paragraph preceding Definition \ref{def:Major residues} implies that we may identify $\hat{G}$ with $G$, and $H^\perp$ with $N_m(\BD)$.

\medskip 

(b) If $G$ and $H$ are as in the end of the proof of part (a), then the above discussion implies that $|N_m(\BD)|=|H^\perp|=|\hat{G}|/|\hat{H}|=|G|/|H|$, where we used that a finite abelian group is isomorphic to its own character group. Since $|G|=\norm{\BD}_{\CP}=\prod_{p\in\CP}p^{2k_p}$ and $|H|=|R_m(\BD)|=\prod_{p\in\CP}p^{k_p}$, the proof of the lemma is complete.
\end{proof}

Next is a basic result regarding the multiplicative relations between the various sets of major residues.

\begin{lemma}
\label{lem:XmDD' to XmD} Let $p$ be a prime, let $D,D' \in \CR_p$, and let $B$ be a residue modulo $D$.
\begin{enumerate}[label=(\alph*)]
	\item The projection $\F_p[T]/(DD') \to \F_p[T]/(D)$ given by reduction modulo $D$ restricts to a surjection $R_m(DD') \to R_m(D)$.
	\item We have $B \in N_m(D)$ if, and only if, $BD' \in N_m(DD')$.
\end{enumerate}
\end{lemma}

\begin{proof}
(a) Let $\deg(D)=2k$ and $\deg(D')=2k'$. Assume the convention that $\CR_p^{\ms{sh}}(-1)=\{0\}$. Then, for every $C \in R_m(DD')$, there exists a polynomial $F' \in \CR_p^{\ms{sh}}(k+k'-1)$ such that $C \equiv T^{m-k-k'+1}F' \bmod{DD'}$. So $C \equiv T^{m-k-k'+1}F' \bmod{D}$. By Lemma~\ref{lem:image coincidence}, there is a polynomial $F \in \CR_p^{\ms{sh}}(k-1)$ such that $F' \equiv T^{(k+k'-1)-(k-1)}F \equiv T^{k'}F \bmod{D}$. So $C \equiv T^{m-k+1}F \bmod{D}$.

Lastly, we prove that the map $R_m(DD') \to R_m(D)$ is surjective. If $k=0$, then the claim is trivial because $R_m(D)=\{0\}$ in this case. Assume now that $k\ge1$ and consider $C\in R_m(D)$. We thus have $C\equiv T^{m-k+1}F\bmod D$ for some $F\in \CR_p^{\ms{sh}}(k-1)$. The polynomial $F'=T^{k'}F$ is $(k'+k-1)$-shifted reciprocal (see remark (c) following Definition \ref{def:ell-shifted reciprocal}). Hence, if we let $C'$ be the class of $T^{m-k-k'+1}F'$ modulo $DD'$, then $C'\in R_m(DD')$ and $C'\equiv C\bmod D$. This proves the claim that the map $R_m(DD')\to R_m(D)$ is surjective.

\medskip 

(b) Let $B\in N_m(D)$. For each $C' \in R_m(DD')$, we must show that $\psi_{DD'}(C'BD') =0$. Indeed, we have $\psi_{DD'}(C'BD') = \psi_{D}(C'B)$ by Lemma~\ref{lem:psi-compatibility}. Moreover, part (a) implies the existence of $C \in R_m(D)$ such that $C' \equiv C \bmod{D}$. For such $C$, we have $\psi_{D}(C'B)=\psi_{D}(CB)$. Finally, the assumption that $B\in N_m(D)$ yields $\psi_{D}(CB) = 0$.

Conversely, assume that $B$ is a residue modulo $D$ such that $BD'\in N_m(DD')$. For each $C \in R_m(D)$, we must show that $\psi_D(CB) \equiv0\bmod 1$. Indeed, by part (a), there exists $C'\in R_m(DD')$ such that $C'\equiv C\bmod D$. We thus have $\psi_D(CB)\equiv \psi_D(C'B)\equiv \psi_{DD'}(C'BD')\bmod 1$, where the latter equality modulo $1$ follows again from Lemma~\ref{lem:psi-compatibility}. Finally, the expression $\psi_{DD'}(C'BD')$ vanishes because of the assumption that $BD'\in N_m(DD')$. 
\end{proof}

We will also use the following more amenable characterisation of the major residues:
\begin{lemma}
\label{lem:minor residues aux}
Let $p$ be a prime, let $m\in\Z_{\ge0}$, let $k\in\Z_{\ge1}$, and let $D\in\CR_p(k)$. A residue $B \bmod{D}$ lies in $N_m(D)$ if, and only if,
\begin{equation*}
\psi_D^{(m,j)}(B) \equiv0\bmod 1 \quad \text{for}\ j=0,1,\dots,k-1. 
\end{equation*}
\end{lemma}
\begin{proof}
By definition, $B \in N_m(D)$ if, and only if, $\psi_D(CB) \equiv 0\bmod 1$ for all $C \in R_m(D)$. Recall that $R_m(D)$ is the image of $T^{m-k+1} \CR_p^{\ms{sh}}(k-1)$ mod $D$. For $j=0,1,\dots,k-1$, let $U_j\in\F_p[T]$ be such that $U_j\equiv T^{m-j}+1_{j\ge1} T^{m+j}\mod D$. Now, $R_m(D)$ is an $\F_p$-vector space spanned by the classes of $U_j\mod D$ for $j=0,1,\ldots, k-1$, which are linearly independent. We thus conclude that $\psi_D(CB)\equiv0\mod1$ for all $C\in R_m(D)$ if, and only if,
\begin{equation*}
\psi_D \big( (c_0U_0 +\cdots + c_{k-1} U_{k-1})B \big)  \equiv0\bmod 1
\end{equation*}
for all $c_0, \ldots, c_{k-1} \in \F_p$. By $\F_p$-linearity of $\psi_p$, this is equivalent to having $\psi_D( U_j B) \equiv \psi_D^{(m,j)}(B) \equiv 0\bmod1$ for all $j \in \{0,1,\ldots,k-1\}$. 
\end{proof}

\subsection*{Good representatives}
\label{sec:good representatives}

Consider the quotient 
\[
\CL_m(\BD) \coloneqq (\F_\CP[T]/(\BD))/N_m(\BD)
\] 
which is an $\F_\CP$-module as well. In fact, in view of part (b) of Lemma \ref{lem:character groups}, for each $p\in\CP$ the component
	\[
	\CL_m(D_p) \coloneqq (\F_p[T]/(D_p))/N_m(D_p)
	\] 
of $\CL_m(\BD)$ is a vector space of dimension $\deg(D_p)/2$ over $\F_p$. Fix such $p$ and write $D = D_p$. The goal of this subsection is to establish the existence of a set of `good' representatives for the quotient $\CL_m(D)$ of $N_m(D)$ inside $\F_p[T]/(D)$. More precisely, we find such representatives for all $D$ at once. 

For the following definition, recall the notation $(B,D)_{\Br}$ from Definition~\ref{def:reciprocal gcd}; it denotes the monic reciprocal polynomial of maximal degree dividing both $B$ and $D$.

\begin{definition}[Good representatives]
\label{def:representatives}
A collection of sets
\begin{equation*}
L_m \coloneqq \big\{L_m(D) : D \in \CR_p \big\} 
\end{equation*}
is a \emph{system of good representatives} if it satisfies the following properties for all $D \in \CR_p$:
\begin{enumerate}
	\item $L_m(D) \subset \F_p[T]/(D)$ is a complete set of representatives for the cosets of $N_m(D)$ in $\F_p[T]/(D)$;
    \item if $B \in L_m(D)$, then $\deg{(B,D)_{\Br}}$ is maximal among all $B' \in B+N_m(D)$;
    \item if $B \in L_m(D)$, then $B/(B,D)_{\Br} \in L_m(D/(B,D)_{\Br})$.
\end{enumerate}
\end{definition}

\begin{lemma}
\label{lem:Lm exists}
A system $L_m$ of good representatives exists.
\end{lemma}
\begin{proof}
Any coset of $N_m(D)$ can certainly be represented by some $B \in \F_p[T]/(D)$ that satisfies property (2). Thus we can adhere to properties (1) and (2) of Definition~\ref{def:representatives}, and it remains to show that we can achieve property (3) simultaneously. To prove the claim, we will construct each layer
\begin{equation*}
L_{m,k} \coloneqq \big\{L_m(D) \subset \F_p[T]/(D) : D \in \CR_p(k) \big\} \subset L_m
\end{equation*}
separately, by induction on $k$. 

For $k=0$ we only have $D=1$ with the residue $B = 0$, and indeed $L_m(1) = \{0\}$ satisfies all properties. Next, suppose we have constructed layers $L_{m,j}$ for all $j < k$. Let $D \in \CR_p(k)$ and suppose we have a set $L_m'(D)$ satisfying properties (1) and (2), that is to say, we know that: (1) $L_m'(D)$ is a complete set of representatives for the cosets of $N_m(D)$ in $\F_p[T]/(D)$; (2) if $B' \in L_m'(D)$, then $\deg{(B',D)_{\Br}}$ is maximal among all $B''\in B'+N_m(D)$. We now show how to construct $L_m(D)$ that also satisfies property (3).

Let $B' \in L_m'(D)$ and write $K = (B',D)_{\Br}$. We will find some $B\equiv B'\bmod{N_m(D)}$ for which properties (2) and (3) are both satisfied, so we will use this $B$ as the choice in $L_m(D)$ of the representative of the class of $B'$.
\begin{itemize}
	\item If $K=1$, property (3) is trivially true, so we simply take $B=B'$. 
	\item If $K \neq 1$, set $B'=GK$ and $D = HK$. Since $\deg{H} < \deg{D}$, there exists $N \in N_m(H)$ such that $G+N \in L_m(H)$ by the induction hypothesis and property (1) applied to $H$. We claim that we may take $B=B'+NK$. Indeed, we have $NK \in N_m(D)$ by Lemma~\ref{lem:XmDD' to XmD}(b), and thus $B\equiv B'\bmod {N_m(D)}$ as needed. Moreover, we have $(B'+NK,D)_{\Br} = K=(B',D)_{\Br}$, which is maximal within the coset of $B'$, and thus of $B$, modulo $N_m(D)$. Lastly, we have $B/K=G+N \in L_m(H)$, so $B$ satisfies property (3) as well.
\end{itemize}
This completes the inductive construction of $L_m$. 
\end{proof}
The proof of Lemma~\ref{lem:Lm exists} does not indicate whether such a system $L_m$ is unique. In fact, it is not: for example, if $D$ is irreducible, $L_m(D)$ can be any complete set of representatives containing $0$ as properties (2) and (3) are trivial. However, the proof does show that our freedom in constructing $L_m$ lies at most in the choice of the \emph{minimal residues}:

\begin{definition}[Minimal residues]
\label{def:minimal residues}
If $G \in L_m(H)$ and $(G,H)_{\Br} = 1$, then we call $G$ a \emph{minimal residue} modulo $H$. In this case, we also say that $G/H$ is a \emph{minimal} or \emph{reduced} fraction, that $(G,H)$ is a \emph{minimal pair}, and that the equivalence class of $G$ in $\CL_m(H)$ is a \emph{minimal class}.
\end{definition}
We will consider $L_m$, and hence the sets $L_m(D)$, to be fixed from now on, for all $p$. We then also adopt the notational convention 
\[
L_m(\BD) = \prod_{p \in \CP} L_m(D_p)
\]
for all $r$-tuples $\BD=(D_p)_{p\in\CP}$ such that $D_p\in \CR_p$ for all $p\in\CP$.

\subsection*{Crossing residues}
\label{sec:crossing residues}

This section provides crucial input for the $L^1$ bounds that will be established in \S \ref{sec:L_1_bounds}.

By Lemma~\ref{lem:XmDD' to XmD}(b), there are multiplicative relations between the sets $N_m(D)$ as $D \in \F_p[T]$ varies: if $B_0 \in N_m(D_0)$, then $B_0D_1 \in N_m(D_0D_1)$. (The subscripts no longer refer to indices of the tuple $\BD$.) Swapping the roles of $D_0$ and $D_1$, and recalling that $N_m(D_0D_1)$ is an abelian group, we obtain the map
\begin{equation}
\label{def:Psi}
\Psi\colon N_m(D_0) \times N_m(D_1) \to N_m(D_0D_1), \quad (B_0,B_1) \mapsto B_0D_1-B_1D_0 .
\end{equation}
We may ask to what extent $N_m(D_0D_1)$ can be generated from $N_m(D_0)$ and $N_m(D_1)$ alone:

\begin{lemma}
\label{lem:NmH NmH' NmHH'}
Let $p$ be a prime, and let $D_0 ,D_1\in\CR_p$ be coprime. Then the map $\Psi$ defined in \eqref{def:Psi} is a bijection.
\end{lemma}
\begin{proof}
Since the cardinalities of $N_m(D_0) \times N_m(D_1)$ and $N_m(D_0D_1)$ are finite and equal by Lemma~\ref{lem:character groups}(b), it suffices to show that the map $\Psi$ is injective. Let $(B_0, B_1), (B_0', B_1') \in N_m(D_0) \times N_m(D_1)$ and suppose $B_0D_1-B_1D_0 = B_0'D_1-B_1'D_0$. Then $(B_0-B_0')D_1 = (B_1-B_1')D_0$. Since $D_0$ and $D_1$ are coprime, we find that $D_0$ divides $B_0-B_0'$ and $D_1$ divides $B_1-B_1'$. Hence $B_0 = B_0'$ and $B_1 = B_1'$.
\end{proof}

\begin{definition}[Crossing residues]
\label{def:crossing residues}
Let $p$ be a prime, let $H_0, H_1 \in \CR_p$, and consider residues $G_0 \bmod{H_0}$ and $G_1 \bmod{H_1}$. We then say that the pairs $(G_0,H_0)$ and $(G_1,H_1)$ are (a pair of) \emph{crossing residues} if $G_0H_1-G_1H_0 \in N_m(H_0H_1)$. 
\end{definition}

The next statement, a supplement to Lemma~\ref{lem:NmH NmH' NmHH'}, says that there are no nontrivial minimal crossing residues.

\begin{lemma}
\label{lem:crossing residues}
Let $p$ be a prime, let $H_0, H_1 \in \CR_p$, and consider minimal residues $G_0 \in L_m(H_0)$ and $G_1 \in L_m(H_1)$. If $G_0H_1-G_1H_0 \in N_m(H_0H_1)$, then $H_0 = H_1$ and $G_0 = G_1$.
\end{lemma}
\begin{proof}
Write $K = (H_0,H_1)_{\Br}=(H_0,H_1)$ and fix a factorisation $K=K_0K_1$ such that $H_0/K_0$ and $H_1/K_1$ are both reciprocal and coprime to each other (cf.~Proposition \ref{prop:factorisation options}). Lemma~\ref{lem:XmDD' to XmD}(b) yields that
\begin{equation*}
G_0 \frac{H_1}{K} - G_1 \frac{H_0}{K} \in N_m\Big(\frac{H_0H_1}{K}\Big) = N_m\Big(\frac{H_0}{K_0} \frac{H_1}{K_1}\Big).
\end{equation*}
Hence by Lemma~\ref{lem:NmH NmH' NmHH'}, there are $\bar{G}_0 \in N_m(H_0/K_0)$ and $\bar{G}_1 \in N_m(H_1/K_1)$ such that
\begin{equation*}
G_0 \frac{H_1}{K} - G_1 \frac{H_0}{K} = \bar{G}_0 \frac{H_1}{K_1} - \bar{G}_1 \frac{H_0}{K_0} = \bar{G}_0 K_0 \frac{H_1}{K} - \bar{G}_1 K_1 \frac{H_0}{K}. 
\end{equation*}
Therefore $(G_0-\bar{G}_0K_0) H_1/K = (G_1-\bar{G}_1K_1)H_0/K$. Since $H_1/K$ and $H_0/K$ are also coprime, the polynomial $H_0/K$ divides $G_0-\bar{G}_0K_0$. 
On the other hand, we have $\bar{G}_0K_0 \in N_m(H_0)$ by Lemma~\ref{lem:XmDD' to XmD}(b), so $G_0$ and $G_0 - \bar{G}_0 K_0$ belong to the same class in $\CL_m(H_0)$. Since $G_0 \in L_m(H_0)$ is minimal, $G_0$ and $H_0$ are reciprocally coprime. Hence the second property of Definition~\ref{def:representatives} of the system $L_m$ guarantees that $G_0-\bar{G}_0K_0$ and $H_0$ are also reciprocally coprime. So $H_0/K$, itself being reciprocal, can only divide $G_0-\bar{G}_0K_0$ if $K = H_0$. 

By the same argument applied to $H_1$ instead of $H_0$, we find that $K=H_1$. In particular, $H_0=H_1$. Hence $G_0, G_1 \in L_m(H_0)$ and $(G_0 - G_1)H_1 \in N_m(H_0H_1)$. Using Lemma \ref{lem:XmDD' to XmD}(b), we conclude that $G_0-G_1\in N_m(H_0)$. Thus, $G_0$ and $G_1$ represent the same class in $\CL_m(H_0)$, implying that $G_0 = G_1$.
\end{proof}

\section{Fourier analysis on $\CR_{\CP}$}
\label{sec:Fourier transform}

In this section, we develop the necessary tools from Fourier analysis for reciprocal polynomials in order to establish Proposition \ref{prop:Delta bound}.

\begin{lemma}[Fourier inversion]
\label{lem:Fourier inversion}
Let $\CP$ be a finite set of primes, let $\BD=(D_p)_{p\in\CP}$ be such that $D_p\in \CR_p$ for each $p$, let $m\in\Z_{\ge0}$, and let $\BC \in R_m(\BD)$. Then
\begin{equation*}
\frac{1}{\norm{\BD}_{\CP}^{1/2}} \sum_{\BB \in L_m(\BD)} e(\psi_{\BD}(\BC\BB) ) = 1_{\BC \equiv \mathbf{0} \bmod{\BD}}.
\end{equation*}
\end{lemma}
\begin{proof}
On the one hand, we have
\begin{equation*}
\sum_{\BB \bmod{\BD}} e(\psi_{\BD}(\BC\BB)) = \sum_{\BB \in L_m(\BD)} \sum_{\BB' \in N_m(\BD)} e(\psi_{\BD}(\BC(\BB+\BB'))).
\end{equation*}
Since $e(\psi_{\BD}(\BC\BB)) = e(\psi_{\BD}(\BC(\BB+\BB')))$ for any $\BB' \in N_m(\BD)$, and we know from Lemma \ref{lem:character groups}(b) that $N_m(\BD)$ contains $\|\BD\|_{\CP}^{1/2}$ elements, we find that
\begin{equation}\label{eq: orthogonality lemma}
\frac{1}{\norm{\BD}_{\CP}} \sum_{\BB \bmod{\BD}} e(\psi_{\BD}(\BC\BB)) = \frac{1}{\norm{\BD}_{\CP}^{1/2}} \sum_{\BB \in L_m(\BD)} e(\psi_{\BD}(\BC\BB)).
\end{equation}

On the other hand, the orthogonality relations for characters mod $p$ imply that the sum on the left-hand side of \eqref{eq: orthogonality lemma} is equal to the indicator function $1_{\forall \BB \bmod{\BD} \colon \psi_{\BD}(\BC\BB) = 0} = 1_{\BC \equiv \Bzero \bmod{\BD}}$.
\end{proof}

Let $A \in \Z[T]$ be a random reciprocal polynomial selected according to $\mb{P}_{\CR(m)}$ and denote by $\BA$ the induced tuple in $\CR_{\CP}(m)$.  Fixing $\BD \in \CR_{\CP}(\Bk)$ and $\BC \in R_m(\BD)$, and replacing $\BC$ in Lemma~\ref{lem:Fourier inversion} by $\BA-\BC$, we find
\begin{align}
\label{eq:FI-2}
\mb{P}_{\CR_{\CP}(m)}\big(\BA \equiv \BC \bmod \BD\big) &= \mb{E}_{\CR_{\CP}(m)} \left[ \frac{1}{\norm{\BD}_{\CP}^{1/2}} \sum_{\BB \in L_m(\BD)}  e(\psi_{\BD}((\BA-\BC)\BB)) \right] \nonumber \\
&= \frac{1}{\norm{\BD}_{\CP}^{1/2}} \sum_{\BB \in L_m(\BD)}  e(\psi_{\BD}(-\BC\BB)) \mb{E}_{\CR_{\CP}(m)} \big[ e(\psi_{\BD}(\BA\BB)) \big].
\end{align}

The remainder of this section follows \cite[\S 4]{BKK} to establish a first step in the proof of Proposition~\ref{prop:Delta bound}, namely the reduction to \eqref{eq:small delta reduction} below.

\begin{lemma}
\label{lem:Fourier expectation}
Let $m\in\Z_{\ge 1}$, let $\BD=(D_p)_{p\in\CP}$ be such that $D_p\in \CR_p$ for each $p$, and let $\BB \in \F_\CP[T]/(\BD)$. We then have
\begin{equation}
\label{eq:Fourier expectation}
\mb{E}_{\BA \in \CR_{\CP}(m)} \big[ e(\psi_{\BD}(\BA\BB)) \big] = e( \psi_{\BD}(\BB + T^{2m}\BB))  \prod_{j=0}^{m-1} \hat{\mu}_j\big(\psi_{\BD}^{(m,j)}(\BB)\big).
\end{equation}
\end{lemma}
\begin{proof} 
We argue similarly to the proof of \cite[Lemma~4.1]{BKK}. Let $A\in\CR(m)$ be a polynomial sampled according to the measure $\mb{P}_{\CR(m)}$, and let $A_p$ be each reduction modulo $p$. We may thus write $A=a_{m}+a_{m-1}T+\cdots+a_{0}T^m+a_{1}T^{m+1}+\cdots+a_{m}T^{2m}$, where $a_{m}=1$. As a consequence, we have
\[
e(\psi_{\BD}(\BA \BB)) = e(a_0 \psi_{\BD}(T^m\BB)) \prod_{j=1}^{m} e(a_j \psi_{\BD}(T^{m-j}\BB + T^{m+j}\BB) ) .
\]
Note that $\psi_{\BD}(T^{m-j}\BB + T^{m+j}\BB) = \psi_{\BD}^{(m,j)}(\BB)$ for $j=1,\dots,m-1$. It is crucial that $m-j\ge0$ here, so that we can take $U_p=T^{m-j}+1_{j\ge1} T^{m+j}$ for all $p$ in the definition of $\psi_{\BD}^{(m,j)}$. Thus
\begin{equation}
	\label{eq:FI-3}
	e(\psi_{\BD}(\BA \BB)) = e\big(a_0 \psi^{(m,0)}_{\BD}(\BB)\big) \prod_{j=1}^{m} e\big(a_j \psi^{(m,j)}_{\BD}(\BB) \big)  .
\end{equation}
We always have $a_m=1$, and $a_0,a_1,\dots,a_{m-1}$ are independent of each other. Hence, taking the expectation on both sides yields the claimed \eqref{eq:Fourier expectation}.
\end{proof}

Given $\BD$ and $\BB$ as in the above lemma, let $\sigma_{\BD}(m; \BB)$ be the absolute value of the right-hand side in \eqref{eq:Fourier expectation}, that is
\begin{equation*}
\sigma_{\BD}(m; \BB) \coloneqq \bigg\vert\prod_{j=0}^{m-1} \hat{\mu}_j\big(\psi_{\BD}^{(m,j)}(\BB)\big)\bigg\vert. 
\end{equation*}
Assume the notation of Lemma \ref{lem:Fourier expectation}. By \eqref{eq:FI-2} and Lemma~\ref{lem:Fourier expectation}, followed by an application of the triangle inequality, we have
\begin{equation*}
\max_{\BC \in R_m(\BD)} \Bigg\vert \mb{P}_{\BA \in \CR_{\CP}(m)}\big(\BA \equiv \BC \bmod \BD\big) - \frac{1}{\norm{\BD}_{\CP}^{1/2}} \Bigg\vert \le \frac{1}{\norm{\BD}_{\CP}^{1/2}} \sum_{\substack{\BB \in L_m(\BD) \\ \BB \neq \mathbf{0}}} \sigma_{\BD}(m; \BB).
\end{equation*}
We call a reciprocal tuple $\BD \in \CR_{\CP}$ \emph{admissible} if either $2 \not \in \CP$, or $2 \in \CP$ and $T^2+1 \nmid D_2$. Note that for reciprocal $\BH, \BK$, the product $\BH \BK$ is admissible if and only if $\BH$ and $\BK$ are both admissible.

Let $k\in\Z_{\ge1}$, recall the definition of $\Delta^\SR_\CP(m; k)$ from \eqref{eq:Delta def}, and recall also Lemma \ref{lem:reciprocal Euclid} and Definition \ref{def:RmD}. Using the above inequality, we find that
\begin{equation}
\label{eq:Deltastar}
\Delta^\SR_\CP(m; k) \le \sum_{\substack{\BD \in \CR_{\CP} \\ \deg{\BD} \le 2k \\ \BD \text{ admissible}}} \frac{1}{\norm{\BD}_{\CP}^{1/2}} \sum_{\substack{\BB \in L_m(\BD) \\ \BB \neq \mathbf{0}}} \sigma_{\BD}(m; \BB).
\end{equation}
Write $L_m(\Bk)$ for the set of all tuples $(\BB, \BD)$ with admissible $\BD \in \CR_\CP(\Bk)$ and $\BB \in L_m(\BD)$. Then
\begin{equation}
\label{eq:RED-1}
\sum_{\substack{\BD \in \CR_{\CP} \\ \deg{\BD} \le 2k \\ \BD \text{ admissible}}} \frac{1}{\norm{\BD}_{\CP}^{1/2}} \sum_{\substack{\BB \in L_m(\BD) \\ \BB \neq \Bzero}} \sigma_{\BD}(m; \BB) = \sum_{\substack{\Bk : k_p \le k \\ \forall p \in \CP}} \frac{1}{\prod_{p \in \CP} p^{k_p}} \sum_{\substack{(\BB,\BD) \in L_m(\Bk) \\ \BB \neq \Bzero}} \sigma_{\BD}(m; \BB).
\end{equation}
For each $\Bl \in \Z_{\ge 0}^r$, let $L_m(\Bk, \Bl)$ be the subset of $L_m(\Bk)$ consisting of all $(\BB, \BD)$ for which the corresponding minimal pair $(\BG, \BH) = (\BB/(\BB, \BD)_{\Br}, \BD/(\BB, \BD)_{\Br})$ lies in $L_m(\Bl)$. Then
\begin{equation}
\label{eq:RED-2}
\sum_{\substack{(\BB,\BD) \in L_m(\Bk) \\ \BB \neq \Bzero}} \sigma_{\BD}(m; \BB) = \sum_{\substack{\Bl \neq \Bzero \\ \ell_p \le k \, \forall p \in \CP}} \sum_{(\BB, \BD) \in L_m(\Bk, \Bl)} \sigma_{\BD}(m; \BB),
\end{equation}
where we used property (3) of Definition \ref{def:representatives}. The set $L_m(\Bk,\Bl)$ is nonempty only if $\ell_p \le k_p$ for all $p \in \CP$. Conversely, if $\ell_p \le k_p$ for all $p \in \CP$, and $\BG/\BH$ is a minimal fraction in $L_m(\Bl)$, then the pair $(\BG\BK, \BH\BK)$ lies in $L_m(\Bk,\Bl)$ if, and only if, $\BK \in \CR_{\CP}(\Bk-\Bl)$ (which is nonempty) and $\BK$ is admissible. Therefore
\begin{align}
\label{eq:RED-3}
\sum_{(\BB, \BD) \in L_m(\Bk, \Bl)} \sigma_{\BD}(m; \BB) &= \sum_{\substack{\BK \in \CR_{\CP}(\Bk-\Bl) \\ \BK \text{ admissible}}} \sum_{\substack{(\BG, \BH) \in L_m(\Bl) \\ (\BG, \BH) \text{ minimal} }} \sigma_{\BH\BK}(m; \BG\BK) \nonumber \\
&\le \prod_{p \in \CP} p^{k_p-\ell_p}  \sum_{\substack{(\BG, \BH) \in L_m(\Bl) \\ (\BG, \BH) \text{ minimal}}} \sigma_{\BH}(m; \BG),
\end{align}
where we used the second identity of Lemma \ref{lem:psi-compatibility} and where we bounded the number of admissible $\BK$ by the number of all $\BK \in \CR_{\CP}(\Bk-\Bl)$. With the definition
\begin{equation}
\label{eq:small delta def}
\delta^{\SR}_{\CP}(m; \Bl) \coloneqq \frac{1}{\prod_{p \in \CP} p^{\ell_p}}  \sum_{\substack{(\BG, \BH) \in L_m(\Bl) \\ (\BG, \BH) \text{ minimal}}} \sigma_{\BH}(m; \BG),
\end{equation}
the relations \eqref{eq:RED-1}, \eqref{eq:RED-2} and \eqref{eq:RED-3} imply
\begin{equation}
\label{eq:RED-4}
\sum_{\substack{\BD \in \CR_{\CP} \\ \deg{\BD} \le 2k \\ \BD \text{ admissible}}} \frac{1}{\norm{\BD}_{\CP}^{1/2}} \sum_{\substack{\BB \in L_m(\BD) \\ \BB \neq \Bzero}} \sigma_{\BD}(m; \BB) \leq \sum_{\substack{\Bk : k_p \le k \\ \forall p \in \CP}} \sum_{\substack{\Bl \neq \Bzero \\ \ell_p \le k_p \, \forall p \in \CP}} \delta^{\SR}_{\CP}(m; \Bl).
\end{equation}
From \eqref{eq:Deltastar} and \eqref{eq:RED-4}, we conclude
\begin{align*}
\Delta^\SR_\CP(m; k) \le  \sum_{\substack{\Bk : k_p \le k \\ \forall p \in \CP}} \sum_{\substack{\Bl \neq \Bzero \\ \ell_p \le k_p \, \forall p \in \CP}} \delta^{\SR}_{\CP}(m; \Bl) \le (k+1)^{2r} \max_{\substack{\Bl = (\ell_p)_{p \in \CP} \neq \Bzero \\ 0 \le \ell_p \le k \, \forall p \in \CP}} \delta^{\SR}_{\CP}(m; \Bl).
\end{align*}
To prove Proposition~\ref{prop:Delta bound}, it thus suffices to show
\begin{equation}
\label{eq:small delta reduction}
\max_{\substack{\Bl = (\ell_p)_{p \in \CP} \neq \Bzero \\ 0 \le \ell_p \le  \gamma m + m^{0.88} \, \forall p \in \CP}}  \delta^{\SR}_{\CP}(m; \Bl) \ll_r m^{-2r} e^{-m^{1/10}}.
\end{equation}

\section{\texorpdfstring{$L^{\infty}$}{L infinity} bounds}
\label{sec:L infinity bounds}

In this section, we prove a bound on $\sigma_{\BD}(m; \BB)$ that is useful when the tuple $\BD$ has a component of small degree. 

\begin{lemma}
\label{lem:L infinity bounds}
Let $m\in\Z_{>0}$, let $\mu_0, \mu_1, \ldots, \mu_{m-1}$ be probability measures on the integers, let $\CP$ be a set of $r$ primes whose product is $P$, and let $\beta \in [0,1]$ be such that
\begin{equation*}
\abs{\hat{\mu}_j(s/P)} \leq \beta \quad \text{ for all } s = 1, 2, \ldots, P-1 \text{ and } j = 0, 1, \ldots, m-1.
\end{equation*}
Let $\BD \in \CR_{\CP}$ be admissible and such that $\deg(D_p)\le 2(m+1)$ for all $p\in\CP$, let $\BB \not \in N_m(\BD)$ be a minor residue class modulo $\BD$, and let $p\in\CP$ and $k \in [1,m+1] \cap \Z$ be such that $2k = \deg{D_p}$ and $B_p \not \in N_m(D_p)$.
Then
\begin{equation*}
\sigma_{\BD}(m; \BB) \leq \beta^{\lfloor (m+k)/(2k)\rfloor}.
\end{equation*}
Furthermore, such $p$ and $k$ exist.
\end{lemma}
\begin{proof} First of all, let us note that for $j=0,1,\dots,m$, we may take $U_p=T^{m-j}+1_{j\ge1}T^{m+j}$ in the definition of $\psi^{(m,j)}$. Hence, 
\begin{equation}
	\label{eq:psi^{(m,j)} when j<m}
\psi^{(m,j)}_{\BD}(\BB) \equiv \sum_{p\in\CP} \frac{\res((T^{m-j}+1_{j\ge1}T^{m+j})B_p/D_p)}{p} \mod 1
\quad\text{for}\ j=0,1,\dots,m-1.
\end{equation}

Next, we deal with the very last claim of the existence of suitable $p$ and $k$. Indeed, since $\BB \not \in N_m(\BD)$, there exists $p \in \CP$ such that $B_p \not \in N_m(D_p)$. Select any such $p$ and set $D \coloneqq D_p$ and $B \coloneqq B_p$ and $k \coloneqq (\deg{D})/2$. Then $k>0$, because for constant $D$ there is no $B$ available with $B \not \in N_m(D)$. 

Given $p$, $k$ as above, we prove the claimed bound on $\sigma_{\BD}(m; \BB)$. For $k=m+1$, we have $\lfloor (m+k)/(2k)\rfloor=0$, so the result holds trivially. Let us now assume that $k\le m$. Throughout, we let $B=B_p$ and $D=D_p$. Since $B \not \in N_m(D)$, Lemma~\ref{lem:minor residues aux} implies that $\psi^{(m,i)}_D(B) \not\equiv0\mod 1$ for some $i\in\Z\cap[0,k)$. This also implies that $\psi^{(m,i)}_{\BD}(\BB)\not\equiv0\mod1$ for the same index $i\in\Z\cap[0,k)$. In particular, we have
\begin{equation}
	\label{eq:L^infty - initial segment}
	\prod_{i=0}^{k-1}  \big| \hat{\mu}_i\big(\psi_{\BD}^{(m,i)}(\BB)\big) \big| \leq \beta.
\end{equation}

Now, if $m \ge k > m-2k$, then $\lfloor (m+k)/(2k) \rfloor = 1$ and the conclusion of the lemma holds by the above. If $k \le m-2k$, then select $j\in\Z\cap[k,m-2k]$ and suppose
\begin{equation}
\label{eq:L infinity bounds lacunary}
\res\left(\frac{(T^{m-i}+T^{m+i})B}{D}\right) \equiv 0\bmod p \, \text{ for all } i\in \{j,j+1,\ldots, j+2k-1\}.
\end{equation}
In the notation of Lemma~\ref{lem:technical Euclid J>0}, by $\F_p$-linearity and admissibility, it follows that $\res(CB/D) \equiv 0\bmod p$ for any $C \in \CR_p(j,k) \subset \CR_p^{\ms{sh}}(m)$.
Write $R = C \brmod{D}$ (cf.~Definition \ref{def:reciprocal remainder}). Since $R$ and $C$ differ by a multiple of $D$ modulo $p$, we have 
\begin{equation}
\label{eq:residue reciprocal residue}
\res\left(\frac{RB}{D}\right) \equiv \res\left(\frac{CB}{D}\right) \equiv 0\bmod p.
\end{equation}
Lemma~\ref{lem:technical Euclid J>0} implies that the map $\CR_p(j,k) \to T^{m-k+1} \CR_p^{\ms{sh}}(k-1)$ given by $C \mapsto C \brmod{D}$ is surjective. Hence \eqref{eq:residue reciprocal residue} holds for any $R \in T^{m-k+1} \CR_p^{\ms{sh}}(k-1)$. Using relation \eqref{eq:psi^{(m,j)} when j<m} and Lemma~\ref{lem:minor residues aux}, this contradicts the assumption $B \not \in N_m(D)$.
It follows that there exists a $j$ in any interval of $[k,m]\cap\Z$ of length $2k$ such that \eqref{eq:L infinity bounds lacunary} fails to hold. That is, there is an $i$ in any interval of $[k,m-1]\cap\Z$ of length $2k$ such that
\begin{equation*}
\res\left(\frac{(T^{m-i}+T^{m+i})B}{D}\right) \not\equiv0\bmod p .
\end{equation*}
Hence there are at least $\lfloor (m-k)/(2k) \rfloor $ values of $i \in [k,m-1] \cap \Z$ such that
\begin{equation*}
\big| \hat{\mu}_i\big(\psi_{\BD}^{(m,i)}(\BB) \big) \big|  \leq \beta.
\end{equation*}
Together with \eqref{eq:L^infty - initial segment}, this completes the proof of the lemma.
\end{proof}

\begin{remark}
The condition $T^2+1 \nmid D_2$ in Lemma~\ref{lem:L infinity bounds} is the origin of the same requirement in the definition of $\Delta^\SR_\CP$, see \eqref{eq:Delta def}.
\end{remark}

\section{$L^1$ bounds}
\label{sec:L_1_bounds}
The purpose of this section is to establish bounds of `$L^1$ shape' that will allow us to prove \eqref{eq:small delta reduction}, and thus Proposition \ref{prop:Delta bound}. Before we proceed, recall the notion of minimal residues (cf.~Definition~\ref{def:minimal residues}).

\begin{definition}
Let $p$ be a prime and $H \in \CR_p$. We denote by $L_m^*(H)$ the set consisting of nonzero minimal residues $G \in L_m(H)$. Furthermore, for every $\ell\ge0$, let $L_m^*(\ell) \coloneqq \{(G,H) : H \in \CR_p(\ell) \text{ admissible}, G \in L_m^*(H)\}$. Similarly, if $\CP$ is a finite set of primes and $\BH \in \CR_{\CP}$, then we let 
\[
L_m^*(\BH) = \big\{(\BG,\BH): (G_p,H_p)\in L_m^*(H_p)\ \forall p\in\CP\big\}.
\]
\end{definition}

Motivated by  \eqref{eq:small delta def}, we prove the following general bound. 

\begin{lemma}
\label{lem:L1 bound} Let $p$ be a prime, let $\nu\in\Z_{\ge0}$, let $f_0, f_1, \ldots , f_{\nu-1}\colon \R/\Z \to \R_{\ge 0}$ be functions, let $m\in\Z_{\ge0}$ and let $\ell \in\Z_{\ge \nu/2}$. Then
\begin{equation}
\label{eq:L1 bound}
\sum_{H \in \CR_p(\ell)} \sum_{G \in L^*_m(H)}  \prod_{j=0}^{\nu-1} f_{j}\big(\psi_H^{(m,j)}(G)\big) \le p^{2\ell-\nu} \prod_{j=0}^{\nu-1} \Bigg( \sum_{\xi=0}^{p-1} f_{j}(\xi/p) \Bigg).
\end{equation}
\end{lemma}

\begin{proof}
Let $X$ denote the set of pairs $(G,H)$ such that $H \in \CR_p(\ell)$ and $G \in L^*_m(H)$. We then define the function $\Upsilon:X\to (\R/\Z)^{2\ell}$ by  $\Upsilon(G,H) = \big(\psi^{(m,j)}_H(G)\big)_{j=0}^{2\ell-1}$. We claim that $\Upsilon$ is injective. Indeed, let us consider two pairs $(G,H),(G',H')$ in $X$ such that $\Upsilon(G,H)=\Upsilon(G',H')$, that is to say $\psi_{H'}^{(m,j)}(G') = \psi_{H}^{(m,j)}(G)$ for $j=0,1,\dots,2\ell-1$. Hence, using Lemma \ref{lem:psi-compatibility} and the linearity of the residue function, we find that
\[
\psi_{HH'}^{(m,j)}(G'H-GH') \equiv \psi_{HH'}^{(m,j)}(G'H) - \psi_{HH'}^{(m,j)}(GH')
\equiv \psi_{H'}^{(m,j)}(G') - \psi_{H}^{(m,j)}(G)\equiv0\mod 1
\]
for $j=0,1,\dots,2\ell-1$. Using Lemma \ref{lem:minor residues aux}, we find that $G'H-GH' \in N_m(HH')$. Thus Lemma~\ref{lem:crossing residues} establishes that $H=H'$ and $G=G'$ as needed. 

Since we have established that $\Upsilon$ is injective, we find that, for any choice of $\xi_0,\dots,\xi_{\nu-1}\in\F_p$, there exist at most $p^{2\ell-\nu}$ choices of $(G,H)\in X$ with $\psi_H^{(m,j)}(G)=\xi_j/p$ for $j=0,1,\dots,\nu-1$. Hence, 
\[
\sum_{(G,H)\in X}  \prod_{j=0}^{\nu-1} f_{j}\big(\psi_H^{(m,j)}(G)\big) \le p^{2\ell-\nu}  
 \sum_{\xi_0,\dots,\xi_{\nu-1}\in\F_p}  \prod_{j=0}^{\nu-1} f_{j} (\xi_j/p)  = p^{2\ell-\nu} \prod_{j=0}^{\nu-1} \Bigg( \sum_{\xi=0}^{p-1} f_{j}(\xi/p) \Bigg),
\]
where we used that $f_j\ge0$ for all $j$. This completes the proof of the lemma.
\end{proof}

\begin{lemma}
\label{lem:bounded FT to bounded delta}
Let $\CP$ be a finite set of primes with $P \coloneqq \prod_{p \in \CP} p$. Suppose $\gamma \ge 1/2$ and $\alpha \ge 0$ are such that
\begin{equation*}
\sum_{k=0}^{Q-1} \abs{\hat{\mu}_j(k/Q + \ell/R)} \le \alpha Q^{1-\gamma}
\end{equation*}
for all $j = 0, 1, \ldots, m-1$ and all $Q, R, \ell \in \Z$ with $QR = P$ and $Q>1$. Let $\Bl \in \Z_{\ge 0}^r$ and $L = \max\{\ell_p : p \in \CP\}$. Then 
\begin{equation*}
\delta^{\SR}_{\CP}(m; \Bl) \le P^{\max \{0, \, L - \gamma m\}} \alpha^{\min\{m,\, 2L\}}.
\end{equation*}
\end{lemma}
\begin{proof}
The needed bound follows immediately by translating the corresponding proof \cite[Lemma~6.3]{BKK} to our setting. In that lemma, there is a parameter $s$; here, we only consider the corresponding case where $s=1$, which makes the proof here slightly easier. In particular, the proof in \cite{BKK} starts by expanding their formula (4.7) (which corresponds to our \eqref{eq:small delta def}), then removing some of the terms using the trivial bound $\abs{\hat{\mu}_j} \le 1$ and applying H\"older. We do not need this step. More precisely, we need to stick with the case corresponding to $s=1$ because we cannot carry out the H\"older procedure, see Remark~\ref{rem:why not s>1} below. Note also that the letter $m$ here corresponds to the parameter $(n-1)/s$ in the notation of \cite{BKK} (which is denoted by $m$ in \cite{BKK}).
\end{proof}

\begin{remark}
\label{rem:why not s>1}
Taking $s>1$ in \cite[Theorem~7]{BKK}, we find that the authors proved results that are, in their strongest form, of type ``a standard polynomial $B$ of degree $n$ has no divisors of degree $\leq \theta n$'' for some absolute $0 < \theta < 1/2$. As commented on in the introduction and in Remark~\ref{rem:thmCthm3}, we do not prove such results here. Indeed, the case $s>1$ in \cite{BKK} requires using H\"older in their Lemma~6.3, which is the analogue of our Lemma~\ref{lem:bounded FT to bounded delta}. After applying H\"older, they can simply shift the residues and obtain their inequality (6.2). In contrast, such shifting of residues $\BG \bmod{\BH}$ would not be possible after applying H\"older to the expression for our $\delta^{\SR}_{\CP}(m; \Bl)$ (roughly \eqref{eq:small delta def}) because such shifts do not always land in $L_m^*(\BH)$. This \emph{can} be circumvented by changing in the definition of our $\Delta^\SR_\CP$ (see \eqref{eq:Delta def}) that the sum only extends over $\BD$ for which $D_p$ is coprime with $T^{2i}+1$ for all $i = 1, 2, \ldots, m$ and all $p \in \CP$. Unfortunately, this leads to further problems in the anatomy (\S \ref{sec:anatomy}).
\end{remark}

We now come to the proof of Proposition~\ref{prop:Delta bound}.

\begin{proof}[Proof of Proposition~\ref{prop:Delta bound}]
We follow the proof of \cite[Proposition~2.3]{BKK}. Recall that it suffices to prove \eqref{eq:small delta reduction}. We denote $L = \max\{\ell_p : p \in \CP\}$.

First suppose $L \le (m/\log{m})^{1/2}/(2P)$.
Using the assumptions with $Q=p \in \CP$ and $\ell = 0$, we find that $\sum_{k=0}^{p-1} \abs{\hat{\mu}_j(k/p)} \le \sqrt{p}$. Combining this with \cite[Equation~(2.8) and Lemma~3.6]{BKK}, we deduce that $\abs{\hat{\mu}_j(k/P)} \le e^{-1/P^2}$ for all $k$ not divisible by $P$ and all $j=0,\ldots,m-1$. Thus we may apply the $L^{\infty}$ bound from Lemma~\ref{lem:L infinity bounds} with $\beta = e^{-1/P^2}$ to obtain
\begin{align*}
\delta^{\SR}_{\CP}(m; \Bl) &\le \prod_{p \in \CP} p^{\ell_p} \max_{(\BG, \BH) \in L^*_m(\Bl)} \sigma_{\BH}(m; \BG) \le P^{L} e^{-\lfloor (m+L)/(2L)\rfloor/P^2} \\ 
&\ll \exp(L \log{P} - m/(2LP^2)).
\end{align*}
Eliminating $L$ and then $P$ using the assumption $P \le m^{1/4}$ shows
\begin{equation*}
L \log{P} - \frac{m}{2LP^2} \le \frac{m^{1/2}(\log{P}-2\log{m})}{2P(\log{m})^{1/2}} \le -\frac{7}{8} m^{1/4} (\log{m})^{1/2},
\end{equation*}
which is $\le -m^{1/5}$ since $m \ge P^4 \ge 16$. Hence
\begin{equation*}
\delta^{\SR}_{\CP}(m; \Bl) \ll e^{-m^{1/5}},
\end{equation*}
proving \eqref{eq:small delta reduction} in this regime. 

In the range
\begin{equation*}
(m/\log{m})^{1/2}/(2P) \le L \le \gamma m
\end{equation*}
we apply the $L^{1}$ bound from Lemma~\ref{lem:bounded FT to bounded delta}, where we may use $\alpha = 1-m^{-1/10}$ as input by our assumptions. This yields
\begin{equation*}
\delta^{\SR}_{\CP}(m; \Bl) \le P^{\max \{0, \, L - \gamma m\}} \alpha^{\min\{m,\, 2L\}} \le \alpha^L \le \exp(-Lm^{-1/10}),
\end{equation*}
which, using $P \leq m^{1/4}$, is $\le \exp(-m^{3/20}/(4\log{m})^{1/2}) \ll \exp(-m^{1/8}) \ll_r m^{-2r} e^{-m^{1/10}}$.

This leaves the range
\begin{equation*}
\gamma m \le L \le \gamma m + m^{0.88},
\end{equation*}
for which we employ the $L^1$ bound from Lemma~\ref{lem:bounded FT to bounded delta} again. In this case we find
\begin{equation*}
\delta^{\SR}_{\CP}(m; \Bl) \le P^{\max \{0, \, L - \gamma m\}} \alpha^{\min\{m,\, 2L\}} \le P^{m^{0.88}} \alpha^m \le \exp(m^{0.88}\log{m}-m^{0.9}) \ll \exp(-m^{0.89}),
\end{equation*}
proving \eqref{eq:small delta reduction} in this range as well.
\end{proof}

\section{Anatomy and large degree factors}
\label{sec:anatomy}

In this section, we prove Proposition~\ref{prop:large degree factors}, which shows how bounds on $\Delta^\SR_\CP$ translate to irreducibility.

Recall Definition~\ref{def:induced_measure_on_trace_pols}, where we defined a measure $\mb{P}_{\SR,\CM(m)}$ on the set $\CM(m)$ of monic degree $m$ polynomials in $\Z[T]$. This is the measure induced by $\mb{P}_{\CR(m)}$ after taking the trace polynomial. It induces further probability measures $\mb{P}_{\SR, \CM_p(m)}$ and $\mb{P}_{\SR,\CM_\CP(m)}$ in the same manner as the measures $\mb{P}_{\CR_p(m)}$ and $\mb{P}_{\CR_\CP(m)}$ are induced by $\mb{P}_{\CR(m)}$. Set
\begin{equation}
\label{eq:Delta for trace polynomials}
\Delta_{\SR,\CP}(m;k) \coloneqq 
\sum_{\substack{\BD \in \CM_\CP \\ \deg{\BD} \le k \\ T \nmid D_p \, \forall p \in \CP}} \max_{\BC \bmod{\BD}} \Bigg\vert \mb{P}_{\SR,\BA \in  \CM_{\CP}(m)} \big(\BA \equiv \BC \bmod{\BD}\big) - \frac{1}{\norm{\BD}_{\CP}} \Bigg\vert.
\end{equation}

\begin{lemma}
\label{lem:Delta correspondence} 
Let $m\ge k\ge1$ be two integers. For any choice of probability measure $\mb{P}_{\CR(m)}$, we have $\Delta_{\SR,\CP}(m;k) \leq \Delta^\SR_\CP(m;k)$. In particular, the conditions of Proposition~\ref{prop:Delta bound} imply
\begin{equation*}
\Delta_{\SR,\CP}\big(m;m/2 + m^{0.88}\big) \ll_r \exp\big(-m^{1/10}\big).
\end{equation*}
\end{lemma}
\begin{proof}
For all polynomials $A,C,D\in \F_p[T]$ such that $\deg(A)=m$ and $\deg(C),\deg(D)\le m$ and $T \nmid D$, we have $A\equiv C\bmod D$ if, and only if, $A^{\SR,m}\equiv C^{\SR,m}\bmod{D^{\SR}}$ by Lemma \ref{lem:trace pol iso of vs}. Thus, if we let $\tilde{\BA}=\BA^{\SR,m}=\BA^{\SR}$, $\tilde{\BC}=\BC^{\SR,m}$ and $\tilde{\BD}=\BD^{\SR}$, then we find
\begin{equation*}
\mb{P}_{\SR,\BA \in \CM_{\CP}(m)} \big(\BA  \equiv \BC \bmod{\BD}\big) = \mb{P}_{\tilde{\BA} \in \CR_\CP(m)}\Big(\tilde{\BA} \equiv \tilde{\BC} \bmod{\tilde{\BD}}\Big).
\end{equation*}
Consequently,
\begin{align*}
\Delta_{\SR,\CP}(m;k) &= \sum_{\substack{\BD \in \CM_\CP \\ \deg{\BD} \le k \\ T \nmid D_p \, \forall p \in \CP}} \max_{\BC \bmod{\BD}} \Bigg\vert \mb{P}_{\tilde{\BA} \in \CR_\CP(m)}\Big(\tilde{\BA} \equiv \tilde{\BC} \bmod{\tilde{\BD}}\Big) - \frac{1}{\norm{\BD}_{\CP}} \Bigg\vert \\ 
&= \sum_{\substack{\tilde{\BD} \in \CR_\CP \\ \deg{\tilde{\BD}} \le 2k \\ T^2+1 \nmid \tilde{D}_p \, \forall p \in \CP}} \max_{\BC \bmod{\tilde{\BD}_{\SR}}}
	\Bigg\vert \mb{P}_{\tilde{\BA} \in \CR_\CP(m)}\Big( \tilde{\BA} \equiv  \tilde{\BC}\bmod{\tilde{\BD}} \Big) - \frac{1}{\norm{\tilde{\BD}}_{\CP}^{1/2}} \Bigg\vert.
\end{align*}
Since the map $\F_\CP[T]/(\tilde{\BD}_\SR) \to R_m(\tilde{\BD})$ sending $\BC \mapsto \tilde{\BC} \bmod{\tilde{\BD}}$ is bijective when $\deg(D_p)\le m$ (cf.~Lemma \ref{lem:reciprocal Euclid} and Definition \ref{def:RmD}), we conclude that the only difference between $\Delta_{\SR,\CP}(m;k)$ and $\Delta^\SR_\CP(m;k)$ (see \eqref{eq:Delta def}) is that in the latter we require $T^2+1 \nmid D_2$ whereas in the former we demand $T^2+1 \nmid D_p$ for all $p \in \CP$. Hence $\Delta_{\SR,\CP}(m;k) \leq \Delta^\SR_\CP(m;k)$. 
\end{proof}

We will show that this suffices to prove that $A_\SR$ is irreducible with high probability. Indeed, as we will argue below, the main result of \cite[\S\S 8--10]{BKK} is roughly the following. Given a probability measure $\mb{P}_{\CM(m)}$ on $\CM(m)$, a finite set of primes $\CP$ and integers $m,k\ge1$, let
\begin{equation}
	\label{eq:Delta general def}
\Delta_{\CP}(m;k) \coloneqq 
	\sum_{\substack{\BD \in \CM_\CP \\ \deg{\BD} \le k \\ T \nmid D_p \, \forall p \in \CP}} \max_{\BC \bmod{\BD}} \Bigg\vert \mb{P}_{\BA \in  \CM_{\CP}(m)} \big(\BA \equiv \BC \bmod{\BD}\big) - \frac{1}{\norm{\BD}_{\CP}} \Bigg\vert.
\end{equation}
In \cite{BKK}, the authors used the measure $\mb{P}_{\CM(m)}$ such that the coefficient $a_j$ of $T^j$ is sampled independently from the rest according to a probability measure $\mu_j$ on $\Z$ satisfying various conditions. They showed for this measure that if $\Delta_{\CP}(m;k)$ is small for $k$ slightly larger than $m/2$, then $A$ does not have any factors of degree in $[m^{1/10},m/2]$ (see Lemma 9.4 and \S10 in \cite{BKK}), which is the analogous result to Proposition \ref{prop:large degree factors}. It turns out that this proof uses very little about the specific structure of the measure $\mb{P}_{\CM(m)}$ with each coefficient sampled independently. We will thus be able to adapt it easily to the measure $\mb{P}_{\SR,\CM(m)}$. To stress this point, we will prove a more general result.

Recall that 
\[
\lambda_0 \coloneqq \frac{1}{4-4\log{2}}.
\]
We then have the following general lemma.

\begin{lemma}[Large degree factors in {\cite[Proposition~2.2]{BKK}}]
	\label{lem:large degree factors in BKK}
	Let $\epsilon \in (0, 1/100]$, let $m \in \Z_{\ge 1}$, and let $\mb{P}_{\CM(m)}$ be a probability measure on the set $\CM(m)$ of monic polynomials of degree $m$ satisfying the following conditions:
	\begin{enumerate}
		\item There is a set $\CP$ of four primes such that $\Delta_{\CP}(m; m/2 + m^{\lambda_0 + \epsilon}) \leq m^{-30}$;
		\item $\mb{P}_{A\in \CM(m)}\big(T^{1+\lceil 4m^{\epsilon/200} \log{m} \rceil} \mid A_p\big) \leq m^{-4}$ for all $p \in \CP$.
	\end{enumerate}
	Then there are constants $c, C > 0$ depending at most on $\epsilon$ such that
	\begin{equation*}
		\mb{P}_{A\in\CM(m)}\Big(A \textup{ has a factor of degree in } \big[m^{1/10}, m/2\big]\Big) \leq C m^{-c}.
	\end{equation*}
\end{lemma}
\begin{proof}
	We make some observations about the use of independence in \cite[\S\S 8--9]{BKK}. First, everything in their \S 8 is already formulated in the required generality. In their \S 9, we encounter the sequence of probability measures $\mu_0, \mu_1, \ldots$ in each of their Lemmas 9.1-9.4. The proofs of their Lemmas 9.1-9.3 do not use that the measure $\mb{P}_{\CM(m)}$ is induced by these, so the lemmas may as well be stated for a general measure $\mb{P}_{\CM(m)}$. This leaves only their Lemma 9.4, a statement having as one of its hypotheses the anti-concentration inequality
	\begin{equation}
		\label{eq:RHS of 9.15}
		\tag{$\star$}
		\sup_{1 \leq j < m} \sum_{a \equiv 0 \bmod{p}} \mu_j(a) \leq 1-\delta \quad \forall p \in \CP.
	\end{equation}
	To be precise, this is the right-hand side of their (9.15). However, the assumption \eqref{eq:RHS of 9.15} only serves the purpose of showing, in their notation,
	\begin{equation}
		\label{eq:9.4 new}
		\mb{P}_{A \in \CM(m)}\Big(T^{1 + \lceil r \delta^{-1} \log{m} \rceil} \mid A_p\Big) \leq m^{-r} \quad \forall p \in \CP;
	\end{equation}
	Indeed, this is established in and around their (9.17) and (9.18). Thus we may take \eqref{eq:9.4 new} as an assumption in \cite[Lemma~9.4]{BKK} instead of \eqref{eq:RHS of 9.15}. The remainder of the proof of their Lemma 9.4 does not make use of the independence of the coefficients of $A$.
	
	Now, observe that condition (1) in Lemma~\ref{lem:large degree factors in BKK} is the same as condition (b) in \cite[Proposition~2.2]{BKK} in the case $\theta = 1/2$. Our condition (2) is implied by the latter half of their condition (c), and in particular suffices for our purposes, which are to be able to apply the revised version of \cite[Lemma~9.4]{BKK} where \eqref{eq:RHS of 9.15} is replaced by \eqref{eq:9.4 new}. Their condition (a) as well as the first half of their condition (c) are only used to rule out factors of small degree by means of \cite[Proposition~2.1]{BKK}, and thus we do not need them here. Hence the result follows by the proof of \cite[Proposition~2.2]{BKK}, given in their \S 10.
\end{proof}

This brings us to the proof of Proposition~\ref{prop:large degree factors}. We first prove the following, more general version.

\begin{proposition}[Large degree factors alternative]
	\label{prop:large degree factors alternative}
	Let $m \in \Z_{\ge 1}$ and $\epsilon \in (0, 1/100]$. Let $\mb{P}_{\CR(m)}$ be any probability measure on $\CR(m)$ satisfying the following:
	\begin{enumerate}
		\item There is a set $\CP$ of four primes such that $\Delta^\SR_\CP(m;m/2 + m^{\lambda_0 + \epsilon}) \leq m^{-30}$.
		\item $\mb{P}_{A\in\CR(m)}\big((T^2+1)^{1+ \lceil 4m^{\epsilon/200} \log{m} \rceil} \mid A_p\big) \leq m^{-4}$ for all $p \in \CP$.
	\end{enumerate}
	Then there are constants $c, C > 0$ depending at most on $\epsilon$ such that
	\begin{equation*}
		\mb{P}_{A \in \CR(m)}\Big(A \textup{ has a reciprocal divisor } D \in \CR(k) \textup{ with } k \in \big[m^{1/10}, m/2\big]\Big) \le Cm^{-c}.
	\end{equation*}
\end{proposition}
\begin{proof}
	In light of Lemma \ref{lem:Delta correspondence}, the conditions in Proposition~\ref{prop:large degree factors alternative} imply those of Lemma~\ref{lem:large degree factors in BKK} for the probability measure $\mb{P}_{\SR,\CM(m)}$. Thus there are constants $c, C > 0$ depending at most on $\epsilon$ such that
	\begin{equation*}
		\mb{P}_{\SR,B \in \CM(m)}\Big(B \textup{ has a factor of degree in } \big[m^{1/10}, m/2\big] \Big) \leq C m^{-c}.
	\end{equation*}
	Making the change of variables $A=B^{\SR}$, and noticing that $B$ has a factor of degree $k$ if and only if $A$ has a reciprocal factor in $\CR(k)$, this is equivalent to
	\begin{equation*}
		\mb{P}_{A \in \CR(m)}\Big(A \textup{ has a reciprocal divisor } D \in \CR(k) \textup{ with } k \in \big[m^{1/10}, m/2\big]\Big) \leq C m^{-c},
	\end{equation*}
	which was to be shown.
\end{proof}

When the coefficients of the reciprocal polynomial $A$ are induced by a sequence of probability measures $\mu_j$ --- as we assume everywhere else in the paper --- we may replace condition (2) of Proposition~\ref{prop:large degree factors alternative} by a (perhaps more appealing) condition on the $\mu_j$. This is entirely similar to the result that \eqref{eq:RHS of 9.15} implies \eqref{eq:9.4 new}. In particular, to prove Proposition~\ref{prop:large degree factors}, it suffices to show that its condition (2) implies condition (2) of Proposition~\ref{prop:large degree factors alternative}. To prove this, we first give a general lemma.

\begin{lemma}
	\label{lem:div_by_fixed_reciprocal}
	Let $\mu_0,\dots,\mu_{m-1}$ be probability measures on $\Z$, and let $\mb{P}_{\CR(m)}$ be the measure defined by \eqref{eq:P on reciprocals}. Let $p$ be a prime and $\delta > 0$, and suppose that 
	\begin{equation}
		\label{eq:concentration mod p}
		\sup_{0 \leq b \leq p-1} \sum_{a \equiv b \bmod{p}} \mu_j(a) \leq 1-\delta \quad \text{ for all } j \in \{0, 1, \ldots, m-1\}.
	\end{equation}
	Then, for every $D\in\CR_p(k)$ with $k\in\Z\cap[1,m]$, we have 
	\[
	\mb{P}_{A\in\CR(m)}\big(D \mid A_p\big) \leq e^{-\delta k}.
	\] 
	In particular, if we let $v = \lceil c \delta^{-1} \log{m} \rceil$ with $c>0$, then 
	\begin{equation*}
		\mb{P}_{\SR,B\in\CM(m)} \big(T^v \mid B_p\big)\le m^{-c}.
	\end{equation*}
\end{lemma}
\begin{proof}
	This is essentially \cite[Lemma~39]{BV} or \cite[Lemma~4.1]{Konyagin}, adapted to the reciprocal setting. For each $j \in [k,m-1]\cap\Z$, pick $a_j$ according to the measure on $\F_p$ induced by $\mu_{j}$ and set $A' = \sum_{j =k}^m a_j(T^{m-j}+T^{m+j}) \in \F_p[T]$. Write $R \coloneqq A' \brmod{D} \in T^{m-k+1}\CR_p^{\ms{sh}}(k-1)$ for the polynomial with $R = A' \bmod{D}$. Conditionally on the choice of $a_k,\dots,a_{m-1}$, we have that 
	\begin{align*}
		D|A_p\quad &\iff\quad a_0 T^m + a_1(T^{m+1}+T^{m-1}) + \ldots + a_{k-1}(T^{m-k+1} + T^{m+k-1}) \equiv - A' \bmod{D} \\
		&\iff\quad a_0 T^m + a_1(T^{m+1}+T^{m-1}) + \ldots + a_{k-1}(T^{m-k+1} + T^{m+k-1})=-R .
	\end{align*}
	We thus see that the coefficients $a_0,\dots,a_{k-1}$ are fixed. Using our assumption \eqref{eq:concentration mod p} proves that the probability that $D|A_p$ conditionally on the choice of $a_k,\dots,a_{m-1}$ is $\le (1-\delta)^k \le e^{-\delta k}$. This proves the first claim. 
	
	To see the second claim, note that
	\[
	\mb{P}_{\SR,B\in\CM(m)} \big(T^v \mid B_p\big)
	=
	\mb{P}_{A\in \CR(m)} \big((T^2+1)^v \mid A_p\big) \le m^{-c},
	\]
	where we made the change of variables  $A=B^{\SR}$. Thus the proof is complete.
\end{proof}

Proposition~\ref{prop:large degree factors} is now a simple corollary of Proposition~\ref{prop:large degree factors alternative}:

\begin{proof}[Proof of Proposition~\ref{prop:large degree factors}] 
It suffices to show that the conditions of Proposition \ref{prop:large degree factors alternative} are met. To this end, it suffices to show that condition (2) of Proposition~\ref{prop:large degree factors} implies condition (2) of Proposition~\ref{prop:large degree factors alternative} for the measure $\mb{P}_{\SR,\CM(m)}$. Indeed, applying Lemma~\ref{lem:div_by_fixed_reciprocal} with $\delta = m^{-\epsilon/200}$ and $c=4$ yields, for any $p \in \CP$, the inequality $\mb{P}_{\SR,B\in\CM(m)}(T^{\lceil 4 m^{\epsilon/200} \log{m} \rceil} \mid B_p) \leq m^{-4}$. This is stronger than condition (2) of Proposition~\ref{prop:large degree factors alternative}, thus completing the proof.
\end{proof}

\section{The hyperoctahedral group and its subgroups}
\label{sec:hyperoctahedral_group}

We have now proven Propositions~\ref{prop:small degree factors}--\ref{prop:large degree factors}, which together yield that $A$ is irreducible with high probability under a wide range of choices of the probability measures $\mu_j$ --- essentially \eqref{eq:A is irreducible}. In the remainder of the paper, we study the Galois group of $A$.

In this section, we describe the generic Galois group of a squarefree reciprocal polynomial $A \in \CR(m)$, which is the \emph{hyperoctahedral group} denoted $\CC_2 \wr \CS_m$. Furthermore, we describe the `large' subgroups of $\CC_2 \wr \CS_m$, in a sense made precise by Lemma~\ref{lem:subgroups projecting onto S_m or A_m} and Theorem~\ref{thm:Luczak--Pyber}. The latter is a \L uczak--Pyber theorem for the group $\CC_2 \wr \CS_m$, which is not used in this article but might be of independent interest.

\subsection*{The hyperoctahedral group}
The zeros of a squarefree reciprocal polynomial $A \in \CR(m)$ are $2m$ distinct algebraic integers that may be labelled as
\begin{equation*}
\alpha_1, \alpha_2, \ldots, \alpha_m \textup{ and } \alpha_{-1} \coloneqq \alpha_1^{-1}, \alpha_{-2} \coloneqq \alpha_2^{-1}, \ldots, \alpha_{-m} \coloneqq \alpha_m^{-1}.
\end{equation*}
 An element $\sigma$ of the Galois group $\CG_A \leq \CS_{2m}$ of $A$ is in particular a field automorphism and thus $\sigma(\alpha_{-j}) = \sigma(\alpha_j)^{-1}$ for all $j$. In other words, the Galois group $\CG_A$ preserves the \emph{block system}
\begin{equation*}
\CB_1 = \{\alpha_1, \alpha_{-1}\}, \ldots, \CB_m = \{\alpha_m, \alpha_{-m}\}
\end{equation*}
in the sense that for all $i$ and all $\sigma$, there exists $j$ such that $\sigma(\CB_i) = \CB_j$. This shows that if $A$ is irreducible, then the action of the permutation group $\CG_A$ on the set of zeros of $A$ is \emph{imprimitive}. (In general, if a transitive permutation group $G \leq \CS_n$ preserves some partitioning of $\{1, 2, \ldots, n\}$ into $k$ sets of equal size for some $2 \leq k \leq n/2$, then $G$ is called \emph{imprimitive}.) 

The group $\CG_A$ naturally lies in the \emph{hyperoctahedral group} $\CC_2 \wr \CS_m$, which is an example of a \emph{permutational wreath product}. In general, if $G$ is an abstract group and $K \leq \CS_m$ is a permutation group, then $G \wr K$ is the semidirect product $G^m \rtimes K$ under the automorphism of $G^m$ induced by $K$ through permutation of the $m$ copies of $G$. Explicitly, it is the group with elements $G^m \times K$ and product
\begin{equation*}
((\epsilon_i)_i, \sigma) \cdot ((\epsilon'_i)_i, \sigma') = ((\epsilon_i \epsilon'_{\sigma^{-1}(i)})_i, \sigma \sigma').
\end{equation*}
The group $\CC_2 \wr \CS_m$ is also known as the Coxeter group of type $B_m$ or as the signed symmetric group. The last name comes from its action on the set $\{-m, \ldots, -1, 1, \ldots, m\}$ of $2m$ signed letters (the set of indices of the $\alpha_j$) given by the explicit formula
\begin{equation}
\label{eq:action_formula}
((\epsilon_i)_i, \sigma).k = \sign(k) \epsilon_{\sigma(\abs{k})}\sigma(\abs{k}).
\end{equation}
This induces the action of $\CG_A$ on the zeros of $A$; it is clear that this action preserves the pairs $\{k, -k\}$. 
Moreover, the formula \eqref{eq:action_formula} gives a procedure to determine the cycle decomposition of an element $((\epsilon_i)_i, \sigma)$, viewed as an element of $\CS_{2m}$. For example, we have $\CC_2 \wr \CS_4 \ni ((-1,1,-1,1),(12)) = (1 \, 2 \, {-1} \, {-2})(3 \, {-3}) \in \CS_8$. 

We refer to \cite{VV} for a proof that a generic reciprocal polynomial of degree $2m$ has Galois group $\CC_2 \wr \CS_m$.

\subsection*{Subgroups of the hyperoctahedral group}

For $m \ge 2$, consider the following proper subgroups of $\CC_2 \wr \CS_m$: 
\begin{enumerate}
    \item The index-$2$ subgroup 
    \begin{equation}
    \label{eq:G1}
			G_1 \coloneqq \{((\epsilon_i)_i, \sigma) \in \CC_2 \wr \CS_m ~:~ \prod_i \epsilon_i = 1\} = (\CC_2 \wr \CS_m) \cap \CA_{2m};
    \end{equation}
    \item The index-$2$ subgroup
    \begin{equation}
    \label{eq:G2}
    	G_2 \coloneqq \{((\epsilon_i)_i, \sigma) \in \CC_2 \wr \CS_m ~:~ \sign(\sigma) \prod_i \epsilon_i = 1\}
    \end{equation}
     which, as a set, is $((\CC_2^m \times \CA_m) \cap \CA_{2m}) \cup ((\CC_2^m \times (\CS_m \setminus \CA_m)) \cap (\CS_{2m} \setminus \CA_{2m}))$; 
    \item The index-$2$ subgroup 
    \begin{equation}
    \label{eq:G3}
			G_3 \coloneqq \{((\epsilon_i)_i, \sigma) \in \CC_2 \wr \CS_m ~:~ \sign(\sigma) = 1\} = \CC_2 \wr \CA_m;
    \end{equation}
    \item The index-$4$ subgroup arising as the intersection of $G_1$, $G_2$, and $G_3$,
    \begin{equation}
    \label{eq:G4}
    	G_4 \coloneqq \{((\epsilon_i)_i, \sigma) \in \CC_2 \wr \CS_m ~:~ \sign(\sigma) = \prod_i \epsilon_i = 1\} = (\CC_2 \wr \CA_m) \cap \CA_{2m}.
    \end{equation}
\end{enumerate}

\begin{remark}
The groups $G_2$ and $G_1$ are not permutation isomorphic, but they \emph{are} isomorphic when $m\ge 3$ is odd: the map is $G_2 \to G_1, ((\epsilon_i)_i, \sigma) \mapsto ((\epsilon_i \sign(\sigma))_i, \sigma)$.
\end{remark}

The group $\CC_2 \wr \CS_m$ comes with the projection map 
\begin{align*}
\proj \colon \CC_2 \wr \CS_m& \to \CS_m, \\
((\epsilon_i)_i, \sigma)&\mapsto \sigma.
\end{align*}
The purpose of this subsection is to show that any proper subgroup $H \leq \CC_2 \wr \CS_m$ that projects onto $\CS_m$ or $\CA_m$ is either one of the `large' groups $G_1$, $G_2$, $G_3$ and $G_4$, or is conjugate to a subgroup of
\begin{equation}
\label{eq:G5}
G_5 \coloneqq \{((\epsilon_i)_i, \sigma) \in \CC_2 \wr \CS_m ~:~ \epsilon_1 = \epsilon_2 = \ldots = \epsilon_m\} \cong \CC_2 \times \CS_m.
\end{equation}

The reason to restrict to subgroups projecting onto $\CA_m$ or $\CS_m$ is because we prove in \S \ref{sec:galois_group_of_trace_pol} that $\proj(\CG_A)$ is with high probability either $\CA_m$ or $\CS_m$.

We start with the following standard result on the subgroups of a semidirect product.

\begin{lemma}
\label{lem:subgroups_semidirect}
Let $U$ be an abelian group and $K \leq \CS_m$. Consider the semidirect product $U^m \rtimes K$, where $K$ acts on the $m$ copies of $U$ by permutation. For any subgroup $H \leq U^m \rtimes K$ with $\proj(H) = K$, the set $V = U^m \cap H =  \{((\epsilon_i)_i, \id) \in H\}$ is a $K$-invariant subgroup of $U^m$. Furthermore, given a $K$-invariant subgroup $V$ of $U^m$, the set of subgroups $H \leq U^m \rtimes K$ with $U^m \cap H = V$ and $\proj(H) = K$ is, up to conjugation by $U^m$, in bijective correspondence with the first cohomology group $\CH^1(K, U^m/V)$; given a cocycle $\xi \in Z^1(K, U^m/V)$, the corresponding subgroup is $H_{\xi} = \{(u, \sigma) : \sigma \in K, \, uV = \xi(\sigma)\}$.
\end{lemma}
\begin{proof}
Since $H$ is closed under conjugation by itself and $U$ is abelian, the group $V$ contains the element
\begin{equation*}
((\epsilon_i')_i, \sigma)\cdot ((\epsilon_i)_i, \id) \cdot ((\epsilon_i')_i, \sigma)^{-1} = ((\epsilon_{\sigma^{-1}(i)})_i, \id)
\end{equation*} 
for any $((\epsilon_i)_i, \id), ((\epsilon_i')_i, \sigma) \in H$. Hence $V$ is a $K$-invariant subgroup of $U^m$. This proves the first claim. For the second claim, we refer the reader to \cite[Section 17]{Aschbacher} or \cite[Lemma~3.3]{ABD}.
\end{proof}

To apply Lemma~\ref{lem:subgroups_semidirect}, we study the invariant subgroups of $\CC_2^m$ by the permutation actions of $\CA_m$ and $\CS_m$. The following lemma is probably well-known.

\begin{lemma}
\label{lem:invariantsubspaceAn}
For $m \in \Z_{\geq 5}$, the only invariant subgroups $V \subset \CC_2^m$ by the permutation action of $\CA_m$ are
\begin{itemize} 
\item $V_0=\{(1,\dots,1)\}$,
\item $V_1=\{(-1,\dots,-1),(1,\dots,1)\}$,
\item $V_{m-1}=\{\epsilon\in \CC_2^m: \prod_i \epsilon_i = 1\}$,
\item $V_m=\CC_2^m$.
\end{itemize}
The same is true if $\CA_m$ is replaced by $\CS_m$.
\end{lemma}
\begin{proof}
For $\CS_m$, a proof of this result is given in \cite[Lemma~4.2]{BBF}. The proof for $\CA_m$ we give here proceeds in similar vein. If an invariant subgroup $V \subset \CC_2^m$ is different from $V_0$ and $V_1$, then there exists $\epsilon\in V$ and distinct $i$, $j$, and $k$ such that $\epsilon_i = \epsilon_j = -\epsilon_k$. Consider $\sigma = (ijk) \in \CA_m$. Then $\epsilon' \coloneqq \epsilon \cdot \sigma(\epsilon)$ has all entries equal to $1$ except for the $i$-th and the $k$-th. Permuting $\epsilon'$ by permutations that are products of two transpositions, we can form all tuples in $\CC_2^m$ with only two entries equal to $-1$. Such tuples generate $V_{m-1}$. If $V \neq V_{m-1}$, then $V=V_m$ since $V_{m-1}$ is of codimension $1$ in $V_m$.
\end{proof}

\begin{lemma}
	\label{lem:subgroups projecting onto S_m or A_m}
	Let $m \geq 6$. Let $H \le \CC_2 \wr \CS_m$ be a proper subgroup that projects onto $K \in \{\CS_m, \CA_m\}$.
	\begin{itemize}
		\item If $[\CC_2 \wr \CS_m: H] < 2^{m-1}$ and $K = \CS_m$, then $H \in \{G_1, G_2\}$. 
		\item If $[\CC_2 \wr \CS_m: H] < 2^{m-1}$ and $K = \CA_m$, then $H \in \{G_3, G_4\}$.
		\item If $[\CC_2 \wr \CS_m: H] \ge 2^{m-1}$, then $H$ is $\CC_2^m$-conjugate to a subgroup of $V_1 \times \CS_m = G_5$.
	\end{itemize}
\end{lemma}
\begin{proof}
	We will apply Lemma~\ref{lem:subgroups_semidirect} with $U = \CC_2$, $K \in \{\CS_m, \CA_m\}$ and $V \in \{V_0, V_1, V_{m-1}, V_m\}$. 
	Start by observing that $\proj \colon H \to K$ is surjective and has kernel $V$, so that $|H| = |V| |K|$ and
	\begin{equation}
		\label{eq:index-H-from-V}
		[\CC_2\wr\CS_m:H]
		= \frac{2^m m!}{|V| |K|}.
	\end{equation}
	In particular, we have $[\CC_2 \wr \CS_m : H] \geq 2^{m-1}$ if and only if $V \in \{V_0, V_1\}$.
	
	First, consider the case that $K = \CS_m$ and $ [\CC_2\wr\CS_m:H]<2^{m-1}$, so that $V \in \{V_{m-1}, V_m\}$. The classification of subgroups of $\CC_2\wr\CS_m$ that project onto $\CS_m$ in \cite[Theorem~3.4 and Remark~3.6]{ABD} shows that, up to $\CC_2^m$-conjugacy, the only proper such subgroups with $V\in\{V_{m-1},V_m\}$ are $G_1$ and $G_2$ (which correspond to the $G_1$ and $G_2$ in \cite{ABD}). Since $G_1$ and $G_2$ are normal in $\CC_2\wr\CS_m$, it follows that $H \in \{G_1, G_2\}$.
	
	Next, we consider the case that $K = \CA_m$ and $ [\CC_2\wr\CS_m:H]<2^{m-1}$, so that again $V \in \{V_{m-1}, V_m\}$. If $V = V_{m-1}$, then
	\begin{equation*}
		\CH^1(\CA_m, \CC_2^m/V_{m-1}) = \mr{Hom}(\CA_m, \CC_2^m/V_{m-1}) = \mr{Hom}(\CA_m, \CC_2),
	\end{equation*}
	since the action of $\CA_m$ on $\CC_2^m/V_{m-1}$ is trivial and $\CC_2^m/V_{m-1} \cong \CC_2$. Hence $\CH^1(\CA_m, \CC_2^m/V_{m-1})$ is trivial, because $\CA_m$ is simple by the assumption $m \geq 6$ (so any homomorphism must have kernel $\CA_m$ or $1$, and only the former is possible in this case). The corresponding map $\xi$ is hence the map sending everything to $1$, which gives, up to conjugacy, the group $H = \{((\epsilon_i)_i, \sigma) \in \CC_2 \wr \CA_m : \prod \epsilon_i = 1\} = G_4$. For $V=V_m$, we find
	\begin{equation*}
		\CH^1(\CA_m, \CC_2^m/V_{m}) = \CH^1(\CA_m, 1) = 1
	\end{equation*}
	which gives, up to conjugacy, the full group $H = \CC_2 \wr \CA_m = G_3$. Since both $G_3$ and $G_4$ are normal in $\CC_2 \wr \CS_m$, the second assertion of the lemma statement follows as well.
	
	We now prove the third assertion, which requires a bit more work. Suppose that $[\CC_2\wr\CS_m:H]\ge2^{m-1}$, so that $V\in\{V_0,V_1\}$. When $K=\CS_m$, the result again follows from \cite[Theorem~3.4 and Remark~3.6]{ABD}. It remains to consider $K=\CA_m$. We treat the two possibilities $V=V_0$ and $V=V_1$ separately.
	
	For $V = V_0$, we apply Shapiro's lemma as in \cite[Theorem~3.4]{ABD}. First observe that in this action, the permutation module $\CC_2^m$ is isomorphic to $\mr{Ind}_{\CA_{m-1}}^{\CA_m}\CC_2$. Indeed, $\CA_m$ acts transitively on $\{1,\ldots,m\}$ with point stabilizer $\CA_{m-1}$. Now $\mr{Ind}_{\CA_{m-1}}^{\CA_m}\CC_2 \cong \mr{Coind}_{\CA_{m-1}}^{\CA_m}\CC_2$ by \cite[Chapter~III, Proposition~5.9]{Brown}. Hence we may apply Shapiro's lemma \cite[Chapter~III, Proposition~6.2]{Brown} to obtain
	\[
	\CH^1(\CA_m,\CC_2^m) \cong \CH^1(\CA_{m-1},\CC_2) = \mr{Hom}(\CA_{m-1},\CC_2) = 1,
	\]
	where the last step follows because $\CA_{m-1}$ is simple for $m\ge6$. An application of Lemma~\ref{lem:subgroups_semidirect} yields that $H$ is $\CC_2^m$-conjugate to $V_0\times\CA_m$, which is a subgroup of $G_5$.
	
	Finally, suppose that $V=V_1$. We claim that
	\begin{equation}
		\label{eq:H1-Am-C2m-mod-V1}
		\CH^1(\CA_m,\CC_2^m/V_1)=1.
	\end{equation}
	It then follows with a similar argument as above that $H$ is $\CC_2^m$-conjugate to $V_1 \times \CA_m$, which is also a subgroup of $G_5$, thus proving the lemma.
	
	It thus remains to prove \eqref{eq:H1-Am-C2m-mod-V1}, for which we will apply inflation-restriction. First note that
	\begin{equation}
		\label{eq:H1-Sm-C2m-mod-V1}
		\CH^1(\CS_m,\CC_2^m/V_1)=1,
	\end{equation}
	which again follows from \cite[Remark~3.6]{ABD}. Next, we show that
	\begin{equation}
		\label{eq:Am-invariants-C2m-mod-V1}
		(\CC_2^m/V_1)^{\CA_m}=1.
	\end{equation}
	Indeed, suppose that $\epsilon V_1\in\CC_2^m/V_1$ is fixed by $\CA_m$. Then, for every $\sigma \in\CA_m$, we have $\epsilon^{-1}\sigma(\epsilon)\in V_1$. If $\sigma$ is a $3$-cycle, the sign vector $\epsilon^{-1}\sigma(\epsilon)$ has all coordinates equal to $1$ outside the three coordinates moved by $\sigma$. Since $m\ge6$, it cannot equal $(-1,\ldots,-1)$. Hence
	\[
	\epsilon^{-1}\sigma(\epsilon)=(1,\ldots,1),
	\]
	so $\sigma(\epsilon)=\epsilon$. Since $\CA_m$ is generated by $3$-cycles, the vector $\epsilon$ is fixed by $\CA_m$. As $\CA_m$ acts transitively on the coordinates, all coordinates of $\epsilon$ are equal, and therefore $\epsilon\in V_1$.
	
	Now apply the inflation-restriction five-term exact sequence \cite[Proposition~1.6.7]{NSW} to
	\[
	1 \longrightarrow \CA_m \longrightarrow \CS_m \longrightarrow \CC_2 \longrightarrow 1
	\]
	with coefficients in $\CC_2^m/V_1$. In particular, this gives the sequence
	\begin{equation*}
	\begin{aligned}
		1 \longrightarrow\;&
		\CH^1\big(\CC_2,(\CC_2^m/V_1)^{\CA_m}\big)
		\longrightarrow
		\CH^1(\CS_m,\CC_2^m/V_1)
		\\
		& \qquad \qquad \qquad \qquad \qquad \qquad \longrightarrow\;
		\CH^1(\CA_m,\CC_2^m/V_1)^{\CC_2}
		\longrightarrow
		\CH^2\big(\CC_2,(\CC_2^m/V_1)^{\CA_m}\big).
	\end{aligned}
	\end{equation*}
	By \eqref{eq:Am-invariants-C2m-mod-V1}, the first and last groups are trivial, while \eqref{eq:H1-Sm-C2m-mod-V1} shows that the second group is trivial. Consequently,
	\begin{equation}
		\label{eq:H1-Am-C2m-mod-V1-C2}
		\CH^1(\CA_m,\CC_2^m/V_1)^{\CC_2}=1.
	\end{equation}
	Since $\CC_2^m/V_1$ is an elementary abelian $2$-group, we also have that $\CH^1(\CA_m,\CC_2^m/V_1)$ is an elementary abelian $2$-group. Suppose that $1 \neq z \in \CH^1(\CA_m,\CC_2^m/V_1)$. Let $s$ denote the nontrivial element of $\CC_2$. Then $s$ does not fix $z$ by \eqref{eq:H1-Am-C2m-mod-V1-C2}. Thus $y\coloneqq z^{-1}s(z) \in  \CH^1(\CA_m,\CC_2^m/V_1)$ is nontrivial. Therefore $s(y) = s(z)^{-1} s^2(z) = s(z)^{-1} z$ since $s$ has order two; and $s(z)^{-1} z = z s(z)^{-1}$ since $\CH^1(\CA_m,\CC_2^m/V_1)$ is abelian; and, finally, we have $z = z^{-1}$ and $s(z)^{-1} = s(z)$, so that $s(y) = y$, since each element of $\CH^1(\CA_m,\CC_2^m/V_1)$ has order at most $2$. But this contradicts \eqref{eq:H1-Am-C2m-mod-V1-C2}.
\end{proof}

The next result is a consequence of the preceding discussion and the classical \L uczak--Pyber theorem; even if it is not needed anywhere else in the paper, we include it as we could not find it recorded anywhere else.
\begin{theorem}[\L uczak--Pyber for $\CC_2 \wr \CS_m$]
\label{thm:Luczak--Pyber}
Let $m \ge 6$. Denote by $\CT_m$ the union of all proper subgroups $H \leq \CC_2 \wr \CS_m$ with $H \neq G_1, G_2, G_3, G_4$ that are transitive on the set $\{-m, \ldots, -1, 1, \ldots, m\}$ by the action given in \eqref{eq:action_formula}. Then there is an absolute constant $c>0$ such that
\begin{equation*}
\# \CT_m/\#(\CC_2 \wr \CS_m) \ll m^{-c}.
\end{equation*}
\end{theorem}
\begin{proof}
Let $H$ be a subgroup of $\CC_2 \wr \CS_m$. Then $\proj(H)$ is a subgroup of $\CS_m$ and comes with the action on the set of pairs $\CY = \{\{j, -j\} : j \in \Z \cap [1,m]\}$ given by $\sigma \{j, -j\} \coloneqq \{\sigma(j), -\sigma(j)\}$ for any $\sigma \in \proj(H)$. This coincides with the action of $H$ given by \eqref{eq:action_formula} on $\CY$. In particular, if the action of $H$ on $\{-m, \ldots, -1, 1, \ldots, m\}$ given by \eqref{eq:action_formula} is transitive, then the action of $H$ on $\CY$ is transitive, so the action of $\proj(H)$ on $\CY$ is transitive as well --- i.e., $\proj(H)$ is a transitive subgroup of $\CS_m$. 

Now write
\begin{equation*}
\CT_m = \CT_m' \cup \CT_m''
\end{equation*}
where $\CT_m'$ is the union of all proper transitive $H \leq \CC_2 \wr \CS_m$ with $\proj(H) \not \in \{\CS_m, \CA_m\}$, and $\CT_m''$ is the union of all proper transitive $H \leq \CC_2 \wr \CS_m$ with $\proj(H) \in \{\CS_m, \CA_m\}$ and $H \neq G_1, G_2, G_3, G_4$. Then the classical Łuczak--Pyber theorem \cite{LP} implies the existence of an absolute constant $c_0>0$ such that $\proj(\CT_m')$ is of size $\ll m!/m^{c_0}$. Since the map $\proj \colon \CC_2 \wr \CS_m \to \CS_m$ is $2^m$-to-one, this implies $\CT_m'$ is of size $\ll 2^m m!/m^{c_0}$. 

To bound the size of $\CT_m''$, observe that by Lemma~\ref{lem:subgroups projecting onto S_m or A_m} it satisfies
\begin{equation*}
	\CT_m'' \subset \CU_m \coloneqq \bigcup_{b \in \CC_2^m} b G_5 b^{-1}.
\end{equation*}
Let $\delta \in \{\pm 1\}$. Conjugating $(\delta (1, \ldots, 1), \sigma) \in G_5$ by some $b = ((b_i)_i, \id) \in \CC_2 \wr \CS_m$ yields
\begin{equation}
	\label{eq:conjugate}
	((b_i)_i, \id) \cdot (\delta (1, \ldots, 1), \sigma) \cdot ((b_i)_i, \id)^{-1} = (\delta (b_i b_{\sigma^{-1}(i)}^{-1})_i, \sigma).
\end{equation}
We now need to bound the number of elements of $\CC_2 \wr \CS_m$ that can be realised by conjugation as in \eqref{eq:conjugate}, or equivalently, the number of pairs $(\epsilon, \sigma) \in \CC_2^m \times \CS_m$ such that $\epsilon = \delta (b_i b_{\sigma^{-1}(i)}^{-1})_i$. That is at most twice the number of pairs $(\epsilon, \sigma)$ with $\epsilon = (b_i b_{\sigma^{-1}(i)}^{-1})_i$ for some $b \in \CC_2^m$. To bound the latter, take such $b$, fix $\sigma$, let $C$ be a cycle in the disjoint cycle decomposition of $\sigma$, and consider the product
\begin{equation*}
	\prod_{i \in C} b_i b_{\sigma^{-1}(i)}^{-1}.
\end{equation*}
Since $\sigma$ permutes the elements within the cycle $C$, we have
\begin{equation*}
	\prod_{i \in C}  b_{\sigma^{-1}(i)}^{-1} = \prod_{i \in C} b_i^{-1} = \prod_{i \in C} b_i.
\end{equation*}
Therefore 
\begin{equation}
	\label{eq:prodei}
	\prod_{i \in C} \epsilon_i = \prod_{i \in C} b_i b_{\sigma^{-1}(i)}^{-1} = 1.
\end{equation}
Thus the cardinality of $\CU_m$ is bounded by twice the number of $(\epsilon, \sigma) \in \CC_2 \wr \CS_m$ with $\prod_{i \in C} \epsilon_i = 1$ for every cycle $C$ (in the disjoint cycle decomposition) of $\sigma$. Fixing $\sigma$, denote by $c(\sigma)$ the number of cycles of $\sigma$. For any cycle $C = (c_1 \, c_2 \, \ldots \, c_k)$, the condition \eqref{eq:prodei} implies that $\epsilon_{c_k}$ is fixed once arbitrary $\epsilon_{c_1}, \ldots, \epsilon_{c_{k-1}}$ are given. Hence there are $2^{m - c(\sigma)}$ possible choices for $\epsilon$. Summing this expression over all $\sigma$, and plugging back in the extra factor $2$ coming from the global sign $\delta$, yields the bound
\begin{equation}
	\label{eq:CUm bound}
\# \CU_m \leq 2 \cdot 2^m \cdot \sum_{\sigma \in \CS_m} 2^{-c(\sigma)} = 2^{m+1} \sum_{k=0}^m \left[ {m \atop k} \right] 2^{-k}
\end{equation}
where $\left[ {m \atop k} \right]$ denotes the number of $\sigma \in \CS_m$ with exactly $k$ disjoint cycles. The $\left[ {m \atop k} \right]$ are known as the \emph{unsigned Stirling numbers of the first kind}. Their ordinary generating function satisfies
\begin{equation*}
	\sum_{k=0}^m \left[ {m \atop k} \right] x^{k} = x(x+1)\cdot \ldots \cdot (x+m-1),
\end{equation*}
see \cite[\S A.8]{FS}. Combining this with \eqref{eq:CUm bound} yields
\begin{equation*}
	\# \CT_m'' \leq \# \CU_m \leq 2^{m+1} \frac{1 \cdot 3 \cdot \ldots \cdot (2m-1)}{2^{m}} \ll (2m-1)!! \ll \frac{2^m m!}{\sqrt{m}}
\end{equation*}
where we used Stirling's asymptotic to obtain the last bound.
Taking $c = \min \{c_0, 1/2\}$, this concludes the proof.
\end{proof}

\section{The Galois group of $A_{\SR}$}
\label{sec:galois_group_of_trace_pol}

For squarefree $A$, the projection $\proj(\CG_A)$ conveniently equals the Galois group of $A_\SR$, as we show in the following lemma.

\begin{lemma}
\label{lem:Galois_gp_trace_pol}
Suppose $A$ is squarefree. Then $\CG_{A_\SR} = \proj(\CG_A)$.
\end{lemma}
\begin{proof}
Let $B = A_{\SR}$. The trace polynomial $B$ has zeros $\alpha_1+\alpha_{-1}, \ldots, \alpha_m+\alpha_{-m}$, which are distinct because $A$ is squarefree and the equation $u+1/u = v+1/v$ has only the solutions $u=v$ and $u=1/v$. Hence $B$ is squarefree as well. Denote the splitting fields of $B$ and $A$ by $K_B$ and $K_A$. These are Galois extensions of $\Q$ with $K_B \subset K_A$. A standard result in Galois theory now says that any automorphism $\sigma \in \CG_B$ can be extended to an element of $\CG_A$, and conversely that the restriction of an element of $\CG_A$ to $K_B$ gives an element of $\CG_B$.
Suppose $\tau \coloneqq ((\epsilon_i)_i, \sigma) \in \CG_A$. Then the restriction $\tau \vert_{K_B} \in \CG_B$ maps $\alpha_k + \alpha_{-k}$ to $\alpha_{\sigma(k)} + \alpha_{-\sigma(k)}$. In particular, since the zeros of $B$ are invariant under the maps $\alpha_j \mapsto \alpha_{-j}$, the action of $\tau$ on the zeros of $B$ only depends on $\sigma$. So we may identify $\tau \vert_{K_B}$ with $\proj(\tau) = \sigma$. Conversely, if $\sigma \in \CG_B$, then any extension of $\sigma$ to $\CG_A$ is of the form $((\epsilon_i)_i, \sigma)$ for some $(\epsilon_i)_i \in \CC_2^m$.
\end{proof}

In the remainder of this section, we show that $\CG_{A_\SR}$ is $\CA_m$ or $\CS_m$ with high probability. Here we use the results from \cite{BKK}, but since the coefficients of $A_{\SR}$ are not independent random variables, we have to make a careful analysis of the use of independence there, just as we did in \S \ref{sec:anatomy}. We start with the following lemma. For this, we recall the notation $\Delta_p$ defined in \eqref{eq:Delta general def} (with the usual convention that $\Delta_p=\Delta_{\{p\}}$). 

\begin{lemma}
\label{lem:replacement_2.4}
Let $\mb{P}_{\CM(m)}$ be any probability measure on the set $\CM(m)$ of monic polynomials of degree $m$. Fix a prime $p$.
Sample $B$ according to $\mb{P}_{\CM(m)}$ and let $B_p \coloneqq (B \bmod{p}) \in \F_p[T]$ be its reduction modulo $p$.
Furthermore, fix a real number $\epsilon>0$ with the properties
\begin{enumerate}
  \item $\Delta_p\big(m, m/2 + m^{\lambda_0+\epsilon}\big) \le m^{-10}$,
  \item For $v > (\log{m})^3$, we have $\mb{P}_{B\in\CM(m)}\big(T^v\mid B_p\big) \leq m^{-1}$.
\end{enumerate}
Then there are constants $c, C > 0$ depending at most on $\epsilon$ such that
\begin{equation*}
\mb{P}_{B\in \CM(m)}\Big(\CG_B \not\in \{\CA_m, \CS_m\} \textup{ and } B \textup{ is irreducible}\Big) \le Cm^{-c}.
\end{equation*}
\end{lemma}
\begin{proof}
This is an adaptation of \cite[Proposition~2.4]{BKK}. Denote their sequence of probability measures by $\mu_0', \ldots, \mu_{m-1}'$. The difference between that proposition and our lemma is that we replaced the condition
\begin{equation}
\label{eq:mild anti-concentration mod p}
\sup_{1 \le j < m} \sum_{a \equiv 0 \bmod{p}} \mu_j'(a) \le 1-\frac{1}{(\log{m})^2}
\end{equation}
by our condition (2) and dropped the condition that $\mb{P}_{\CM(m)}$ is induced by a sequence of measures $\mu_0', \ldots, \mu_{m-1}'$ --- that is, that the coefficients of the polynomial $B \in \CM(m)$ are \emph{independent} random variables. To show that the conclusion of the lemma continues to hold under this weakened condition (2), we analyze the use of \eqref{eq:mild anti-concentration mod p} in \cite{BKK}. Note that this is the right-hand side of their Assumption (11.3), which is used in the following places:
\begin{itemize}
  \item Proof of Lemma~11.1(b): Here, the assumption is only used to ensure that the conditions of \cite[Lemma~9.4]{BKK} are met with $\epsilon = 1/100$, $\delta=1/(\log{n})^2$, $\theta=1/2$, and $\CP = \{p\}$. However, in the proof of Lemma~\ref{lem:large degree factors in BKK}, we showed that the conditions of \cite[Lemma~9.4]{BKK} can be weakened to a condition that holds under our condition (2).
  \item Proof of Lemma~11.5: The assumption is used only to prove our condition (2).
  \item Statement of Lemma~11.6: The assumption is also included here. But Lemma~11.6 follows immediately from combining Proposition~12.1, which does not use this assumption, with Lemma~11.1, which we already discussed.
\end{itemize}
This finishes the proof.
\end{proof}

Now we adapt Lemma~\ref{lem:replacement_2.4} to the setting of reciprocal polynomials.

\begin{lemma}
\label{lem:Galois gp of trace pol}
Consider probability measures $\mu_0, \mu_1, \ldots, \mu_{m-1}$ on the integers. Fix a prime $p$ and a real number $\epsilon>0$ with the properties
\begin{enumerate}
  \item $\Delta^{\SR}_p(m, m/2 + m^{\lambda_0+\epsilon}) \le m^{-10}$,
  \item $\sup_{0 \le j < m} \sup_{0 \le b < p} \sum_{a \equiv b \bmod{p}} \mu_j(a) \le 1-\frac{1}{(\log{m})^2}$.
\end{enumerate}
Then there are constants $c, C > 0$ depending at most on $\epsilon$ such that
\begin{equation*}
\mb{P}_{A \in \CR(m)}\Big(\CG_{A_{\SR}} \not\in \{\CA_m, \CS_m\} \textup{ and } A \textup{ is irreducible}\Big) \le Cm^{-c}.
\end{equation*}
\end{lemma}
\begin{proof}
Our aim is to show this using Lemma~\ref{lem:replacement_2.4} for the probability measure $\mb{P}_{\SR,\CM(m)}$. In light of Lemma \ref{lem:Delta correspondence}, the first condition in Lemma~\ref{lem:replacement_2.4} follows from our condition (1) here. In addition, Lemma~\ref{lem:div_by_fixed_reciprocal}, which we may apply with $\delta = (\log{m})^{-2}$ and $c=1$, yields $\mb{P}_{\SR,B\in \CM(m)}(T^{\lceil (\log{m})^3 \rceil} \text{ divides } B_p) \leq 1/m$. Lastly, if $A$ is irreducible then $B=A_{\SR}$ is irreducible, so the conclusion of the lemma follows from Lemma~\ref{lem:replacement_2.4}.
\end{proof}

\section{The Galois group of $A$}
\label{sec:galois_group_of_A}

This section contains the proof of Proposition~\ref{prop:Galois}, which roughly states that if the coefficients of $A$ are selected according to measures $\mu_j$ that are not too pathological, then
\[
\mb{P}_{A \in \CR(m)}\Big( \CG_A \not \in \{G_2, G_3, \CC_2 \wr \CS_m\} \text{ and } A \text{ is irreducible} \Big) \leq Cm^{-c}
\]
where $c,C>0$ are appropriate constants. More specifically, we assume there exist parameters $B\ge1$, $\alpha\in(0,1/2)$ and $\epsilon\in(0,1/100]$ such that:
\begin{enumerate}
	\item (bounded support) each $\mu_j$ is supported on $[-H,H]$ with $H \le Bm^{1/2-\alpha}$;
	\item (mod-$p$ approximate equidistribution) $\Delta^{\SR}_p(m, m/2 + m^{\lambda_0+\epsilon}) \le m^{-10}$;
	\item (mod-$p$ anti-concentration) $\sup_{0 \le j < m} \norm{\bar{\mu}_{j,p}}_{\infty} \le 1-(\log{m})^{-2}$.
\end{enumerate}
The constants $C$ and $c$ are allowed to depend on $B,\alpha,\epsilon$.

Now, let 
\[
\CI(m) = \Big\{ A\in \CR(m) \,: \, A\ \text{is irreducible}\Big\} 
\]
and define the events 
\begin{align*}
\CE_1 &= \Big\{ A\in \CI(m) : \CG_A\in \{G_2,G_3, C_2\wr \CS_m\} \Big\},\\
\CE_2 &= \Big\{ A\in \CI(m) : \CG_A\in \{G_1,G_4\} \Big\},\\
\CE_3 &=\Big\{ A\in \CI(m) : \CG_A\ \mbox{is $C_2^m$-conjugate to a subgroup of $G_5$}\Big\}.
\end{align*}
In this notation, Proposition~\ref{prop:Galois} is equivalent to proving that
\begin{equation}
		\label{eq: prop Galois target}
		\mb{P}_{\CR(m)}\big( \CI(m)\setminus \CE_1  \big) \leq Cm^{-c} 
\end{equation}
assuming that the measures $\mu_0,\mu_1,\dots,\mu_{m-1}$ satisfy conditions (1)-(3) above.

Combining the results of \S \ref{sec:hyperoctahedral_group} and \S \ref{sec:galois_group_of_trace_pol} (see Lemmas \ref{lem:Galois_gp_trace_pol}, \ref{lem:subgroups projecting onto S_m or A_m} and \ref{lem:Galois gp of trace pol}), which are applicable in view of conditions (2) and (3), we find that
\begin{equation}
	\label{eq: prop Galois step 1}
	\mb{P}_{\CR(m)}\Big( \CI(m)\setminus \big( \CE_1\cup\CE_2\cup\CE_3\big)  \Big) \ll_\epsilon m^{-c_1}
\end{equation}
for some $c_1=c_1(\epsilon)>0$. 

In addition, note that $G_4 \leq G_1 \leq \CA_{2m}$. Together with Lemma~\ref{lem:square-disc-alt-gp}, this implies that 
\[
\CE_2  \subseteq \{\Delta(A) \text{ is a nonzero square}\} .
\]
Thus, Proposition~\ref{prop:square_discriminant}, in conjunction with conditions (1) and (3), implies that
\begin{equation}
	\label{eq: prop Galois step 2}
	\mb{P}_{\CR(m)}\big(\CE_2\big)  \ll_\alpha B m^{-\alpha/2} .
\end{equation}

Combining \eqref{eq: prop Galois step 1} and \eqref{eq: prop Galois step 2}, we have reduced the proof of Proposition~\ref{prop:Galois} to showing that 
\begin{equation}
	\label{eq: prop Galois step 3}
\mb{P}_{\CR(m)} (\CE_3)  \ll m^{-c_2} 
\end{equation}
for some appropriate $c_2>0$, under the assumption of conditions (1)-(3). We prove this estimate in the next subsection (cf.~Lemma \ref{lem:Galois_C2xSm_2}).

\subsection*{Proof of \eqref{eq: prop Galois step 3}}
The idea is that the group $G_5$ is very small compared to $\CC_2 \wr \CS_m$ --- of index $2^{m-1}$ in the latter, to be precise; one should thus expect it, or one of its $C_2^m$-conjugates, to be unlikely to be the Galois group of $A$. The next lemma forms the basis of our proof that this intuition is correct.

\begin{lemma}
\label{lem:Galois_C2xSm_1}
Let $m\in\Z_{\ge1}$ and let $A \in \CR(m)$ be squarefree. In addition, let $p$ be a prime, let $d\in\Z_{\ge1}$ and let $I\in\CR_p(2d)$ be irreducible over $\F_p$ and such that $I \mid A_p$ but $I^2 \nmid A_p$. Then $\CG_A$ is not $\CC_2^m$-conjugate to a subgroup of $G_5$.
\end{lemma}
\begin{proof}
The proof is inspired by \cite[Lemma~11.3]{BKK}. Let $F$ be the splitting field of $A$ over $\Q$ and write $\CO_F$ for the ring of integers in $F$. Consider a prime ideal $\mfp$ of $\CO_F$ lying over $p$. 
For $x \in \CO_F$, denote by $\bar{x}$ its reduction mod $\mfp$ to the field $\CO_F/\mfp$. 
Recall that the Frobenius automorphism $\phi_p \colon \CO_F/\mfp \to \CO_F/\mfp$, defined as $\phi_p(\bar{x}) = \bar{x}^p$, lifts to an element $\phi \in \CG_A$.

Write $\bar{\alpha}_1, \ldots, \bar{\alpha}_{2d}, \bar{\alpha}_{-1} = (\bar{\alpha}_1)^{-1},\ldots, \bar{\alpha}_{-2d}=(\bar{\alpha}_{2d})^{-1} \in \CO_F/\mfp$ for the zeros of $I$. Since $I$ is irreducible and divides $A_p$ just once, the action of $\phi_p$ on the $\bar{\alpha}_j$ is transitive and the orbit of each $\bar{\alpha}_j$ is of length $4d$. 

Without loss of generality we may label the zeros of $I$ so that
\begin{equation*}
	\phi(\bar{\alpha}_j) =
	\begin{cases}
		\bar{\alpha}_{j+1} \quad\text{for}\ j=1,2,\dots,2d-1, \\
		\bar{\alpha}_1^{-1} \quad \text{for} j =2d.
	\end{cases} 
\end{equation*}
Hence, if we identify $\CS_{2m}$ with the group of permutations of the symbols $\pm1,\pm2,\dots,\pm m$, with the convention that $\pm j$ corresponds to $\bar{\alpha}_{\pm j}$ for $j=1,\dots,2d$, then we find that $\phi$ contains the cycle $(1\,2\, \ldots 2d \, -\!1\, -\!2\, \ldots -\!2d)$ in its disjoint cycle decomposition. Writing $\phi = ((\epsilon_i)_i, \sigma) \in \CC_2 \wr \CS_m$, this implies that $\sigma$ contains the cycle $\tau = (1\, 2 \, \ldots \, 2d)$ in its disjoint cycle decomposition, and that
\begin{equation}
	\label{eq:cycle-prod-epsilon}
	\prod_{i \in \tau} \epsilon_i = 1^{2d-1} \cdot (-1) = -1.
\end{equation}
On the other hand, if $\CG_A$ is $\CC_2^m$-conjugate to a subgroup of $G_5$, we may take an element $((b_i)_i, \id) \in \CC_2 \wr \CS_m$ such that $((b_i)_i, \id)\cdot \phi \cdot ((b_i)_i, \id)^{-1} = ((b_i \epsilon_i b_{\sigma^{-1}(i)})_i, \sigma)$ lies in $\CC_2 \times \CS_m$. Thus there exists $\delta \in \{\pm 1\}$ such that $b_i \epsilon_i b_{\sigma^{-1}(i)} = \delta$ for all $i$. Since the cycle $\tau$ is of even length $2d$, we obtain
\begin{equation}
\label{eq:cycle-prod-conjugate}
1 = \delta^{2d} = \prod_{i \in \tau} b_i \epsilon_i b_{\sigma^{-1}(i)} = \Big( \prod_{i \in \tau} \epsilon_i \Big) \Big( \prod_{i \in \tau} b_i b_{\sigma^{-1}(i)} \Big).
\end{equation}
Since $\sigma$ permutes the elements within the cycle $\tau$, we have
\begin{equation*}
\prod_{i \in \tau} b_i b_{\sigma^{-1}(i)} = \Big( \prod_{i \in \tau} b_i \Big) \Big( \prod_{i \in \tau} b_{\sigma^{-1}(i)} \Big) = \Big( \prod_{i \in \tau} b_i \Big)^2 = 1.
\end{equation*}
Combining this with \eqref{eq:cycle-prod-conjugate} we obtain
\begin{equation*}
\prod_{i \in \tau} \epsilon_i  = 1
\end{equation*}
which contradicts \eqref{eq:cycle-prod-epsilon}. This concludes the proof.
\end{proof}

Our next result about Galois groups, Lemma \ref{lem:Galois_C2xSm_2},	necessitates the following auxiliary lemma.

\begin{lemma}[Prime Polynomial Theorem for reciprocal polynomials]
	\label{lem:PPT for reciprocals}
	Suppose $m \in \Z_{\ge 1}$ is an integer and $p$ is a prime and let $S_p(2m)$ denote the number of monic irreducible reciprocal polynomials of degree $2m$ over $\F_p$. Then, denoting by $\mu$ the M\"obius function, we have
	\begin{equation}
		\label{eq:reciprocal PNT}
		S_p(2m) = 
		\begin{cases}
			\frac{1}{2m} (p^m-1) &\text{if } p>2 \text{ and } m = 2^s \text{ for some } s, \\
			\frac{1}{2m} \sum_{d \mid m,\,  2\nmid d} \mu(d) p^{m/d} &\text{otherwise.}
		\end{cases}
	\end{equation}
	In particular, 
	\begin{equation}
		\label{eq:reciprocal PNT simple}
    S_p(2m) > \frac{p^m}{2m} - \frac{p^{m/3}}{m}.
	\end{equation}
\end{lemma}
\begin{proof}
The expression for $S_p(2m)$ in \eqref{eq:reciprocal PNT} comes from \cite[Theorem~3]{MG}. Combining that with the trivial lower bound $\mu(d) \ge -1$ for $d \ge 3$ yields \eqref{eq:reciprocal PNT simple}.
\end{proof}

\begin{lemma}
\label{lem:Galois_C2xSm_2} Let $m\in\Z_{\ge 3}$ and $p$ be a prime. 
Suppose $\Delta^{\SR}_p(m;m/2) \le m^{-9}$. Then, for every fixed $\epsilon > 0$, we have
\begin{equation*}
\mb{P}_{A \in \CR(m)}\Big(A \textup{ is irreducible and } \CG_A \textup{ is conjugate to a subgroup of } G_5\Big) \ll_{\epsilon} m^{-1/4+\epsilon}.
\end{equation*}
\end{lemma}
\begin{proof}
Define the set
\begin{equation*}
\CI = \bigcup_{1\le d \le m/(2 \log{m})} \big\{I \in \CR_p(2d) : I\ \text{irreducible over}\ \F_p\big\} .
\end{equation*}
Observe that $\CI$ is the set of candidate polynomials $I$ in the statement of Lemma~\ref{lem:Galois_C2xSm_1} in a restricted degree range. Hence it suffices to prove that for any $\epsilon>0$, there exists a constant $C_{\epsilon}>0$ such that
\begin{equation}
\label{eq:I not I^2}
\mb{P}_{A \in \CR(m)}\Big(\exists I \in \CI : I \mid A_p, \, I^2 \nmid A_p\Big) \ge 1-C_\epsilon m^{-1/4+\epsilon}.
\end{equation}

Since each $I\in\CI$ is irreducible and the map $(-)_{\SR}$ is injective, the polynomials $I_{\SR}$ with $I\in\CI$ are irreducible and distinct. Now, consider
\begin{equation*}
f \colon \F_p[T] \setminus \{0\} \to \Z_{\ge 0}, \qquad f(G) = \sum_{\substack{I \in \CI \\ I \mid G^{\SR}, \, I^2 \nmid G^{\SR}}} 1.
\end{equation*}
The function $f$ is \emph{additive}: that is, $f(G_1G_2) = f(G_1) + f(G_2)$ whenever $G_1$ and $G_2$ are coprime (because we must then also have that $G_1^{\SR}$ and $G_2^{\SR}$ are coprime --- see Lemma \ref{lem:coprimality in terms of trace polynomials}). With this notation, \eqref{eq:I not I^2} is equivalent to showing that 
\begin{equation}
	\label{eq:I not I^2 rephrased}
	\mb{P}_{\SR,B\in \CM(m)}\big( f(B_p)=0 \big) \le C_\epsilon m^{-1/4+\epsilon}. 
\end{equation}

Let 
\[
\ell=\big\lfloor m/\log m\big\rfloor 
\]
and, for $G \in \CM(m)$, write 
\begin{equation*}
G_p^{S(\ell)} = \prod_{\substack{J^v \| G_p,\ J \text{ irreducible} \\ \deg{J} \leq \ell, \, J \neq T }} J^v
\end{equation*}
for the \emph{$\ell$-smooth part} of $G_p$ (see \cite[Equation (9.2)]{BKK}). We have 
\begin{equation}
	\label{eq:I not I^2 reduced to smooth part}
	\mb{P}_{\SR,B\in \CM(m)}\big( f(B_p)=0 \big) \le \mb{P}_{\SR,B\in \CM(m)}\Big( f\big(B_p^{S(\ell)}\big)=0 \Big) .
\end{equation}
Having set up notation, observe that $f(J^v) \in \{0,1\}$ when $J\in \F_p[T]$ is monic and irreducible and $v \in \Z_{\ge 1}$. So we may apply \cite[Lemma~9.2(a)]{BKK} with parameters $\theta = 1/2$, $C_1 = 3$, $t = 0$, and with their $m$ being our $\ell$. In addition, note that $\Delta_{\SR,p}(m; m/2) \le \Delta^{\SR}_p(m;m/2) \leq m^{-9}$ by Lemma~\ref{lem:Delta correspondence}. In conclusion, 
\begin{equation*}
\mb{P}_{\SR,B\in \CM(m)} \big(f(B_p)=0 \big) \le \mb{P}_{\SR,B\in \CM(m)} \Big(f\big(B_p^{S(\ell)}\big) = 0\Big) \ll  e^{-L} + 1/m,
\end{equation*}
where 
\begin{equation*}
	L \coloneqq \sum_{\substack{\deg{J} \leq \ell \\ f(J) = 1 \\ J\textup{ irreducible }}} \frac{1}{\norm{J}_p} = \sum_{I \in \CI} \frac{1}{\norm{I}_p^{1/2}}. 
\end{equation*}
By the Prime Polynomial Theorem for reciprocals, Lemma~\ref{lem:PPT for reciprocals}, we have
\begin{equation*}
L = \sum_{1 \le d \le \ell/2} \frac{S_p(4d)}{p^{2d}} \ge \sum_{1\le d \le \ell/2} \frac{1-2p^{-4d/3}}{4d} = \frac{\log \ell}{4}+O(1) .
\end{equation*}
Since $\log{\ell} = \log{m} + O(\log{\log{m}})$, the result follows.
\end{proof}

Lemma \ref{lem:Galois_C2xSm_2} clearly establishes the required estimate \eqref{eq: prop Galois step 3} under the assumption of condition (3). This concludes the proof of Proposition~\ref{prop:Galois}. 

\begin{remark}
Condition (2) of Lemma~\ref{lem:Galois gp of trace pol} and condition (3) of Proposition~\ref{prop:Galois} may be replaced by something weaker, just as condition (2) of Lemma~\ref{lem:replacement_2.4} replaced condition \eqref{eq:mild anti-concentration mod p}. The adjusted condition would then bound the probability that $T^2+1$ divides $A_p$ to a large power. In Proposition~\ref{prop:Galois}, we would still require an anti-concentration condition on all measures $\mu_j$ over the integers, so that we may still apply Proposition~\ref{prop:square_discriminant}.
\end{remark}

\begin{remark}
For squarefree $A$, we have $\CG_A \leq G_2$ if and only if $(-1)^m A(1)A(-1)\Delta(A_{\SR})$ is a nonzero square \cite[Lemma~3.8(b)]{ABD} and $\CG_A \leq G_3$ if and only if $\Delta(A_{\SR})$ is a nonzero square.
\end{remark}

\bibliographystyle{amsplain}

\end{document}